%% file: TAMS_final_version.tex
\documentclass[a4paper, notitlepage, 10pt]{amsart}

\usepackage{graphicx}
\usepackage{geometry}
\usepackage{appendix}
\usepackage{color}
\usepackage{tikz-cd}

\definecolor{refblue}{RGB}{0, 0, 153}
\definecolor{citegreen}{RGB}{0, 115, 0}
\definecolor{linkred}{RGB}{191, 26, 61}

\usepackage{leftidx}
\usepackage[all]{xy}
\usepackage[colorlinks,linkcolor=refblue,citecolor=citegreen,urlcolor=linkred,bookmarks=false,hypertexnames=true]{hyperref} 
\usepackage[nameinlink]{cleveref}
\usepackage{amssymb}
\usepackage{setspace}
\usepackage{stmaryrd}
\usepackage{amsthm}
\usepackage{colonequals}

\newcommand{\stackcite}[1]{\cite[\href{https://stacks.math.columbia.edu/tag/#1}{#1}]{stk}}

\theoremstyle{plain} % changed the style for theorems, lemmas, etc.
\newtheorem{theorem}{Theorem}[section]
\newtheorem{lemma}[theorem]{Lemma}
\newtheorem{proposition}[theorem]{Proposition}
\newtheorem{corollary}[theorem]{Corollary}

\newtheorem{maintheorem}{Theorem}

\newtheorem*{untheorem}{Theorem}
\newtheorem*{undefinition}{Definition}

\theoremstyle{definition} % changed the style for definitions, examples, etc.
\newtheorem{definition}[theorem]{Definition}
\newtheorem{example}[theorem]{Example}

\newtheorem{notation}[theorem]{Notation}

\theoremstyle{remark} % changed the style for remarks
\newtheorem{remark}[theorem]{Remark}

\numberwithin{equation}{section}

\usepackage[skip=5pt, indent=10pt]{parskip}

\usepackage{tikz} % TikZ package for diagrams and drawings
\usetikzlibrary{shapes.geometric, arrows} % TikZ libraries loaded

\input{operators.tex}

\usepackage{amsfonts}
\usepackage{mathrsfs}
\usepackage{microtype}
\usepackage{enumitem}
\usepackage{xspace}

\usepackage{fancyhdr}
\author[M.~Moakher]{Mohamed Moakher} 
 \address{Mohamed Moakher, Department of Mathematics, University of Pittsburgh, Pittsburgh, USA}
 \email{mom224@pitt.edu}

\begin{document}

\title{Quantitative Level lowering for weight two Hilbert modular forms}



%    General info
\makeatletter
\@namedef{subjclassname@2020}{\textup{2020} Mathematics Subject Classification}
\makeatother


\setcounter{tocdepth}{1}

\begin{abstract}
    We generalize a result of Ribet and Takahashi on the parametrization of elliptic curves by Shimura curves to the Hilbert modular setting. In particular, we study the behaviour of the parametrization of modular abelian varieties by Shimura curves associated to quaternion algebras $D$ over a totally real field $F$, as we vary $D$. As a consequence, we obtain that on these Shimura curves, the cohomological congruence module is equal to the ring theoretic congruence module even in cases where we do not have multiplicity one, thereby extending results of Manning and Böckle-Khare-Manning. 
\end{abstract}

\maketitle

\tableofcontents
\section{Introduction}
In \cite{RT97}, Ribet and Takahashi established a formula comparing the degrees of modular parametrizations by Shimura curves of the isogeny class of an elliptic curve. More precisely, let $p$ be a prime number, $N$ be a square-free integer, and let $f$ be a weight two newform of level $\Gamma_0(N)$ with $\bbQ$-coefficients. We can associate to $f$ a Galois representation $\rho_f\colon \Gal(\overline{\bbQ}/\bbQ)\to \GL_2(\bbZ_p)$ which corresponds to an isogeny class of elliptic curves $\mathcal{E}$. Consider a  decomposition $N=N^+Q$ such that $Q$ has an even number of prime factors. Then there exists (up to isomorphism) a unique indefinite quaternion algebra $D_Q$ over $\bbQ$ with discriminant $Q$, which gives rise to a Shimura curve $X_0^{Q}(N^+)$. The Eichler--Shimura construction gives us an elliptic curve $E_{Q,f}\in \mathcal{E}$, and a modular parametrization
$$ X_0^Q(N^+)\longrightarrow E_{Q,f}.$$
The degree of this covering map is referred to as the Shimura degree of $E_{Q,f}$ by $X_0^Q(N^+)$ and is denoted by $\delta_Q$. Now suppose $Q=\overline{Q}q_1q_2$ for distinct prime numbers $q_1$ and $q_2$. Then there is an indefinite quaternion algebra $D_{\overline{Q}}$ of discriminant $\overline{Q}$, and applying the same procedure to $\overline{Q}$ produces another integer $\delta_{\overline{Q}}$. The result of Ribet--Takahashi can then be stated as follows.
\begin{untheorem}{\cite[Theorem 1]{RT97}}
    Assume that $\overline{\rho}_f=(\rho_f \bmod p)^{\text{ss}}$ is absolutely irreducible and that $N^+$ is not a prime, then we have the equality
    \begin{equation}\label{RTresult}
\ord_p(\delta_{\overline{Q}})=\ord_p(\delta_{Q})+c_{q_1}+c_{q_2}, 
    \end{equation}
    where $c_{q_i}$ is the highest power of $p$ such that $\rho_f|_{G_{q_i}}\bmod p^{c_{q_i}}$ is unramified. 
\end{untheorem}
Subsequently in \cite{Takahashi01}, Takahashi removed the condition on $N^+$.

The Shimura degree $\delta_D$ is known to be closely related to the congruence number $\eta_f(N^+,Q)$ which measures the congruences between $f$ and other forms of weight two and level $N$ that are $Q$-new. For instance, the order at $p$ of both these numbers agree under certain multiplicity one results (see \cite{ASR}).
This suggests the following conjecture \cite[Conjecture 1.4]{PW11} which can purely be stated in terms of $f$: 
\begin{equation*}\tag{\textbf{PW}}\label{PW}
\begin{split}
     \text{If }Q=\widetilde{Q}q& \text{ with }q \text{ a prime number,} \text{ and if }\overline{\rho}_f \text{ is absolutely irreducible,} \text{ then } \\ &\ord_p(\eta_f(N^+q,\widetilde{Q}))=c_{q} + \ord_p(\eta_f(N^+,Q)). 
\end{split}
\end{equation*}
Pollack and Weston proved this under certain conditions (see \cite[Theorem 6.11]{PW11}), notably requiring $\overline{\rho}_f $ to be ramified at $q$ whenever $q\mid Q$ and $q\equiv \pm 1 \bmod p$. They coined the term \emph{quantitative level lowering} to describe this type of result. Under favorable assumptions ensuring multiplicity one results and that certain Hecke algebras are complete intersection, Kim and Ota (\cite{KO19}) prove this conjecture for modular forms of higher weight.

In \cite{BKM2}, Böckle, Khare, and Manning offered a new and more solid perspective on the result of Ribet and Takahashi. They gave a new proof of one of the inequalities in \eqref{RTresult} using only the Taylor--Wiles--Kisin patching method. They also obtained a stronger result (\cite[Theorem 7.4]{BKM2}) towards $(\textbf{PW})$ in cases where the relevant Hecke algebras are not complete intersection (in fact, the failure of being complete intersection occurs whenever the restriction of $\overline{\rho}_f$ at the decomposition group of $q$ is trivial for some $q\mid Q$ with $q\equiv 1 \bmod p$).

In this article, we generalize the results of Ribet--Takahashi and Böckle--Khare--Manning, concerning classical modular forms, to the weight two Hilbert modular case. Now let $F$ be a totally real field with $[F:\mathbb{Q}]=d$, and fix an odd prime number $p$ which does not ramify in $F$. Let $\mathbf{f}$ be a Hilbert cuspidal newform of parallel weight $2$, level $\mathfrak{n}$ prime to $p$, and nebentypus $\psi$. We fix a finite extension $E$ of $\bbQ_p$ with ring of integers $\mathcal{O}$, uniformizer $\varpi$, and residue field $k$. We assume that $\mathcal{O}$ is large enough so that the Galois representation associated to $\mathbf{f}$ is of the form
$$ \rho_{\mathbf{f}}\colon G_F \longrightarrow \GL_2(\mathcal{O}). $$
Let $Q$ be a finite set of places dividing $\mathfrak{n}$ and such that $\mathbf{f}$ is Steinberg at these places. Then we can associate to $Q$ a quaternion algebra $D_Q$, and a space of automorphic forms $M_Q$ endowed with an action of a localized Hecke algebra $\bbT^{Q}$\footnote{We caution the reader that in this introduction, we use a simplified notation. For a precise definition, see \Cref{setting}, where these are denoted by $M_Q(U^{\mathbf{f}})$ and $\overline{\bbT}{}^Q_{\mathcal{O}}(U^{\mathbf{f}})_{\mathfrak{m}}$.}. Let us also make the following definition.
\begin{undefinition}
    Let $\rho\colon G_F\to \GL_2(\mathcal{O})$ be a continuous Galois representation with residual Galois representation $\overline{\rho}\colonequals(\rho\bmod \varpi)^{\text{ss}}$,  and let $v$ be a finite place of $F$.
    \begin{enumerate}
        \item[1.] We write $c_v(\rho)\in \bbN$ (or simply $c_v$ if $\rho$ is clear from the context) for the highest power $n$ of $\varpi$ such that $\rho|_{G_v}\bmod \varpi^n$ is unramified.
        \item[2.] We say that $v$ is a trivial place for $\rho$ if $\overline{\rho}|_{G_v} $ is unramified, and $\overline{\rho}(\Frob_v)$ is a scalar.
    \end{enumerate}
\end{undefinition}

At the level of generality we consider, the algebra $\bbT^{Q}$ will not be Gorenstein, and so the module $M_Q$ will fail to be free over $\bbT^{Q}$. This occurs whenever $Q$ contains a \emph{trivial} place. 
In fact we prove the following theorem:

\begin{maintheorem}\label{maintheorem0}[\Cref{freenessoverGorenstein,notGorenstein}]
    Assume that $p>2$ and that the residual Galois representation $\overline{\rho}:=(\rho_{\mathbf{f}} \bmod \varpi)^{\text{ss}}$ satisfies the following hypotheses:
\begin{enumerate}
    \item[$(a)$] The image $\overline{\rho}$ is non-exceptional:  $\overline{\rho}(G_F)$ contains a subgroup of $\GL_2(\overline{\bbF_p})$ conjugate to $\SL_2(\bbF_p)$.
    Moreover if $p=5$ and the projective image of $\overline{\rho}$ is isomorphic to $\PGL_2(\bbF_5)$, then the kernel of $\proj\overline{\rho}$ does not fix $F(\zeta_5)$.
    \item[$(b)$] If $v $ is a finite place of $F$ such that $\overline{\rho}|_{G_v}$ is ramified and $q_v\equiv -1 \bmod p$, then either $\overline{\rho}|_{I_v}$ is irreducible, or $\overline{\rho}|_{G_v}$ is absolutely reducible (we exclude vexing primes, see \cite{Dia97}). 
    \end{enumerate} 
        Then the following statements are equivalent:
    \begin{enumerate}
        \item[1.] $\bbT^Q$ is Gorenstein,
        \item[2.]  $M_Q$ is free over $\bbT^Q$,
        \item[3.] $Q$ does not contain any trivial place for $\overline{{\rho}}$.
    \end{enumerate}
\end{maintheorem}
\begin{remark}
    Note that hypothesis $(a)$ implies the usual Taylor--Wiles condition that the restriction $\overline{\rho}|_{G_{F(\zeta_p)}}$ is absolutely irreducible (see \cite[(3.2.3)]{Kis09}). We make this more restrictive hypothesis in order to use Ihara's lemma for Shimura curves over $F$, as established in \cite{MS21}.
\end{remark}

Therefore, when $Q$ contains trivial places, we distinguish between two congruence modules: the \emph{cohomological} congruence module $\Psi_{\lambda_{\mathbf{f}}}(M_Q)$, associated to the space of automorphic forms,  and the \emph{ring-theoretic} congruence module $\Psi_{\lambda_{\mathbf{f}}}(\bbT^Q)$, intrinsic to the Hecke algebra. In practice, the cohomological congruence module can be related to special values of automorphic $L$-functions, while the ring-theoretic one encodes congruences between $\mathbf{f}$ and other forms. Thus, it is interesting to ask whether these
two cyclic modules are isomorphic. Although one always has $\ell_{\mathcal{O}}(\Psi_{\lambda_{\mathbf{f}}}(\bbT^Q))\ge \ell_{\mathcal{O}}(\Psi_{\lambda_{\mathbf{f}}}(M_Q))$, there is no apparent reason for equality to hold in general.

In some cases, the patching method alone suffices to establish an isomorphism between these two modules. Notably, in the minimal case, Manning showed that $\bbT^Q\cong \End_{\bbT^Q}(M_Q)$ (\cite[Theorem 1.2]{Man21}), which implies the desired equality. For non-optimal levels, however, this equality is more challenging to prove, and this is precisely our contribution. Our approach combines the arithmetic geometry ideas of Ribet and Takahashi with the new techniques of Böckle, Khare, and Manning to prove the following result.

\begin{maintheorem}[\Cref{changeofcongmoduleforheckealg,finalresultoverindefinite}]\label{maintheorem}
    Assume that the residual Galois representation $\overline{\rho}=(\rho_{\mathbf{f}} \bmod \varpi)^{\text{ss}}$ satisfies  the same assumptions $(1)$ and $(2)$ of \Cref{maintheorem0}.
    
Let $Q_0=\{t\}$ if $[F:\bbQ]$ is even, and $Q_0=\{t,s\}$ if $[F:\bbQ]$ is odd, where $t$ and $s$ are finite places of $F$ that are not trivial for $\overline{\rho}_{\mathbf{f}}$ and at which $\mathbf{f}$ is Steinberg. 
If $Q$ is a finite set of finite places at which $\mathbf{f}$ is Steinberg and $Q_0\subseteq Q$, then 
$$\Psi_{\lambda_{\mathbf{f}}}(\bbT^Q)\cong \Psi_{\lambda_{\mathbf{f}}}(M_Q),$$ and if $Q=\{v\}\sqcup\overline{Q}$ with $Q_0\subseteq \overline{Q}$, then
$$ \ell_{\mathcal{O}}(\Psi_{\lambda_{\mathbf{f}}}(M_{\overline{Q}}))=\ell_{\mathcal{O}}(\Psi_{\lambda_{\mathbf{f}}}(M_Q))+c_v. $$
\end{maintheorem}
Our theorem can be seen as a direction towards the analogue of conjecture $\textbf{(PW)}$ in the Hilbert modular setting. Even in the case $F=\bbQ$, it yields new results, notably in the case where $K_{\mathbf{f}}\neq \bbQ$. Nevertheless, we note that the assumption on the existence of the set $Q_0$ is \emph{essential} to our method, and so it excludes the case where all the places in $Q$ are trivial for $\overline{\rho}$.

\begin{example}
    Let $F=\mathbb{Q}(\sqrt{13})$ with $\mathcal{O}_F=\mathbb{Z}[a]$, where $a^2=a+3$. Consider the elliptic curve labeled $2.2.13.1$-$1391.1$-$a1$ in \cite{lmfdb} with global minimal model
    $$ E: y^2+xy+y=x^3+(1+a)x^2 -(5+2a)x + (6+4a),$$
    discriminant $\Delta=-81887-62855a=(17-9a)(14+11a)^3$, and conductor $\mathfrak{n}=(1-2a)(10+a)=\mathfrak{q}_1\mathfrak{q}_2$. The associated mod-$3$ representation $\overline{\rho}_E\colon G_F\to \GL_2(\mathbb{F}_3)$ is surjective, and has tame conductor $\mathfrak{q}_2$ since $v_{\mathfrak{q}_1}(\Delta)=3$. Writing $\mathbf{f}_E$ for the associated Hilbert modular form, we find that it is congruent to a form of level $\mathfrak{q}_2$. In fact, the $\varpi$-order of the Tamagawa number $c_{\mathfrak{q}}$ of $E$ agrees with $ c_{\mathfrak{q}}(\rho_E)$, and $c_{\mathfrak{q}_1}=3$. Hence \Cref{maintheorem} gives
    $$ \ell_{\mathcal{O}}(\Psi_{\lambda_{\mathbf{f}_E}}(M_{\mathfrak{q}_2}))=\ell_{\mathcal{O}}(\Psi_{\lambda_{\mathbf{f}_E}}(M_{\mathfrak{q}_1\mathfrak{q}_2}))+\ord_{\varpi}(3).$$
    In order to find trivial primes for $\overline{\rho}_E$, we follow \cite{Rib90} and \cite[p.305]{Ser72}. Suppose we have a prime number $\ell\equiv 1\bmod 3\cdot 13$ that splits as $(\ell)=\mathfrak{q}\mathfrak{q}'$, and suppose  $\tr(\rho_E(\Frob_{\mathfrak{q}}))\equiv \pm 2 \bmod 3$. We must distinguish the two possibilities
    $$ \overline{\rho}_E(\Frob_{\mathfrak{q}})=\pm \id, \quad \quad \text{ and } \quad \quad \overline{\rho}_E(\Frob_{\mathfrak{q}})=\begin{pmatrix}
        \pm 1 & x \\ 0 & \pm 1
    \end{pmatrix}, x\neq 0,  $$
  which, up to a sign, have orders $1$ and $3$, respectively. To do so, we identify $\GL_2(\mathbb{F}_3)/\{\pm 1\}\cong \mathfrak{S}_4$ by viewing $G_F$ as acting on the $x$-coordinates $\{x_1,x_2,x_3,x_4\}$ of the points of $E$ of order $3$. Using the formula $\Delta^{1/3}=b_4-3(x_1x_2-x_3x_4)$, with $b_4\in F$ computed from the coefficients of the Weierstrass model of $E$, we obtain a composition 
    $$ G_F\xrightarrow{\overline{\rho}_E} \GL_2(\mathbb{F}_3)/{\pm 1} \xrightarrow{\sim} \mathfrak{S}_4 \longrightarrow \mathfrak{S}_3\cong \Gal(F(\Delta^{1/3},\mu_3)/F). $$
    Therefore, $\mathfrak{q}$ is trivial if and only if $\Frob_{\mathfrak{q}}(\Delta^{1/3})=\Delta^{1/3}$ if and only if $\Delta^{(\ell-1)/3}\equiv (17-9a) ^{(\ell-1)/3} \equiv 1 \bmod \mathfrak{q}$. Letting $a_0\in \mathbb{Z}$ be congruent to the image of $a$ in $\mathcal{O}_F/\mathfrak{q}=\mathbb{F}_{\ell}$, this condition is equivalent to $(17-9a_0) ^{(\ell-1)/3}\equiv 1 \bmod \ell$. One such example is $\ell=157$ with $\mathfrak{q}=(157,a-90)$.
\end{example}

\subsection*{Strategy of the proof}
We first start by proving \Cref{maintheorem0}. The implication $3.\Rightarrow 2.$ is done using the patching method by upgrading a freeness result at the minimal level. We use Ihara's lemma proven in our context in \cite{MS21} and \cite{Tay06} to give a quantitative version of level raising as in \cite{Diamond97}. We then use an improved version of Wiles' numerical criterion given in \cite{BIK} to conclude. This is the purpose of \Cref{galdef}. To our knowledge, this case of level raising in this level of generality did not appear elsewhere in the literature. For the equivalence $1. \Leftrightarrow 3.$, we introduce a measure of Gorensteinness defect similar to that of \cite{WWE21} and \cite{KW08}. We prove that it remains invariant when
quotienting by regular sequences and that it is multiplicative under completed tensor products. This allows us to conclude in \Cref{notGorenstein}.

By assumption on $Q_0$ and by \Cref{maintheorem0}, we have that $M_{Q_0}$ is free of rank one over $\bbT^{Q_0}$, which gives us that $\ell_{\mathcal{O}}(\Psi_{\lambda_{\mathbf{f}}}(\bbT^{Q_0}))= \ell_{\mathcal{O}}(\Psi_{\lambda_{\mathbf{f}}}(M_{Q_0}))$. Our idea is to propagate this equality to $Q$ by generalizing the method of Ribet and Takahashi to modular abelian
varieties over totally real fields.  As in their setting, we work with the Jacobian $J_Q$ of the Shimura curve associated to an indefinite quaternion algebra $D_Q$.  We write $T^\vee_{\mathfrak{m}}J_Q$ for the localized contravariant $p$-adic Tate module of $J_Q$, which satisfies $T^\vee_{\mathfrak{m}}J_Q\cong M_Q^2$. We also extract from $J_{Q_0}$ an auxiliary abelian variety with a Hecke action, that we denote by $J^{\min}$, and that satisfies $T^\vee_{\mathfrak{m}}J^{\min}\cong(\bbT^Q)^2$.

Since we do not assume that $K_{\mathbf{f}}=\bbQ$, we can no longer work with Shimura degrees. Indeed, $\mathbf{f}$ might have congruences with forms that are Galois conjugate to it, yet geometrically,
$\mathbf{f}$ and its Galois conjugates are indistinguishable. Therefore, we introduce in \Cref{lambdashimuradeg} a modified version of the Shimura degree, which we call the $\lambda$-Shimura degree $\delta_Q$. 

Similarly to what is done in \cite{Khare03} (except that we allow $K_{\mathbf{f}}$ to be ramified at $p$), we relate the $\lambda$-Shimura degree to the size of the cohomological congruence module by working with the localized character group of the toric part of $J_{Q}$ at a place $t\in Q_0$, which we denote by $\mathcal{X}_t(J_{Q})$. A key result in this setting is \Cref{tatemoduletochargrp}, which is proved following an argument of \cite{Helm}. Conditionally on the freeness of $\mathcal{X}_t(J^{\min})$ over $\bbT^Q$, this theorem allows us to show that
$$ \Psi_{\lambda_{\mathbf{f}}}(M_Q)\cong \Psi_{\lambda_{\mathbf{f}}}(\mathcal{X}_t(J_Q)). $$
Moreover, by a study of the $p$-adic uniformization of abelian varieties in \Cref{factsonAB}, we
also deduce in \Cref{shdegtocongmodule} that 
$$\delta_Q=\ell_{\mathcal{O}}\big(\Psi_{\lambda_{\mathbf{f}}}(\mathcal{X}_t(J_Q))\big)=\ell_{\mathcal{O}}(\Psi_{\lambda_{\mathbf{f}}}(M_Q)).$$
 
The main ingredient in the proof
is Ribet’s exact sequence realizing the Jacquet-Langlands correspondence on the
character groups of the toric parts of Jacobians, thereby allowing us to compare $\lambda$-Shimura degrees. This is explained in \Cref{proofofmaintheorem}, where we show using a delicate level-raising argument that $\mathcal{X}_t(J^{\min})$ is free of rank one over $\bbT^Q$ (see \Cref{freeToricModule}), and where we ultimately obtain the inequality
$$\ell_{\mathcal{O}}(\Psi_{\lambda_{\mathbf{f}}}(M_{\overline{Q}}))\ge\ell_{\mathcal{O}}(\Psi_{\lambda_{\mathbf{f}}}(M_Q))+c_{v}+c_{w},$$
with $Q=\overline{Q}\sqcup\{v,w\}$. 
Working purely on the ring-theoretic side, the method of \cite{BKM2} gives us that
$$\ell_{\mathcal{O}}(\Psi_{\lambda_{\mathbf{f}}}(\bbT^{\overline{Q}}))=\ell_{\mathcal{O}}(\Psi_{\lambda_{\mathbf{f}}}(\bbT^Q))+c_{v}+c_{w},$$
from which we deduce the main theorem by induction. Finally, we note that in order to get the result for $\widetilde{Q}= Q\setminus\{v\}$ corresponding to a definite quaternion algebra, we prove
in \Cref{chargroupvnotinQ} that $M_{\widetilde{Q}}\cong \mathcal{X}_v(J_{\overline{Q}})$, analogously to statements found in the works of Kohel \cite{Kohel} and Takahashi \cite{Takahashi}.

\subsection*{Acknowledgements} This work forms part of the author's PhD thesis, supervised by Stefano Morra, whom I thank for his guidance and encouragement throughout the project. The project itself grew out of an interest in a series of papers by Chandrashekhar Khare and his collaborators, and I am grateful to C. Khare for the stimulating conversations we had during his visit to Paris 13 in the summer of 2022, and for his encouragement to pursue these questions further. A substantial part of this research was carried out at Oxford University in the spring of 2023, at the invitation of James Newton, to whom I am grateful both for this opportunity and for his guidance during my stay. I also thank Gebhard Böckle and Carl Wang-Erickson for their careful reading of earlier versions of this work, and Matteo Tamiozzo and Aleksander Horawa for valuable conversations. Finally, I thank the anonymous referee for numerous suggestions that greatly improved the quality of the exposition.

\subsection*{Notation}
We will use the following notation:
\begin{enumerate}
\item  $F$ is a totally real field of degree $d$ over $\mathbb{Q}$, with ring of integer $\mathcal{O}_F$. We fix an algebraic closure $\overline{F}$, and write $G_F\colonequals \Gal(\overline{F}/F)$ for the absolute Galois group.
\item $p$ is an odd prime number unramified in $F$.
\item $E$ is a finite extension $\bbQ_p$ with ring of integers $\mathcal{O}$, uniformizer $\varpi$, and residue field $k$.
    \item If $v$ is a finite place of $F$, we write $F_v$ for the completion of $F$ with respect to $v$, $k_v$ for its residue field, and $\varpi_v$ for a fixed choice of a uniformizer. We write $q_v$ for the cardinality of the residue field of $F_v$.
    \item If $v$ is a finite place of $F$, we write $G_v\colonequals \Gal(\overline{F_v}/F_v)$ for the absolute Galois group of $F_v$, and $I_v\subset G_v$ for its inertia subgroup. We fix a lift $\sigma_v\in I_v$ of a generator of the $\bbZ_p$-quotient of $I_v$, and a lift $\Frob_v\in G_v$ of the arithmetic Frobenius.
     We also fix an embedding $\overline{F}\hookrightarrow \overline{F_v}$, inducing an embedding $G_v\hookrightarrow G_F$.
    \item We set  $\widehat{\mathcal{O}}_F=\mathcal{O}_F\otimes \widehat{\mathbb{Z}}$, $\mathbb{A}_{F}^\infty=\widehat{\mathcal{O}}_F\otimes \mathbb{Q}$, $\mathbb{A}_{F,\infty}= F\otimes_{\mathbb{Q}} \mathbb{R}$, and let $\mathbb{A}_F=\mathbb{A}_{F}^\infty\times \mathbb{A}_{F,\infty}$ be the ring of adeles of $F$.
    \item  We let $J_F$ be the set of embeddings of $F$ into $\mathbb{R}$, and we fix an isomorphism $\iota\colon \mathbb{C}\xrightarrow{\sim}\overline{\mathbb{Q}}_p$. Using $\iota$ we have a decomposition:
$$ J_F=\bigsqcup_{v|p} J_v $$
where an element $\tau\in J_F$ lies in $J_v$ if and only if the composition $\iota\circ \tau$ induces the $p$-adic place $v$ on $F$.
\item We write $\mathbb{A}_{F,p}=F\otimes_{\mathbb{Q}}\mathbb{Q}_p\cong \prod_{v|p}F_v$. If $\tau\in J_v$, then $\tau$ extends to a $\mathbb{Q}_p$-linear embedding $\tau\colon  F_v\hookrightarrow \overline{\mathbb{Q}}_p$ for which we use the same symbol. In this way, we identify $J_v$ with $\Hom_{\mathbb{Q}_p}(F_v,\overline{\mathbb{Q}}_p)$.
\item $\epsilon_p\colon G_F\rightarrow \mathbb{Z}_p^\times$ is the $p$-adic cyclotomic character. 

\item If $M$ is a finite torsion $\mathcal{O}$-module, we let $\ell_{\mathcal{O}}(M)$ denote its length. And if $x\in \mathcal{O}$, we let $\ord_{\varpi}(x)\colonequals\ell_{\mathcal{O}}(\mathcal{O}/(x))$.
\item $\CNL_\mathcal{O}$ is the category of complete noetherian local $\mathcal{O}$-algebras with residue field $k$. 
\end{enumerate}
 
\section{Shimura curves and spaces of quaternionic forms}\label{Section2}
\subsection{Setting}\label{setting}
Let $\mathbf{f}$ be the Hilbert cuspidal newform in the introduction, which is of parallel weight $2$, level $\mathfrak{n}$ prime to $p$, and nebentypus $\psi\colon \mathbb{A}_F^\times/F^\times \to \mathbb{C}^\times$ (in the sense of \cite[\S 12.3]{Nek06}).  We write $\psi_{\iota}\colonequals \iota\circ \psi$, and assume that $\psi_{\iota}$ takes values in $\mathcal{O}$.

Let $S$ be a finite set of places of $F$. For a commutative ring $A$, we write $\mathbb{T}^{S}_A$ for the polynomial algebra over $A$ in the indeterminates $T_v, S_v$ for all finite places $v\not\in S$ of $F$. Given a subset $S'\subseteq S$, we also write $\overline{\bbT}{}^{S,S'}_A$ for the $\bbT^{S}_A$-polynomial algebra $\mathbb{T}^{S}_{A}\left[\mathbf{U}_v,v\in S'\right]$. 
We will omit the subscript $-_A$ if $A=\bbZ$. 

We assume that $\mathcal{O}$ and $S$ are sufficiently large (and $S'$ is chosen appropriately) so that the Hilbert modular form $\mathbf{f}$ induces an augmentation of the universal Hecke algebra
$$ \lambda_{\mathbf{f}}\colon\overline{\mathbb{T}}{}^{S,S'}_{\mathcal{O}} \longrightarrow \mathcal{O}. $$
We will denote by the same symbol the corresponding augmentation on its subsequent quotients.

The $p$-adic Galois representation associated to $\mathbf{f}$
$$ \rho_{\mathbf{f}}: G_F \longrightarrow \GL_2(\mathcal{O}) $$
satisfies the following properties (see \cite{Car}, \cite[Lemma 3.4.2]{Kis09}, \cite[\S 12.4]{Nek06}):
\begin{itemize}
    \item $\rho_{\mathbf{f}}$ is unramified outside $\mathfrak{n}p$.
    \item For a finite place $v \nmid p\mathfrak{n}$, the characteristic polynomial of $\rho_{\mathbf{f}}(\Frob_v)$ is
    $$\mathrm{char}\big(\rho_{\mathbf{f}}(\Frob_v),X\big)=X^2-\lambda_{\mathbf{f}}(T_v)X+ \psi_\iota(\varpi_v)q_v.$$
    \item For $v\mid p$, $\rho_{\mathbf{f}}|_{G_{v}}$ is crystalline with $\tau$-Hodge-Tate weights $\{0,1\}$ for $\tau\in J_v$.
\end{itemize}
We assume that the residual Galois representation $\overline{\rho}=(\rho_{\mathbf{f}} \bmod \varpi)^{\text{ss}}$ satisfies the assumptions of \Cref{maintheorem0}, and that $k$ contains the eigenvalues of all the elements in the image of $\overline{\rho}$. We associate to $\overline{\rho}$ the maximal ideal  $\mathfrak{m}\subset \mathbb{T}^{S}_{\mathcal{O}}$ given by $\mathfrak{m}=\lambda_{\mathbf{f}}^{-1}(\varpi\mathcal{O})$. We can then see $\overline{\rho}$ as a representation $\overline{\rho}\colon G_F\to \GL_2(\bbT^{S}_{\mathcal{O}}/\mathfrak{m})$.
\subsection{Spaces of quaternionic forms}

 We write $\tau_1,\dots,\tau_d$ for the distinct real embeddings of $F$. For a finite set $Q$ of finite places of $F$, let $D_Q$ be the quaternion algebra over $F$ ramified exactly at $Q\cup \{\tau_2,\dots,\tau_d\}$ if $|Q| \equiv d-1 \bmod 2$, and at $Q\cup \{\tau_1,\dots,\tau_d\}$ otherwise. We say that $Q$ is of \emph{indefinite type} in the first case, and of \emph{definite type} in the second case.

 We fix a maximal order $\mathcal{O}_{D_Q}\subseteq D_Q$ and an isomorphism $$\mathcal{O}_{D_Q}\otimes_{\mathcal{O}_F}\prod_{v\nmid Q\infty} \mathcal{O}_{F_v}\cong \prod_{v\nmid Q\infty} M_2(\mathcal{O}_{F_v}). $$  We write $G_{Q}$ for the reductive group over $\mathcal{O}_F$ whose functor of points is $G_{Q}(R)=(D_Q\otimes_{\mathcal{O}_F}R)^\times$. Thus, we obtain for each finite place $v\not\in Q$ of $F$ an isomorphism $G_{Q}(\mathcal{O}_{F_v})\cong \GL_2(\mathcal{O}_{F_v})$.
 
 Let $v$ be a finite place of $F$. We define for each $n\ge 1$ the open compact subgroups
\begin{equation*}%\label{subgroup}
    \Gamma_1^1(v^n)\subset \Gamma_1(v^n)\subset \Gamma_0(v^n) \subset \GL_2(\mathcal{O}_{F_v}),
\end{equation*}
with
\begin{itemize} \vspace{-6pt}
\item  $\ \Gamma_0(v^n)= \left\{\begin{pmatrix} *& * \\ 0 & * \end{pmatrix} \bmod \varpi^n_v \right\}, $ \quad \quad \quad $\bullet \ \Gamma_1(v^n)= \left\{\begin{pmatrix} * & * \\ 0 & 1 \end{pmatrix} \bmod \varpi^n_v \right\}, $
    \item $ \  \Gamma_1^1(v^n)= \left\{\begin{pmatrix} 1& * \\ 0 & 1 \end{pmatrix} \bmod \varpi^n_v \right\}.$
\end{itemize}
By \cite[Lemma 4.11]{DDT}, we may fix a finite place $a$ of $F$ satisfying 
\begin{itemize}
    \item $q_a>4^d$, and $q_a>|\bbN_{F/\bbQ}(\zeta-1)|$ for every non-trivial root of unity $\zeta$ in $F$.
    \item $q_a\not\equiv 1 \bmod p$.
    \item $\rho_{\mathbf{f}}$ is unramified at $a$, and the ratio of the eigenvalues of $\overline{\rho}(\Frob_a)$ is not equal to $q_a^{\pm 1}$ in $k$.
\end{itemize}
We say that $U\subset G_{Q}(\mathbb{A}_F^\infty)$ is \emph{good} if it satisfies the following conditions
\begin{itemize}
    \item $U=\prod_v U_v$ for open compact subgroups $U_v\subset G_{Q}(F_v)$.
    \item For each $v\in Q$, $U_v$ is the unique maximal compact subgroup of $G_{Q}(F_v)$.
    \item $U_a=\Gamma_1^1(a)$.
\end{itemize}
We will write $\mathcal{J}_Q$ for the set of good subgroups $U\subset G_Q(\mathbb{A}^\infty_F)$. For $U\in \mathcal{J}_Q$, we let $\Sigma(U)$ be the set of finite places $v\not \in Q\cup \Sigma_p\cup \{a\}$ such that $U_v\not\cong \GL_2(\mathcal{O}_{F_v})$. 
\\ If $w\not \in \Sigma(U)$, we write 
$$ U_0(w)= U^w \Gamma_0(w) \subset U= U^w\GL_2(\mathcal{O}_{F_w}). $$
We say that $U\in \mathcal{J}_Q$ is a \emph{congruence subgroup} if for every finite place $v$, $U_v$ is isomorphic to $\Gamma_0(v^n)$ or $\Gamma_1(v^n)$ for some integer $n$. 
\subsubsection{}\label{indefinitesetting} Assume that $Q$ is of \emph{indefinite type}. We fix an isomorphism $D_Q\otimes_{F,\tau_1}\mathbb{R}\cong M_2(\mathbb{R})$, and write $X$ for the $G_Q(F\otimes_{\mathbb{Q}}\mathbb{R})$-conjugacy class of the homomorphism $h: \mathbb{S}= \Res_{\mathbb{C}/\mathbb{R}} \mathbb{G}_m\rightarrow (\Res_{F/\mathbb{Q}} G_Q)_{/\mathbb{R}} $ which sends $z=x+iy\in \mathbb{C}^\times=\mathbb{S}(\mathbb{R})$ to the element
$$ \Big(\left(\begin{smallmatrix}x & y \\ -y & x\end{smallmatrix}\right)^{-1}, 1,\cdots,1\Big)\in \prod_{i=1}^d G_Q(F\otimes_{F,\tau_i}\mathbb{R})\cong \GL_2(\bbR)\times \bbH^\times \times \cdots \times \bbH^\times. $$
For each $U\in \mathcal{J}_Q$, there is an associated compact Riemann surface given by
$$ X_Q(U)(\mathbb{C})= G_Q(F)\backslash \big(G_Q(\mathbb{A}_{F,f})/U \times X \big).$$

By the theory of Shimura varieties, there is a canonical projective algebraic curve $X_Q(U)$ over $F$ such that $X_Q(U)(\mathbb{C})$ is its set of complex points for the embedding $\tau_1: F \rightarrow \mathbb{C}$. We write $[\gamma,x]$ for the point in $X_Q(U)(\bbC)$ corresponding to $\gamma\in G_Q(\mathbb{A}_{F,f})$ and $x\in X$. There is a natural right action of the group $G_Q(\mathbb{A}_F^{a,\infty})$ on the projective system $\{ X_Q(U)(\mathbb{C})\}_{U\in \mathcal{J}_Q}$ given by the system of isomorphisms
\begin{align*}
    [\ \cdot \ g]\colon X_Q(U)(\mathbb{C})&\rightarrow X_Q(g^{-1}Ug)(\mathbb{C})
    \\ [\gamma,x]&\mapsto [\gamma g,x].
\end{align*}
We will write $H_Q(U)=H^1_{\text{ét}}(X_Q(U)_{\overline{F}},\mathcal{O})$ for the first étale cohomology group, which is a finite free $\mathcal{O}$-module.

For $U,U'\in \mathcal{J}_Q$ and $g\in G_Q(\bbA_{F}^{a,\infty})$, we define the Hecke correspondence $[UgU']$ by the diagram
\begin{center}
    \begin{tikzcd}
        X_Q(U\cap gU'g^{-1}) \arrow[d,"pr"] \arrow[r,"{[\ \cdot \ g ]}"] & X_Q(g^{-1}Ug \cap U') \arrow[d,"pr'"]
        \\ X_Q(U) \arrow[r,dashed] & X_Q(U')
    \end{tikzcd}
\end{center}
so that it induces a morphism on étale cohomology
$$ [UgU']\colonequals pr_*\circ [\ \cdot \ g]^*\circ pr'^*\colon H_Q(U') \to H_Q(U). $$
Poincaré duality induces a perfect pairing 
\begin{equation}\label{PoincaréPairing}
    (\cdot,\cdot)_U \colon H_Q(U)\times H_Q(U)\longrightarrow \mathcal{O} 
\end{equation}
such that $([UgU]\cdot,\cdot)_U= (\cdot,[Ug^{-1}U]\cdot)_U$.

Let $S=Q\cup \Sigma(U)\cup \Sigma_p$, and let $S'\subseteq\Sigma(U)$ be a set of places $v$ such that $U_v=\Gamma_0(v)$. The Hecke algebra $\mathbb{T}^{S}_{\mathcal{O}}$ acts on $H_Q(U)$, with each Hecke operator $T_v,S_v$ $(v\not\in S)$ acting respectively by the double coset operator
$$ \Big[ \GL_2(\mathcal{O}_{F_v}) \begin{pmatrix} \varpi_v & 0 \\ 0 & 1 \end{pmatrix} \GL_2(\mathcal{O}_{F_v})\Big]  \quad \text{ and, } \quad \Big[ \GL_2(\mathcal{O}_{F_v}) \begin{pmatrix} \varpi_v & 0 \\ 0 & \varpi_v \end{pmatrix} \GL_2(\mathcal{O}_{F_v}) \Big].$$
We define $\bbT^Q_{\mathcal{O}}(U)$ to be the image of $\mathbb{T}^{S}_{\mathcal{O}}$ in $\End_{\mathcal{O}}(H_Q(U))$. Since $U$ is a good subgroup, the algebra $\mathbb{T}^Q_{\mathcal{O}}(U)$ is reduced and $\mathcal{O}$-torsion free. We also write $\overline{\bbT}{}^Q_{\mathcal{O}}(U)$ for the image of $\overline{\bbT}{}^{S,S'}_{\mathcal{O}}$ in $\End_{\mathcal{O}}(H_Q(U))$, where $\mathbf{U}_v$ acts by the double coset operator
$$\Big[ U_v \begin{pmatrix} \varpi_v & 0 \\ 0 & 1 \end{pmatrix} U_v\Big].  $$
The group $G_Q(\bbA_{F,f})$ admits an involution $(-)^{\alt}$ that is given for each finite place $v$ by
\begin{align*}
    D_{Q,v}^\times &\longrightarrow D_{Q,v}^\times
    \\ d_v &\mapsto d_v^{\alt} \colonequals \frac{1}{\Nrd(d_v)}d_v.
\end{align*}
This induces a morphism of Shimura curves 
\begin{align*}
    (-)^{\alt} \colon X_Q(U)(\bbC)& \longrightarrow X_Q(U^{\alt})(\bbC)
    \\ [\gamma, x] &\mapsto[\gamma^{\alt},x].
\end{align*}
By \cite[Lemma 2.2]{Fouquet}, this morphism is defined over a finite abelian extension of $F$, that we denote by $F_U$. If $U\in \mathcal{J}_Q$ is a congruence subgroup, we define the element $w_U=(w_{U,v})_v\in G_Q(\mathbb{A}_{F,f})$ such that $w_{U,v}=\id$ if $v\not\in \Sigma(U)$, and
$$ w_{U,v}=\begin{pmatrix}
    0 & -1 \\ \varpi_v^{n} & 0
\end{pmatrix} $$
if $U_v\cong \Gamma_0(v^n)$ or $ \Gamma_1(v^n)$. The composition 
$$ W_U\colon  X_Q(U)\xrightarrow{(-)^{\alt}} X_Q(U^{\alt})\xrightarrow{[\ \cdot  w_U]} X_Q(U) $$
defines an $F_U$-morphism that we refer to as the Atkin--Lehner operator. We have the following relations
$$ W_U\circ T_v=S_v^{-1}T_v\circ W_U, \quad \quad W_U\circ S_v=S_v^{-1}\circ W_U, \quad \quad W_U\circ \mathbf{U}_v=\mathbf{U}_v^*\circ W_U, $$
where $\mathbf{U}_v^*$ is the operator acting by $[U_v\mathrm{diag}(\varpi_v^{-1},1)U_v]$.
This allows to introduce a modified Poincaré pairing 
\begin{align*}
    \langle\cdot,\cdot\rangle_U \colon H_Q(U)\times & H_Q(U)\longrightarrow \mathcal{O}(-1)
    \\(x&,y)\mapsto (x,W_Uy)_U,
\end{align*}
which is $\overline{\bbT}{}^Q_{\mathcal{O}}(U)$-equivariant. 
%We denote by the same symbol the induced perfect pairing on $M_Q(U)$.
\subsubsection{}\label{Shimura set} Now assume that we are in the \emph{definite case}. For $U\in \mathcal{J}_Q$, we define 
$$ Y_Q(U)\colonequals G_{Q}(F)\backslash G_{Q}(\mathbb{A}_{F,f})/U, $$
which is a finite set. As in the indefinite case, we have an action of $G_Q(\mathbb{A}_{F}^{a,\infty})$ on the projective system $\{Y_Q(U)\}_{U\in \mathcal{J}_Q}$. We define the space of automorphic forms for $G_Q$ of parallel weight $2$ and of level $U$ to be the finite free $\mathcal{O}$-module
$$ H_Q(U)=H^0(Y_Q(U),\mathcal{O})=\left\{ f\colon Y_Q(U)\longrightarrow \mathcal{O}\right\}. $$
For $U,U'\in \mathcal{J}_Q$ and $g\in G_Q(\bbA^{a,\infty}_F)$, we can define a Hecke correspondence
$$ [UgU']\colon H_Q(U')\longrightarrow H_Q(U) $$
given by $f\mapsto \sum_i f(\cdot\ u_ig)$, where $U=\bigsqcup_i u_i (U\cap gU'g^{-1})$. We can also define a perfect pairing
\begin{equation}\label{pairingdefinite}
\begin{split}
        (\cdot,\cdot)_U \colon H_Q(U)\times &H_Q(U) \longrightarrow \mathcal{O}
    \\ (f&,h) \mapsto \sum_{y\in Y_Q(U)} \frac{1}{\omega_y} f(y)h(y),
\end{split}
\end{equation}
where $\omega_y= | y^{-1}D_Q^\times y \cap U|$. By calculations similar to \cite[Proposition 3]{Newton11}, this pairing satisfies the same relations as the Poincaré pairing (\ref{PoincaréPairing}). 

\begin{remark}\label{wy1}
    A priori, the groups $y^{-1}D_Q^\times y \cap U$ need not be trivial. However  by \cite[Lemma 4.2]{Tho16}, an element $\gamma \in y^{-1}D_Q^\times y \cap U$ must lie in $F^\times$. Hence $\gamma$ is a root of unity that moreover satisfies the congruence $\gamma\equiv 1 \bmod \varpi_a \mathcal{O}_{F_a}$. This implies that either $|\mathbb{N}_{F/\bbQ}(\gamma-1)|\ge q_a$, or $\gamma=1$.  Our hypothesis that $q_a> |\mathbb{N}_{F/\bbQ}(\zeta-1)|$ for any non-trivial root of unity $\zeta$ of $F$ then forces $\gamma=1$. Thus we have that these groups are trivial, i.e., $\omega_y=1$.
\end{remark}

Setting $S=Q\cup\Sigma(U)\cup \Sigma_p$ and $S'\subseteq\Sigma(U)$ as before, we get an action of $\overline{\bbT}{}^{S,S'}_{\mathcal{O}}$ on $H_Q(U)$, and we denote the corresponding faithful algebras by $\bbT^Q_{\mathcal{O}}(U)$ and $\overline{\bbT}{}^Q_{\mathcal{O}}(U)$. If $U$ is a congruence subgroup, we also modify the pairing \eqref{pairingdefinite} by setting
$$ \langle f,h\rangle_U\colonequals (f, W_U h)_U.$$
with $W_Uh(y)=h(y^{\alt}w_U)$. We then get a $\overline{\bbT}{}^Q_{\mathcal{O}}(U)$-equivariant pairing.

\subsubsection{} \noindent In both the definite and indefinite case, there exists a Galois representation valued in the corresponding Hecke algebra, as recorded in the following proposition (see \cite[(3.4.1)]{Kis09} and \cite[Proposition 4.7]{Tho16} which uses the geometric Frobenius normalization).
\begin{proposition}\label{Tgaloisrep}
We have the following properties
\begin{itemize}
    \item[(1)] There exists a continuous representation $\rho_{\mathfrak{m}}\colon G_F \rightarrow \GL_2(\mathbb{T}^{Q}_{\mathcal{O}}(U)_{\mathfrak{m}})$ lifting $\overline{\rho}$ and satisfying the following condition: for each finite place $v\not\in S$ of $F$, $\rho_{\mathfrak{m}}|_{G_{v}}$ is unramified and $\rho_{\mathfrak{m}}(\Frob_v)$ has characteristic polynomial
    $$ \mathrm{char}\big(\rho_{\mathfrak{m}}(\Frob_v),X\big)= X^2-T_vX+q_vS_v. $$
    \item[(2)] If $Q$ is of indefinite type, then there exists a finite $\mathbb{T}^Q_{\mathcal{O}}(U)_{\mathfrak{m}}$-module $M_Q(U)$, together with an isomorphism of $\mathbb{T}^{Q}_{\mathcal{O}}(U)_{\mathfrak{m}}[G_F]$-modules
    $$ H_Q(U)_{\mathfrak{m}} \cong \rho_{\mathfrak{m}}(-1)\otimes_{\mathbb{T}^Q_{\mathcal{O}}(U)_{\mathfrak{m}}} M_Q(U). $$
 \end{itemize}
\end{proposition}

If $Q$ is of definite type, then we set $M_Q(U)\colonequals H_Q(U)_{\mathfrak{m}}$. 

In both cases, $M_Q(U)$ is self-dual with respect to the (induced) pairing introduced above.

\subsubsection{}  Later we will need to work with a fixed central character version of these objects which we now define. Consider the finite abelian group $C_U=F^\times \backslash \mathbb{A}^\times_{F,f}/(U\cap \mathbb{A}_{F,f}^\times)$ with its action on $H_Q(U)$, where for a finite place $v\not\in  S$, the image of $\varpi_v\in F^\times_v\subset \mathbb{A}_{F,f}^\times$ acts as $S_v$. In fact, the group algebra $\mathcal{O}[C_U]$ can be identified with a subalgebra of $\overline{\mathbb{T}}{}^{Q}_{\mathcal{O}}(U)$.
Let $C_{U,p}$ be the $p$-Sylow subgroup of $C_U$, then setting $\mathfrak{m}'=\mathfrak{m}\cap \mathcal{O}[C_U]$, we get that $\mathcal{O}[C_U]_{\mathfrak{m'}}\cong \mathcal{O}[C_{U,p}]$. We thus have an embedding $\mathcal{O}[C_{U,p}]\hookrightarrow \overline{\mathbb{T}}{}^{Q}_{\mathcal{O}}(U)_{\mathfrak{m}}$, which turns $M_Q(U)$ into a finite free $\mathcal{O}[C_{U,p}]$-module.

Given a character $\phi_{\iota}\colon C_U\to \mathcal{O}$, we obtain by extension an augmentation $\phi_\iota\colon \mathcal{O}[C_{U,p}]\rightarrow \mathcal{O}$. We write $\mathfrak{J}_\phi$ for its kernel, and define $M_{Q}^\phi(U)$ and $\overline{\mathbb{T}}{}^{Q}_{\mathcal{O},\phi}(U)_{\mathfrak{m}}$ to be the quotients of $M_Q(U)$ and $\overline{\mathbb{T}}{}^{Q}_{\mathcal{O}}(U)_{\mathfrak{m}}$ by $\mathfrak{J}_\phi$. 

\subsection{Level structure}\label{levelstructure}
We introduce the level structures used throughout this paper. To do so, we analyze how the tame conductor of $\rho_{\mathbf{f}}$ degenerates upon reduction modulo $\varpi$. This will also allow us to impose the appropriate deformation conditions on the local Galois deformation rings.

\subsubsection{} We recall a few facts about the Artin conductor of $\rho_{\mathbf{f}}$ (see \cite[\S 1]{Liv89} for more details).  Let $v\nmid p$ be a finite place and consider the ramification filtration $$G_0=I_{v}\supseteq G_1\supseteq \cdots \supseteq G_i \supseteq \cdots$$ of $G_{v}$. Let $\rho: G_{v}\rightarrow \Aut(V)$ be a continuous representation acting on an $E$ (or $k$)-vector space $V$.
Writing $(\rho^{\text{ss}},V^{\text{ss}})$ for its semisimplification, we know by Grothendieck's Monodromy theorem that $\rho^{\text{ss}}$ is trivial on an open compact subgroup of $I_{v}$. The swan conductor of $\rho$ is defined by the formula $$\text{sw}(\rho):=\sum_{i\ge 1} \frac{ \text{codim}((V^{\text{ss}})^{G_i})}{[G_0:G_i]}\in \mathbb{N}.$$ And the Artin conductor of $\rho$ is defined as
$$ \Art(\rho)= 2- \dim(V^{I_{v}})+ \text{sw}(\rho). $$
We write $m_v=\text{Art}(\overline{\rho}|_{G_{v}})$, and let  $\mathfrak{n}_{\overline{\rho}}=\prod_{v\nmid p} \mathfrak{p}_v^{m_v}$ denote the tame conductor of $\overline{\rho}$. Similarly, we write $\mathfrak{c}=\prod_{v\nmid p} \mathfrak{p}_v^{s_v}$, with $s_v=\Art(\psi_{\iota}|_{G_v})$.

 Since $\mathbf{f}$ is a newform, we know by \cite[Théorème A.]{Car} that $\mathfrak{n}=\prod_{v\nmid p} \mathfrak{p}_v^{m_{v}'}$, with $m_{v}'=\Art(\rho_{\mathbf{f}}|_{G_{v}})$.
Let $d_v$ (resp. $d_v'$) denote the dimension of $(\overline{\rho}|_{G_{v}})^{I_{v}}$ (resp. $(\rho_{\mathbf{f}}|_{G_{v}})^{I_{v}}$). For convenience, we set $m_v=m'_v=d_v=d'_v=0$ if $v\mid p$. Since $\text{sw}(\overline{\rho}|_{G_{v}})= \text{sw}(\rho_{\mathbf{{f}}}|_{G_{v}})$ (see \cite[Proposition 1.1]{Liv89}), we obtain $$m_v+d_v=m_v'+d_v'.$$
By the analysis on the degeneracy of the Artin conductor carried out in \cite[\S 1]{Car89} (see also \cite{DT94}), and the hypothesis $(3)$ in \Cref{maintheorem0}, we have that $m_v<m'_v$ only in the following three cases:
\begin{enumerate}
    \item[$\mathbf{(I)}$]  $\rho_{\mathbf{f}}|_{G_{v}}$ is Steinberg in the sense of \Cref{galdef} and $\overline{\rho}|_{G_{v}}$ is unramified.
    \item[$\mathbf{(II)}$] $\rho_{\mathbf{f}}|_{G_{v}}$ is Steinberg or decomposable and $\Art(\psi_{\iota}|_{G_v})=1$, while $\Art(\overline{\psi_{\iota}}|_{G_v})=0$.
    \item[$\mathbf{(III)}$] $\rho_{\mathbf{f}}|_{G_{v}}=\Ind_{H}^{G_v}\xi$, where $H$ is the absolute Galois group of an unramified quadratic extension of $F_v$, and $\xi\colon H\to \mathcal{O}^\times$ is a character that is residually unramified.
\end{enumerate}
In case (\textbf{I}), $m_v=0$ and $m'_v=1$. In case (\textbf{II}), $q_v\equiv 1 \bmod p$, $m_v=0$, and $m'_v\in \{1,2\}$. In case (\textbf{III}), $q_v\equiv-1\bmod p$, $m_v\in \{0,1\}$, $m'_v=2$, and $\det \rho_{\mathbf{f}}|_{G_v}$ is unramified. Moreover, \cite[Proposition 2.1-2.3]{Diamond97} shows that in the last case, $m_{v}=1$ can occur only if $\overline{\rho}|_{G_v}$ is Steinberg.

\subsubsection{} Let $Q$ be a finite set of finite places of $F$ with associated quaternion algebra $D_Q$. \emph{From now on}, all such $Q$ that we will consider will satisfy the following condition
\begin{enumerate}
    \item[$(\mathbf{St})$] $\mathbf{f}$ is Steinberg at $Q$. In other words, for all $ v\in Q$, $m'_v= 1$ and $s_v=0$.
\end{enumerate}
 We define the compact open subgroups $U^{\mathbf{f}},U^{\min}\in \mathcal{J}_Q$ of $G_Q(\bbA_F^\infty)$ such that for $v\not\in Q\cup\{a\}$, we have
\begin{itemize}
    \item $U^{\min}_v=\Gamma_0(v^{m_v})$,
    \item $U^{\mathbf{f}}_v=\Gamma_0(v^{m_v'})$ if $s_v=0$, and $U^{\mathbf{f}}_v=\Gamma_1(v^{m'_v})$ if $s_v\neq 0$.
\end{itemize}
\begin{remark}
    To simplify the exposition, we omit $Q$ from the notation of these subgroups, despite their dependence on it.
\end{remark}

Finally, we take $S'\subseteq \Sigma(U^{\mathbf{f}})$ to be the set $\Sigma_{\uni}$ of finite places $v\not\in Q$ such that $m_v=0$, while $m_v'=1$. Consequently, we will work with the decorated algebra Hecke algebra $\overline{\mathbb{T}}{}^{Q}_{\mathcal{O}}(U^{\mathbf{f}})_{\mathfrak{m}}$ at level $U^{\mathbf{f}}$, and with the anemic Hecke algebra $\mathbb{T}^{Q}_{\mathcal{O}}(U^{\min})_{\mathfrak{m}}$ at level $U^{\min}$.

By the Jacquet--Langlands correspondence, the hypothesis $(\textbf{St})$ implies that the augmentation $\lambda_{\mathbf{f}}$ factors through the algebra $\overline{\mathbb{T}}{}^{Q}_{\mathcal{O}}(U^{\mathbf{f}})_{\mathfrak{m}}$. We denote the induced augmentation by the same symbol
$$ \lambda_{\mathbf{f}}\colon \overline{\mathbb{T}}{}^{Q}_{\mathcal{O}}(U^{\mathbf{f}})_{\mathfrak{m}}\longrightarrow \mathcal{O}. $$

\begin{remark}\label{mult1overE}
    If $\mathcal{O}$ is sufficiently large, then $\overline{\bbT}{}^Q_{\mathcal{O}}(U^{\mathbf{f}})_{\mathfrak{m}}[\tfrac{1}{\varpi}]$ decomposes as a product of copies of $E$, indexed by the Galois conjugacy classes of tuples $(\mathbf{f}',(a_v)_{v\in S'})$. Here $\mathbf{f}'$ ranges over eigenforms of level $U^{\mathbf{f}}$ that are Steinberg at $Q$ and satisfy $\overline{\rho}=(\rho_{\mathbf{f}'}\bmod \varpi$), and $a_v$ denotes an eigenvalue of $\mathbf{U}_v$ acting on $M_Q(U^{\mathbf{f}})[\lambda_{\mathbf{f}'}]$. By the Eichler-Shimura relation (as in \cite[(2.3)]{Tho16} and \cite[\S 2.3]{Car} for the indefinte case and \cite[Lemma 1.3]{Tay06} for the definite case) and multiplicity one, the space $M_Q(U^{\mathbf{f}})[\tfrac{1}{\varpi}]$ admits the same description. It follows that $M_Q(U^{\mathbf{f}})[\tfrac{1}{\varpi}]$ is free of rank one over $\overline{\bbT}{}^Q_{\mathcal{O}}(U^{\mathbf{f}})_{\mathfrak{m}}[\tfrac{1}{\varpi}]$. In fact, the reason for adding the operators $\mathbf{U}_v$ is precisely to ensure this generic multiplicity-one property (after localization, this is only an issue when $\overline{\rho}$ is unramified at $v$ and  $\overline{\rho}(\Frob_v)$ has a multiple eigenvalue).
 \end{remark}

The level-raising argument in \Cref{galdef} requires us to work with a Hilbert modular form that is congruent to $\mathbf{f}$ and has minimal level. This form is provided by the following lemma.

\begin{lemma}\label{formg}
    There exists a Hilbert cuspdial newform $\mathbf{g}$ of parallel weight $2$, level $\mathfrak{n}_{\overline{\rho}}$, and nebentypus $\psi_0$,  such that $\rho_{\mathbf{g}}\bmod \varpi \cong \overline{\rho}$ (up to enlarging $E$, we can assume $K_{\mathbf{g}}\subseteq E$). In particular, this gives an augmentation
\begin{equation*}
    \lambda_{\mathbf{g}}\colon \mathbb{T}^{Q}_{\mathcal{O},\psi_{0}}(U^{\min})_{\mathfrak{m}}\longrightarrow \mathcal{O}. 
\end{equation*}
\end{lemma}
\begin{proof}
    Applying the level lowering results of \cite[Theorem 1.1]{Fuj1}, \cite{Jar99}, and \cite{Jar04}, we get the desired newform $\mathbf{g}$ up to the places $v$ such that $q_v\equiv 1 \bmod p$, in which case $\mathbf{g}$ could be unramified or Steinberg at $v$. We deal with these places using the main theorem of   \cite{Raj01}, which applies at the places $q_v\equiv 1 \bmod p$ and Iwahori level.
\end{proof}
Note that in general the nebentypus  $\psi_0$ of $\mathbf{g}$ differs from the nebentypus $\psi$ of $\mathbf{f}$. 

\subsection{Jacobian of Shimura curves}\label{jacshimuracurves}$ $
Let $Q$ be a finite set of places of $F$ of indefinite type, and $U\in \mathcal{J}_Q$. We consider the Jacobian associated to the Shimura curve $X_Q(U)$,
\begin{equation*}
J_Q(U)\colonequals\Pic^0(X_Q(U)/F),
\end{equation*}
which is an abelian variety defined over $F$. Since $\Pic^0(-/F)$ is both covariant and contravariant for finite flat morphisms of proper smooth curves over $F$, for each $U,U'\in \mathcal{J}_Q$ and $g\in G_Q(\mathbb{A}_F^f)$, we obtain a morphism of abelian varieties
\begin{equation*}
    [UgU']\colon J_Q(U') \longrightarrow J_Q(U).
\end{equation*}
Therefore, the algebra $\bbT^{S}$ (resp. $\overline{\bbT}{}^{S,S'}$) acts on $J_Q(U)$, and we denote by $\bbT^Q(U)$ (resp. $\overline{\bbT}{}^Q(U)$) its image in the endomorphism algebra of $J_Q(U)$.

The Jacobian $J_Q(U)$ carries a canonical principal polarization $\theta\colon \widehat{J_Q}(U) \xrightarrow{\sim} J_Q(U) $, and the corresponding Rosati involution $(-)^\dagger$ satisfies
\begin{equation*}
    [UgU]^\dagger=[Ug^{-1}U],
\end{equation*}
for $g\in G_Q(\mathbb{A}_F^f)$.

If $U$ is a congruence subgroup, then the Atkin--Lehner involution $W_U$ on $X_Q(U)$ induces an involution on $J_Q(U)$, which we denote by the same symbol. Subsequently, we can define a new principal polarization
$$\varphi\colonequals W_U\circ \theta \colon \widehat{J_Q}(U)\longrightarrow J_Q(U) $$
that is defined over a finite extension of $F$. This polarization has the advantage of being $\overline{\bbT}{}^Q(U)$-equivariant.

\section{Galois deformation theory, R=$\bbT$, and Gorensteinness}\label{galdef}
The goal of this section is to prove the freeness of $M_Q(U^{\mathbf{f}})$ over the Hecke algebra $\overline{\bbT}{}^{Q}_{\mathcal{O}}(U^{\mathbf{f}})_{\mathfrak{m}}$ in the case where $Q$ does not contain trivial places. This is the content of \Cref{freenessoverGorenstein}, which is achieved by upgrading a freeness result at the minimal level $U^{\min}$ to the level $U^{\mathbf{f}}$ through a careful study of level raising congruences. Conversely, if $Q$ contains at least a trivial place, we show in \Cref{notGorenstein} that $M_Q(U^{\mathbf{f}})$ is $\emph{not}$ free over $\overline{\bbT}{}^{Q}_{\mathcal{O}}(U^{\mathbf{f}})_{\mathfrak{m}}$ by proving that the latter is \emph{not} Gorenstein.

\subsection{Congruence modules} We introduce congruence modules and the Wiles defect, and collect some of their properties that will be used throughout the paper.

\begin{notation}
Given an $\mathcal{O}$-module $M$, we write $M_{\tors}$ for its $\varpi$-torsion submodule, and $M^{\tf}$ for its $\varpi$-torsion free quotient.  We also set 
$$ M^*\colonequals \Hom_{\mathcal{O}}(M,\mathcal{O}). $$
\end{notation}
Recall that $\CNL_\mathcal{O}$ is the category of complete noetherian local $\mathcal{O}$-algebras with residue field $k$.
\begin{definition}
    Given $R\in \CNL_{\mathcal{O}}$ and an augmentation $\lambda\colon R\to \mathcal{O}$, we define the cotangent module of $R$ with respect to $\lambda$ to be $\Phi_{\lambda,R}\colonequals(\mathfrak{p}_\lambda/\mathfrak{p}_{\lambda}^2)$, where $\mathfrak{p}_\lambda \colonequals \ker(\lambda)$. We also let $\Phi_{\lambda,R}^{\tors}$ to be its torsion part.
    
     Given a map $\pi\colon S\to R$ in $\CNL_{\mathcal{O}}$, we define the relative cotangent module $\Phi_{\lambda,S/R}$ to be the kernel of the surjection $\Phi_{\lambda\circ \pi,S}\twoheadrightarrow \Phi_{\lambda,R}$. 
\end{definition}
\begin{remark}\label{torsioncotangentspace}
    If $S\twoheadrightarrow R$ is surjective with $S$ and $R$ being of relative dimension $d$ and generically formally smooth at $\lambda$, then we can write $\Phi_{\lambda,R}\cong \mathcal{O}^d\oplus \Phi_{\lambda,R}^{\tors} $, and similarly for $S$ (see \cite[Lemma 2.7]{BIK}). Therefore, the cotangent module is the kernel of the surjective map
    $$ \Phi_{\lambda,S/R}=\ker\big( \Phi_{\lambda\circ \pi,S}^{\tors}\twoheadrightarrow \Phi_{\lambda,R}^{\tors}\big).$$
\end{remark}
Let $R$ be a reduced finite flat local $\mathcal{O}$-algebra. This means that the map $R\to R\otimes_{\mathcal{O}}E=R[\tfrac{1}{\varpi}]$ is injective and $R[\tfrac{1}{\varpi}]$ is a product of finite extensions of $E$. Let us fix an augmentation $\lambda\colon R \to \mathcal{O}$, and set
$$I_\lambda\colonequals R[\mathfrak{p}_\lambda].$$
Then we have a decomposition $R[\tfrac{1}{\varpi}]=\mathfrak{p}_\lambda[\tfrac{1}{\varpi}]\oplus I_\lambda[\tfrac{1}{\varpi}]$. We also note that $\Phi_{\lambda,R}$ is of finite length.

We adopt the definition of the congruence module over such an $\mathcal{O}$-algebra $R$ given in \cite{KIM}.
\begin{definition}
    Let $M$ be an $\mathcal{O}$-torsion free finitely generated $R$-module. The \emph{congruence module} of $M$ with respect to $\lambda$ is defined to be
    $$\Psi_\lambda(M)\colonequals \coker\left( M[\mathfrak{p}_\lambda]\longrightarrow (M/\mathfrak{p}_\lambda M)^{\tf}\right). $$
    The \emph{Wiles defect} of $M$ is the integer
    $$ \delta_\lambda(M)\colonequals (\rank_{\mathcal{O}}M[\mathfrak{p}_\lambda])\cdot \ell_{\mathcal{O}}(\Phi_{\lambda,R})-\ell_{\mathcal{O}}(\Psi_\lambda(M)). $$
\end{definition}
\begin{example}\label{congmoduleforR}
    If $M=R$, then $\Psi_{\lambda}(R)=\mathcal{O}/\lambda(I_{\lambda})$.
\end{example}
\begin{remark}\label{upperbound}
    If $\rank_{\mathcal{O}}M[\mathfrak{p}_\lambda]=1$, then the inequality $\ell_{\mathcal{O}}(\Psi_\lambda(M))\le \ell_{\mathcal{O}}(\Psi_\lambda(R))$ always holds (see \cite[Theorem 3.12]{BKM1}).
\end{remark}
The Wiles defect gives us a criterion to prove freeness of modules over Gorenstein rings as illustrated by the following theorem.
\begin{theorem}[{\cite[Theorem 1.2]{BIK}}]\label{freenessoverGorenstein}
    Suppose that $R$ is Gorenstein and that $M$ is an $\mathcal{O}$-torsion free finitely generated $R$-module such that $\rank_{\mathcal{O}}M[\mathfrak{p}_\lambda]=1$. If
    $ \delta_\lambda(M)= \delta_\lambda(R),$
    or equivalently $$\ell_{\mathcal{O}}(\Psi_{\lambda}(M))=\ell_{\mathcal{O}}(\Psi_{\lambda}(R)),$$ and if $M$ has rank at most one at each generic point of $R$, then $M$ is free over $R$.\end{theorem}
\begin{remark}
    Note that the last condition in the statement of the theorem holds for example if $M[\tfrac{1}{\varpi}]$ is free of rank one over $R[\tfrac{1}{\varpi}]$.
\end{remark}
\begin{remark}\label{dualofGorenstein}
    Throughout this paper, we will extensively use the fact that if $R$ is a finite flat local $\mathcal{O}$-algebra, then $R$ is Gorenstein if and only if $R^*=\Hom_{\mathcal{O}}(R,\mathcal{O})$ is a free $R$-module of rank one (see \cite[Theorem 21.15]{Eisen}).
\end{remark}

\begin{lemma}\label{quotientsiso}
    Let $M$ be an $\mathcal{O}$-torsion free finitely generated $R$-module, then we have the following canonical isomorphisms
    $$ (M/\mathfrak{p}_\lambda M)^{\tf} \cong M/(M[I_\lambda]) \quad \text{ and, } \quad (M/I_\lambda M)^{\tf}\cong M/(M[\mathfrak{p}_\lambda]).$$
    We also have canonical isomorphisms
    $$ (M[\mathfrak{p}_\lambda])^*\cong M^*/(M^*[I_\lambda]) \quad \text{ and, } \quad M^*[\mathfrak{p}_\lambda]\cong (M/\mathfrak{p_\lambda}M)^*. $$
\end{lemma}

\begin{proof} The first two isomorphisms follow from the fact that $M[\tfrac{1}{\varpi}]=\mathfrak{p}_\lambda M[\tfrac{1}{\varpi}]\oplus I_\lambda M[\tfrac{1}{\varpi}]$.
For the third isomorphism, let us consider the exact sequence
$$ 0 \to (M/M[\mathfrak{p}_{\lambda}])^* \to M^* \to (M[\mathfrak{p}_\lambda])^*\to 0. $$
We have $\mathfrak{p}_\lambda M^*\subseteq (M/M[\mathfrak{p}_{\lambda}])^* \subseteq M^*[I_\lambda]$. Since $\mathfrak{p}_\lambda M^*$ and $M^*[I_\lambda]$ have the same $\mathcal{O}$-rank, we find that $(M/M[\mathfrak{p}_{\lambda}])^* = M^*[I_\lambda]$ since both are saturated submodules of $M^*$.
\end{proof}
The article \cite{BKM2} extends the notion of Wiles defect to the higher dimensional setting. More specifically, let $R\in \CNL_{\mathcal{O}}$, which is Cohen--Macaulay and flat over $\mathcal{O}$ of relative dimension $d$, and let $\lambda\colon R\to \mathcal{O}$ be an augmentation such that $ R[\tfrac{1}{\varpi}]$ is formally smooth at the point corresponding to $\lambda$ (we retain this notation for the rest of this section). By \cite[Proposition 3.3]{BKM2}, there exists $\widetilde{R}\in \CNL_{\mathcal{O}}$ which is complete intersection, flat and equidimensional over $\mathcal{O}$ of relative dimension $d$, together with a continuous surjection of $\mathcal{O}$-algebras $\pi \colon \widetilde{R}\to R$, such that $\widetilde{R}[\tfrac{1}{\varpi}]$ is formally smooth at $\widetilde{\lambda}\colonequals \lambda\circ \pi$. We set $I=\ker(\pi)$. Given this data, the Wiles defect of $R$ is defined as follows.
\begin{definition}\label{originaldefwilesdefect}
Consider the following quantities
   \begin{itemize}
       \item $D_{1,\lambda}(R)\colonequals \widehat{\Der}^1_{\mathcal{O}}(R,E/\mathcal{O})$ is the $1^{st}$ continuous André--Quillen cohomology group (see \cite[\S 3.3]{BKM2}), 
       \item $c_{1,\lambda}(R)\colonequals \widetilde{\lambda}(\widetilde{R}[I])/\widetilde{\lambda}(\Fitt_0(I))$.
   \end{itemize}
   Then the \emph{Wiles defect} of $R$ at $\lambda$ is defined to be
   $$ \delta_\lambda(R)= \ell_{\mathcal{O}}\left(D_{1,\lambda}(R)\right)- \ell_{\mathcal{O}}\left(c_{1,\lambda}(R)\right), $$
   and is independent of the choice of $\widetilde{R}$.
\end{definition}
We record the following important generalization of Wiles' numerical criterion, which shows that the Wiles defect precisely measures the failure of being complete intersection.
\begin{proposition}[{\cite[Proposition 3.28]{BKM2}}]\label{WDcompleteIntersection}
      $R$ is a complete intersection if and only if $\delta_{\lambda}(R)=0$. In particular, $\delta_{\lambda}(\mathcal{O}[[x_1,\dots,x_m]])=0$ for any $m\ge 1$, and augmentation $\mathcal{O}[[x_1,\dots,x_m]]\to \mathcal{O}$.
\end{proposition}

The Wiles defect enjoys two further important properties that we will use extensively later. The first is its invariance under quotienting by a regular sequence, and the second is its additive property on completed tensor products. For completeness, we record both properties below.
\begin{proposition}[{\cite[Theorem 3.29]{BKM2}}]\label{invunderregseq}
   Assume that $(\varpi,\theta_1,\dots,\theta_i)$ is a regular sequence in $R$, with $\theta_1,\dots,\theta_i\in \ker \lambda$. Let $R_{\theta}=R/(\theta_1,\dots,\theta_i)$ equipped with the induced augmentation $\lambda_{\theta}$. Then 
    $$ \delta_\lambda(R)=\delta_{\lambda_{\theta}}(R_{\theta}). $$
\end{proposition}

\begin{proposition}[{\cite[Proposition 3.30]{BKM2}}]\label{invtensorprod}
    Let $R_1,R_2\in \CNL_{\mathcal{O}}$ be Cohen-Macaulay and flat over $\mathcal{O}$ of relative dimension $d_1$ and $d_2$. Let $\lambda_i\colon R_i\to \mathcal{O}$ be augmentations such that $R_i[\tfrac{1}{\varpi}]$ is formally smooth at the point corresponding to $\lambda_i$. Let $R=R_1\widehat{\otimes}_{\mathcal{O}}R_2$, and $\lambda\colonequals \lambda_1\widehat{\otimes}\lambda_2 \colon R\to \mathcal{O}$, then we have the equality
    $$ \delta_{\lambda}(R)=\delta_{\lambda_1}(R_1)+\delta_{\lambda_2}(R_2).$$
\end{proposition}

For the remainder of this paper, we will not use the formula in \Cref{originaldefwilesdefect} for the Wiles defect. We will instead introduce an alternative one, which is given in the following lemma.

\begin{lemma}\label{formulaWD}
    With the same notation as above, we have
\begin{equation*}
    \delta_{\lambda}(R)=\ell_{\mathcal{O}}\big(\mathcal{O}/\widetilde{\lambda}_\theta(\widetilde{R}_{\theta}[I_\theta])\big)-\ell_{\mathcal{O}}(\Phi_{\lambda,\widetilde{R}_\theta/R_\theta}),
\end{equation*}
where $I_{\theta}=\ker(\widetilde{R}_{\theta}\to R_{\theta}).$
\end{lemma}
\begin{proof}
    Suppose that $(\varpi,\theta_1,\dots,\theta_d)$ is a regular sequence in $R$, with $\theta_1,\dots,\theta_d\in \ker \lambda$. Then the algebra $R_\theta=R/(\theta_1,\dots,\theta_d)$ is finite over $\mathcal{O}$, and we write $\lambda_\theta$ for the induced augmentation. By \Cref{invunderregseq}, we have 
$$ \delta_\lambda(R)=\delta_{\lambda_\theta}(R_\theta)= \ell_{\mathcal{O}}(\Phi_{\lambda_\theta,R_\theta})-\ell_{\mathcal{O}}(\Psi_{\lambda_\theta}(R_\theta)).$$
By \cite[Lemma 3.7]{BKM2}, the elements $\theta_1,\dots,\theta_d$ lift to elements $\widetilde{\theta}_1,\dots,\widetilde{\theta}_d$ such that $(\varpi, \widetilde{\theta}_1,\dots,\widetilde{\theta}_d)$ is a regular sequence on $\widetilde{R}$. Letting $\widetilde{R}_\theta=\widetilde{R}/(\widetilde{\theta}_1,\dots,\widetilde{\theta}_d)$ with its augmentation $\widetilde{\lambda}_\theta$, we get by \Cref{WDcompleteIntersection} that $\delta_{\widetilde{\lambda}_\theta}(\widetilde{R}_\theta)=0$, and so
$$ \delta_{\lambda_\theta}(R_\theta)=\delta_{\lambda_\theta}(R_\theta)-\delta_{\widetilde{\lambda}_\theta}(\widetilde{R}_\theta)=\big(\ell_{\mathcal{O}}(\Phi_{\lambda_\theta, R_\theta})-\ell_{\mathcal{O}}(\Phi_{\lambda_\theta,\widetilde{R}_\theta})\big)-\big(\ell_{\mathcal{O}}(\Psi_{\lambda_\theta}(R_\theta))-\ell_{\mathcal{O}}(\Psi_{\widetilde{\lambda}_\theta}(\widetilde{R}_\theta))\big).$$
By \cite[Lemma A.10]{FKR}, we have that $\ell_{\mathcal{O}}(\Psi_{\lambda_\theta}(R_\theta))-\ell_{\mathcal{O}}(\Psi_{\widetilde{\lambda}_\theta}(\widetilde{R}_\theta))=-\ell_{\mathcal{O}}(\mathcal{O}/\widetilde{\lambda}_\theta(\widetilde{R}_\theta[I_\theta]))$. Therefore, we get the stated formula.
\end{proof}

\subsection{Local deformation rings}\label{localdefring} We introduce the local deformation rings that we will work with, recall their relevant properties, and compute their Wiles defect in the cases not treated in \cite{BKM2}. In particular, \Cref{cotangentuni} addresses the case of the unipotent deformation ring at an unramified augmentation.

We recall that for a finite place $v$ of $F$, $\sigma_v\in I_v$ is a fixed lift of a generator of the $\bbZ_p$-quotient of $I_v$, and $\Frob_v\in G_v$ is a fixed lift of the arithmetic Frobenius. Given $A\in \CNL_{\mathcal{O}}$ with maximal ideal $\mathfrak{m}_A$, the multiplicative group $(1+\mathfrak{m}_A,\times)$ is a pro-$p$ group and $p$ is odd. Therefore, we can define a square root $(\ \cdot \ )^{1/2}\colon (1+\mathfrak{m}_A,\times) \to (1+\mathfrak{m}_A,\times)$.

We now introduce the deformation functors of our fixed residual representation $\overline{\rho}=(\rho_{\mathbf{f}} \bmod \varpi)$ of \Cref{setting}. Let $\psi_0\colon \mathbb{A}_F^\times/F^\times\to \mathbb{C}^\times$ be the nebentypus associated to the Hilbert modular form $\mathbf{g}$ of \Cref{formg}, and let $\psi_{0,\iota}\colonequals \iota\circ \psi_{0}$. For a finite place $v$ of $F$, define the framed deformation functor $
\mathcal{D}^\square_v\colon \CNL_\mathcal{O} \rightarrow \Sets$ by assigning to each $A\in \CNL_{\mathcal{O}}$ the set of all continuous morphisms $\rho_A\colon G_{v}\rightarrow \GL_2(A)$ satisfying $\rho_A \bmod \mathfrak{m}_A=\overline{\rho}|_{G_{v}}$, \emph{and} $\det \rho_A=(\epsilon_p\cdot \psi_{0,\iota}\circ \Art_F^{-1})|_{G_v}$.
This functor is representable by a ring $R^\square_v\in \CNL_{\mathcal{O}}$ admitting a universal lift $\rho^\square\colon G_{v}\rightarrow \GL_2(R_v^\square)$.

We next impose local deformation conditions. The resulting lifting ring is denoted $R_v^{\tau_v}$, where the superscript $\tau_v\in \{\text{fl,min,st,uni,}\square\}$ specifies the condition type, and $\rho_v^{\tau_v}$ denotes the corresponding universal lifting.

We fix a finite set $\Sigma$ of finite places of $F$, decomposed as a union of the following finite sets:
\begin{itemize}
    \item $\Sigma_p$ is the set of places $v\mid p$.
    \item $\Sigma_{\st}= Q$ is a set of finite places $v$ such that $m'_v=1$ and $s_v=0$.
    \item $\Sigma_{\uni}$ is the set of places $v\not\in Q$ such that $m_v=0$, while $m_v'=1$.
    \item $\Sigma_\square$ is the set of places $v$ such that $m_v=0$ but $m'_v=2$.
      \item $\Sigma_{\triangle}$ is the set of places $v$ such that $m_v=1$ but $m'_v=2$ (we must have $q_v\equiv -1\bmod p$).
     \item $\Sigma_{\text{min}}$ is the set of places $v\not\in Q$ such that $m_v=m'_v$.
\end{itemize}
For $v\in \Sigma_p$, the extension $F_v/\mathbb{Q}_p$ is unramified by assumption. Since $\rho_{\mathbf{f}}|_{G_{v}}$ is crystalline with $\tau$-Hodge-Tate weights $\{0,1\}$ for all $v\in \Sigma_p$, the residual representation $\overline{\rho}|_{G_{v}}$ is flat. Therefore, the Fontaine--Laffaille theory applies, and we define
\begin{itemize}
    \item $R_v^{\text{fl}}$ to be the quotient of $R_v^{\square}$ parameterizing flat deformations.
\end{itemize}
For $v\not\in \Sigma_p$,
\begin{itemize}
    \item We let $R_v^{\min}$ be the quotient of $R_v^{\square}$ parameterizing deformations $\rho$ of  $\overline{\rho}|_{G_v}$ such that $\rho\otimes(\det \rho)^{-1/2}$ is minimally ramified in the sense of \cite[Definition 2.4.14]{CHT}. Explicitly, for any $A\in \CNL_{\mathcal{O}}$, a representation $\rho_A\in \mathcal{D}_v^{\square}(A)$ is induced from a map $R_v^{\min}\to A$ if and only if the natural map
    $$ \ker\big(\rho_A\otimes(\det \rho_A)^{-1/2}(\sigma_v)-\id\big)^i\otimes_A k\longrightarrow \ker\big(\overline{\rho}\otimes(\det \overline{\rho})^{-1/2}(\sigma_v)-\id\big)^i $$
    is an isomorphism for all $i\ge 1.$
    If $\overline{\rho}|_{G_v}$ is unramified, we may also write $R^{\unr}_v$ for this ring.
\end{itemize}
For $v\in \Sigma_{\st}\cup\Sigma_{\uni}\cup \Sigma_{\triangle}$, the representation $\overline{\rho}|_{G_{v}}$ is of the form
\begin{align*}
    \begin{pmatrix}
        \epsilon_p \overline{\chi} & *
        \\ 0 & \overline{\chi}
    \end{pmatrix}
\end{align*}
with respect to some basis $e_1,e_2$ of $k^2$ and where the character $\overline{\chi}$ is unramified. We further assume that the basis is chosen so that $*$ is zero whenever $\overline{\rho}|_{G_{v}}$ is split. Let $[\overline{\chi}]$ be the Teichmuller lift of $\overline{\chi}$, and let $\alpha=[\overline{\chi}](\Frob_v)$.
\begin{itemize}
  \item If $\overline{\rho}|_{G_v}$ is unramified, we let $\overline{R}^{\unr}_v=\big(R^{\unr}_v[X]/(\mathrm{char}(\rho^{\unr}(\Frob_v),X))\big)_{\mathfrak{m}_{\unr}}$ ($\mathfrak{m}_{\unr}$ is the maximal ideal containing $X-\alpha$) be the ring parametrizing unramified deformations of $\rho$ with a choice of an eigenvalue of the Frobenius congruent to $\alpha$ modulo $\varpi$. 
\end{itemize}
For $v\in \Sigma_{\st}$, we define the Steinberg quotient $R^{\st}_v$ of $R_v^{\square}$ as follows:
\begin{itemize}
    \item If $\overline{\rho}|_{G_v}$ is ramified, then $R^{\st}_v$ is defined to be $R_v^{\min}$.
    \item If $\overline{\rho}|_{G_v}$ is unramified, then $R^{\st}_v$ is defined to be the maximal reduced $\lambda$-torsion free quotient of $R^\square_v$ characterized by the fact that the $L$-valued points of its generic fibre, for any finite extension $L/E$, corresponds to a representation of the form
    \begin{equation*}
         \begin{pmatrix}
            \epsilon_p \chi & * \\ 0 & \chi
        \end{pmatrix}
    \end{equation*}
    and with the additional condition $\chi(\Frob_v)=\alpha$ in the case $q_v\equiv-1 \bmod p$.
\end{itemize}
For $v\in \Sigma_{\uni}$, we define
\begin{itemize}
    \item The modified unipotent quotient $R_v^{\uni}$ of $R_v^{\square}$ to be the unique reduced quotient characterized by the fact that for any finite extension $L/E$, the $L$-valued points of its generic fibre correspond to pairs $(\rho,a)$ where either $\rho$ is unramified and $a\equiv \alpha\bmod\varpi$ is an eigenvalue of $\rho(\Frob_v)$, or $\rho$ of the form
        \begin{align*}
        \begin{pmatrix}
            \epsilon_p \chi & * \\ 0 & \chi
        \end{pmatrix}
    \end{align*}
    with $\chi$ unramified and satisfying $\chi(\Frob_v)=a$ (with the additional condition that $a\equiv\alpha \bmod \varpi$ when $q_v\equiv -1 \bmod p$).
\end{itemize}
\begin{remark}
    The ring $R^{\uni}_v$ is denoted by $R^{\varphi-\uni}$ in \cite{BKM2}. We adopt this notation since we do not use the ring denoted by $R^{\uni}_v$ in \cite{BKM2}. 
\end{remark}
\begin{proposition}\cite[Proposition 4.6]{BKM2}\label{propertiesoflocaldefrings}
    We have the following properties of local deformation rings:
    \begin{enumerate}
        \item For $v\in \Sigma_p$, there is an isomorphism $R_v^{\mathrm{fl}}\cong \mathcal{O}[[x_1,\dots,x_{3+[F_v:\mathbb{Q}_p]}]]$.
        \item For $v\in \Sigma\setminus \Sigma_p$ there is an isomorphism $R_v^{\min}\cong \mathcal{O}[[x_1,x_2,x_3]]$.
        \item For $v\in \Sigma\backslash \Sigma_p$, the ring $R_v^\square$ is a complete intersection, reduced, and flat over $\mathcal{O}$ of relative dimension $3$.
        \item For $v\in \Sigma_{\st}$, the ring $R_v^{\st}$ is Cohen-Macaulay, flat of relative dimension $3$ over $\mathcal{O}$ and geometrically integral. If $v$ is not a trivial place, then $R_v^{\st}\cong \mathcal{O}[[x_1,x_2,x_3]]$. Otherwise, $R_v^{\st}$ is not Gorenstein.
        \item For each $v\in \Sigma_{\uni}\cup \Sigma_{\min}\cup \Sigma_{\st}$ such that $\overline{\rho}$ is unramified, the ring $R_v^{\uni}$ is Gorenstein, reduced, and flat over $\mathcal{O}$ of relative dimension $3$. It is complete intersection if and only if $v$ is not a trivial place for $\overline{\rho}$.
    \end{enumerate}
\end{proposition}

The following three technical lemmas compute the lengths of the relative cotangent modules for various local deformation rings. These computations are a variation of those carried out in \cite[\S 5]{BKM2}, and will be needed to compute the lenght of the congruence module of the global deformation rings. 
\begin{lemma}\label{cotangentsquare1}
   Let $\lambda\colon R^{\square}_v\to \mathcal{O}$ be an augmentation such that the corresponding representation of $\rho_\lambda\colon G_v\to \GL_2(\mathcal{O})$ is unramified with $\mathrm{char}(\rho_\lambda(\Frob_v),X)=(X-\alpha)(X-\beta)$. Then the module $\Phi_{\lambda,R^{\square}_v/R_{v}^{\unr}}$ is finite over $\mathcal{O}$ and $$\ell_{\mathcal{O}}(\Phi_{\lambda,R^{\square}_v/R_{v}^{\unr}})=\ord_{\varpi}\big((1-q_v)((\alpha+\beta)^2-\alpha\beta q_v^{-1}(1+q_v)^2)\big).$$ 
    \end{lemma}
\begin{proof}
    Let $\mathcal{R}=\mathcal{O}[[a,b,c,d,u,w,t,s]]$, and let $A,B\in \GL_2(\mathcal{R})$ be the matrices
    $$ A\colonequals \begin{pmatrix}
        \alpha+u & t+\gamma \\ s & \beta+w \\ 
    \end{pmatrix}, \quad B\colonequals \begin{pmatrix}
        1+a & b \\ c & 1+d
    \end{pmatrix}, $$
    with $\alpha,\beta\in \mathcal{O}^\times,\gamma\in \mathcal{O}$, and  $\alpha\beta=\epsilon_p\cdot\psi_{0,\iota}(\Frob_v)$. Let $\mathcal{I}^\square_v=(r_1,\dots,r_6)\subset \mathcal{R}$ be the ideal generated by the relations 
    \begin{align*}
        r_1&=\det(A)-\alpha\beta= \beta u+\alpha w+uw-st-\gamma s,
        \\  r_2&=\det(B)-1=a+d+ad-bc,
        \\ \begin{pmatrix}
            r_3 & r_4
            \\ r_5 & r_6
        \end{pmatrix}&=AB-B^{q_v} A.
    \end{align*}
    Then $R^{\square}_v=\mathcal{R}/\mathcal{I}^{\square}_v$, with universal representation is given by $\rho^{\square}_v(\Frob_v)= A \bmod \mathcal{I}^{\square}_v$, and $\rho^{\square}_v(\sigma_v)= B \bmod \mathcal{I}^{\square}_v$. Moreover, we have  $R^{\unr}_v=\mathcal{R}/\mathcal{I}^{\unr}_v$, where $\mathcal{I}^{\unr}_v=(r_1,a,b,c,d)$. 

    For $R\in \CNL_{\mathcal{O}}$ equipped with an augmentation $\lambda\colon R\to \mathcal{O}$, let $\widehat{\Omega}_{\mathcal{R}/\mathcal{O}}$ denote the module of continuous Kähler differentials (see \cite[\S 7.1]{BKM1}). By \cite[Lemma 7.3]{BKM1}, it satisfies $\widehat{\Omega}_{\mathcal{R}/\mathcal{O}}\otimes_{\lambda}\mathcal{O}\cong \Phi_{\lambda,\mathcal{R}}$.
    
     Consider the augmentation $\lambda$ given by sending $a,b,c,d,u,w,t,s$ to $0$. As in the proof of \cite[Proposition 7.9]{BKM1}, the module $\widehat{\Omega}_{\mathcal{R}/\mathcal{O}}\otimes_{\lambda}\mathcal{O}$ is free over $\mathcal{O}$ with basis $\dd a|_\lambda,\dots,\dd s|_{\lambda}$. Let $\Lambda^{\square}$ (resp. $\Lambda^{\unr}$) be the image of $\mathcal{I}^{\square}_v$ (resp. $\mathcal{I}^{\unr}_v$) under the $\mathcal{O}$-linear derivation map $\mathcal{R}\to \widehat{\Omega}_{\mathcal{R}/\mathcal{O}}\otimes_{\lambda}\mathcal{O}$. Then 
    $$ \Phi_{\lambda, R^{\square}_v/R^{\unr}_v}=\Lambda^{\unr}/\Lambda^{\square}. $$
     Let $N=B-\id$ so that $N|_\lambda=0$. We have
    \begin{align*}
        \dd r_1|_\lambda &=  \beta \dd u +\alpha \dd w - \gamma \dd s,
        \\ \dd r_2|_\lambda &=\dd a + \dd d,
        \\ \begin{pmatrix}
            \dd r_3|_\lambda & \dd r_4|_\lambda
            \\ \dd r_5|_\lambda & \dd r_6|_\lambda
        \end{pmatrix} &= \dd \left(AN- \sum_{i=1}^{q_v} \binom{q_v}{i} N^i A \right)|_\lambda =  (A \dd N - q_v \dd N A)|_\lambda
        \\ &=\begin{pmatrix}
            \alpha(1-q_v)\dd a + \gamma \dd c & (\alpha-q_v\beta)\dd b +\gamma \dd d -q_v\gamma \dd a 
            \\ (\beta-q_v\alpha) \dd c & \beta(1-q_v)\dd d - q_v\gamma \dd c
        \end{pmatrix}.
    \end{align*}
    One can verify that $\dd r_6|_\lambda= \alpha^{-1}(\gamma \dd r_5|_\lambda-\beta \dd r_3|_\lambda+ \alpha\beta(1-q_v)\dd r_2|_\lambda)$. Hence $\{\dd r_i|_\lambda\}_{1\le i \le 5}$ forms a basis of $\Lambda^{\square}$. Writing this basis in the basis $(\dd r_1|_\lambda,\dd a|_\lambda,\dd b|_\lambda, \dd d|_\lambda,\dd c|_\lambda)$ of $\Lambda^{\unr}$, we get the following matrix
       $$ \begin{pmatrix}
            1 & 0 & 0 & 0 & 0
            \\ 0 & 1 & 0 & 1 & 0
            \\ 0 & -q_v\gamma & (\alpha-q_v\beta) &  \gamma & 0
            \\ 0 & 0 & 0 & 0 & (\beta-q_v\alpha)
            \\ 0 & \alpha(1-q_v) & 0 & 0 & \gamma
        \end{pmatrix}.$$
        Computing the determinant of this matrix yields
        $$\ell_{\mathcal{O}} (\Lambda^{\unr}/\Lambda^{\square})=\ord_{\varpi}\big((1-q_v)(\beta-q_v\alpha)(\alpha-q_v\beta)\big).$$
   The result then follows from the equality $(\beta-q_v\alpha)(\alpha-q_v\beta)=(q_v+1)^2\alpha\beta-q_v(\alpha+\beta)^2$.
\end{proof}
\begin{lemma}\label{cotangentsquare2}
   Assume that $q_v\not\equiv 1\bmod p$, and that $\overline{{\rho}}|_{G_v}$ is Steinberg and ramified. Let $\lambda\colon R^{\square}_v\to \mathcal{O}$ be an augmentation such that the corresponding representation of $\rho_\lambda\colon G_v\to \GL_2(\mathcal{O})$ is Steinberg. Then the module $\Phi_{\lambda,R^{\square}_v/R_{v}^{\st}}$ is finite over $\mathcal{O}$ and $$ \ell_{\mathcal{O}}(\Phi_{\lambda,R^{\square}_v/R_{v}^{\st}})=\ord_{\varpi}(q_v^2-1).$$ 
    \end{lemma}
\begin{proof}
    By \cite[Lemma 5.4]{Sho16}, the universal representation associated to $R_v^\square$ may be written as 
    \begin{align*}
        \rho_v^{\square}(\Frob_v) &= \begin{pmatrix}
            1 & x \\ y & 1
        \end{pmatrix}^{-1} \begin{pmatrix}
            q_v \alpha (1+u) & 0
            \\ 0 & \alpha (1+w)
        \end{pmatrix} \begin{pmatrix}
            1 & x \\ y & 1
        \end{pmatrix},
        \\ \rho_v^{\square}(\sigma_v) &=  \begin{pmatrix}
            1 & x \\ y & 1
        \end{pmatrix}^{-1} \begin{pmatrix}
            1+a & \gamma+ b
            \\ c & 1+d
        \end{pmatrix} \begin{pmatrix}
            1 & x \\ y & 1
        \end{pmatrix},
    \end{align*}
such that $\alpha,\gamma\in \mathcal{O}^\times$ and $q_v\alpha^2=\epsilon_p\cdot\psi_{0,\iota}(\Frob_v)$. 

As in the previous lemma, let $\mathcal{R}=\mathcal{O}[[a,b,c,d,u,w,x,y]]$, and let $\mathcal{I}^\square_v=(r_1,\dots,r_6)\subseteq \mathcal{R}$ be the ideal generated by the relations    \begin{align*}
        r_1 &=\det(\rho_v^{\square}(\Frob_v))-q_v\alpha^2= \alpha^2q_v (u+w+uw),
        \\  r_2 &=\det(\rho_v^{\square}(\sigma_v))-1= a+d+ ad-c(\gamma+b),
        \\ \begin{pmatrix}
            r_3 & r_4
            \\ r_5 & r_6
        \end{pmatrix} &=\rho_v^{\square}(\Frob_v)\rho_v^{\square}(\sigma_v)-\rho_v^{\square}(\sigma_v)^{q_v} \rho_v^{\square}(\Frob_v),
    \end{align*}
    so that $R_v^\square= \mathcal{R}/\mathcal{I}^\square_v$. The Steinberg deformation ring can be presented as $R^{\st}_v=R^{\min}_v=\mathcal{R}/\mathcal{I}^{\st}_v$, where the ideal $\mathcal{I}^{\st}_v$ is generated by the relations $r_1,\dots,r_6$ together with $r_7=a+d$, the latter arising from the minimally ramified condition $\det(\rho^{\st}_v(\sigma_v)-\id)=0$.
    Now consider the augmentation $\lambda\colon \mathcal{R}\to \mathcal{O}$ sending $x,y,u,w,a,b,c,d$ to $0$. Similarly to the proof of the previous lemma, we compute
\begin{align*}
    \dd r_1|_\lambda &= \alpha^2 q_v (\dd u + \dd w),
    \\ \dd r_2|_\lambda &= \dd a+ \dd d -\gamma \dd c,
    \\ \begin{pmatrix}
        \dd r_3|_\lambda & \dd r_4|_\lambda
        \\ \dd r_5|_\lambda & \dd r_6|_\lambda
    \end{pmatrix} &= \scalebox{0.95}{$\begin{pmatrix}
        q_v(1-q_v)\alpha\dd a-q_v\tbinom{q_v}{2}\alpha\gamma \dd c 
        & q_v\alpha\gamma (\dd u-\dd w) -\tbinom{q_v}{2}\alpha\gamma (\dd a+\dd d) -\tbinom{q_v}{3}\alpha\gamma^2\dd c
        \\ \alpha(1-q_v^2) \dd c 
        & \alpha (1-q_v) \dd d-\tbinom{q_v}{2}\alpha\gamma\dd c
    \end{pmatrix}$}, 
    \\ \dd r_7|_\lambda &= \dd a+\dd d.
\end{align*}
    Therefore, we find that the lattice $\Lambda^{\square}$ defined as in the proof of the previous lemma has basis $(\dd r_1|_{\lambda},\dots,\dd r_5|_{\lambda})$, while the lattice $\Lambda^{\st}$ corresponding to the ideal $\mathcal{I}^{\st}_v$ has basis $ (\dd r_1|_{\lambda},\dots,\dd r_4|_{\lambda},\dd r_7|_{\lambda})$. This shows that $\ell_{\mathcal{O}}(\Lambda^{\st}/\Lambda^{\square})=\ord_{\varpi}(q_v^2-1)$, which gives the desired result.
\end{proof}

\begin{lemma}\label{cotangentuni}
Let $\lambda\colon R^{\uni}_v\to \mathcal{O}$ be an augmentation such that the corresponding representation $\rho_\lambda\colon G_v\to \GL_2(\mathcal{O})$ is unramified with $\mathrm{char}(\rho_\lambda(\Frob_v),X)=(X-\alpha)(X-\beta)$, $\beta\equiv q_v\alpha \bmod \varpi$, and $\beta\alpha^{-1}\not\in \{1,q_v,q_v(q_v-2)^{-2}\}$. We set $\eta$ to be the square-root of $q_v^{-1}\alpha\beta$ that satisfies $\eta\equiv 1 \bmod \varpi$.
\begin{enumerate}
    \item If $q_v\not\equiv 1 \bmod p$, or if $q_v\equiv 1 \bmod p $ and $v$ is not a trivial place for  $\overline{\rho}$, then  
    $$\ell_{\mathcal{O}}(\Phi_{\lambda,R^{\uni}_v/\overline{R}^{\unr}_v})=\ord_{\varpi}(\alpha-\eta).$$
    \item if $q_v\equiv 1 \bmod p$ and $v$ is a trivial place for $\overline{\rho}$, then
    $$ \delta_\lambda(R^{\uni}_v)=2\cdot\min\big(\ord_\varpi(\alpha-\eta),\ord_\varpi(\beta-\eta)\big),$$
    and
    $$ \ell_{\mathcal{O}}(\Phi_{\lambda,R^{\uni}_v/\overline{R}^{\unr}_v})=2\cdot \min(\ord_{\varpi}(\alpha-\eta),\ord_{\varpi}(\beta-\eta))+\ord_{\varpi}(\alpha-\eta). $$
\end{enumerate}
\end{lemma}
\begin{proof} $ $

     $\bullet$ If $q_v\not \equiv 1 \bmod p$, the ring $R^{\uni}_v$ is the quotient of $\mathcal{O}[[x,X,Y,a,b,c]]$ by the relations given by the entries of
     $$ N^2, \ A(1+N)-(1+q_vN)A, \ N(A-\eta(1+x)),\ (A-q_v\eta(1+x)^{-1})N, $$
     with
     \begin{equation*}
         A=\begin{pmatrix}
         q_v\eta(1+x)^{-1} & 0 \\ 0 & \eta(1+x)
     \end{pmatrix}, \quad \quad N=\begin{pmatrix}
         a & b \\ c & -a
     \end{pmatrix},\quad \quad P=\begin{pmatrix}
         1 & X \\ Y & 1
     \end{pmatrix}.
     \end{equation*}
     The universal representation is given by $\rho^{\uni}(\Frob_v)=P^{-1}AP$ and $\rho^{\uni}(\sigma_v)=1+P^{-1}NP$ (see \cite[Lemma 5.4]{Sho16}). The second equation forces $a=0$, and combining the second and third equations gives $c=0$. More generally, we obtain 
     $$ R^{\uni}_v\cong \mathcal{O}[[x,X,Y,b]]/(bx).$$
     From this we see that $\overline{R}^{\unr}_v\cong \mathcal{O}[[x,X,Y]]$.
     With respect to the augmentation $\lambda\colon X,Y,b\mapsto 0, x\mapsto \eta^{-1}\alpha-1$, we find
\begin{equation*}
\Phi_{\lambda,R^{\uni}_v/\overline{R}^{\unr}_v}= (\alpha-\eta)\cdot\mathcal{O}.
\end{equation*}  

$\bullet$ Now suppose that $q_v\equiv 1 \bmod p$, and let $\gamma\in \mathcal{O}$. Consider the ring $\mathcal{R}=\mathcal{O}[[x,y,t,z,a,b,c]]$, and the ideal $\mathcal{I}_{v}^{\uni}$ generated by the entries of the matrices
    $$ N^2,\ N(A-(1+x)),\ (A-\alpha\beta(1+x)^{-1})N,\ AN-q_vNA,\ \det(A)-\alpha\beta, $$
    with
    \begin{equation*}
        A=\begin{pmatrix}
        \alpha\beta(1+x)^{-1}+y & t+\gamma \\ z & 1+x-y
    \end{pmatrix},\quad \quad   N=\begin{pmatrix}
        a & b \\ c & -a
    \end{pmatrix}.
    \end{equation*}
    By \cite[\S 5.1 Case ($\varphi$-uni)]{BKM2} (see also \cite[Lemma 2.4]{Calegari18}), we have that $R^{\uni}_v=\mathcal{R}/\mathcal{I}^{\uni}_v$ with the universal representation given by $\rho_v^{\text{uni}}(\Frob_v)=A$, and $\rho_v^{\text{uni}}(\sigma_v)=N+\id$. By the same argument as in the proof of \cite[Lemma 5.2]{BKM2}, $\mathcal{I}^{\uni}_v$ is generated by 
        \begin{gather*}
        r_1=a(x+1-\eta), \quad r_2=b(x+1-\eta), \quad r_3=c(x+1-\eta),
        \\r_4=\alpha\beta y+(y^2+z(t+\gamma))(1+x)-y(1+x)^2,\quad r_5=a^2+bc,\quad  r_6=\eta az- c(\alpha\beta-\eta^2+\eta y), 
        \\ r_7=ay+c(t+\gamma),\quad 
        r_8=a(\alpha\beta-\eta^2+\eta y)+\eta bz, \quad r_9=by-a(t+\gamma).
    \end{gather*}
We have $\overline{R}^{\unr}_v=\mathcal{R}/\mathcal{I}^{\unr}$, with $\mathcal{I}^{\unr}_v=(a,b,c,r_4)$. The augmentation $\lambda\colon R^{\uni}_v\to \mathcal{O}$ that we consider is given by 
    $$ \lambda(z)=\lambda(a)=\lambda(b)=\lambda(c)=\lambda(t)=\lambda(y)=0, \quad \lambda(x)=\alpha-1. $$
    
    Note that if $\gamma\in \mathcal{O}^\times$, then by straightforward calculations we have that  $\mathcal{I}_v^{\uni}=(r_2,r_7,r_9,r_4)$. From $r_7$ and $r_9$, we see that the variables $a$ and $c$ are determined by $b,y,$ and $t$. Therefore, the ring $R^{\uni}_v=\mathcal{O}[[x,y,t,z,b]]/(r_2,r_4)$ is a complete intersection, and an easy calculation gives that
 \begin{equation*}
\Phi_{\lambda,R^{\uni}_v/\overline{R}^{\unr}_v}= (\alpha-\eta)\cdot\mathcal{O}.
 \end{equation*}
    From now on, we assume $\gamma=0$. To compute the Wiles defect, we set
    \begin{gather*}
        s_1=r_3, \  \quad \ s_2=\eta \cdot r_1+r_8-\eta\cdot r_7, \ \quad \ s_3=\eta\cdot r_5+r_2, \ \quad \ s_4=\eta^2\cdot r_2-r_4.
        \\ \theta_1=t-z,\ \quad \ \theta_2=\eta b-z, \ \quad \ \theta_3=\eta c-x+\alpha-1.
    \end{gather*}
    We write $\widetilde{r}_i,\widetilde{s}_j,\widetilde{\theta}_k$ for the analogues of $r_i,s_j,\theta_k$ obtained by specializing $\alpha,\beta,\eta$ to $1$. Define $\widetilde{R}_{\bbZ}\colonequals\bbZ[x,y,t,z,a,b,c]/(\widetilde{\theta}_1,\widetilde{\theta}_2,\widetilde{\theta}_3,\widetilde{s}_1,\dots,\widetilde{s}_4)$. In \Cref{appendixA}, we showed that $\widetilde{R}_{\bbZ}\otimes K$ 
    is free of rank $16$ over $K$ for any field $K$, and computed an explicit basis. It follows that $\widetilde{R}_{\bbZ}$ is is itself free of rank $16$, with basis 
    \begin{gather*}
               b_{1}= 1,\ b_{2}= y,\ b_{3} = a, \ b_{4}=b,\ b_{5}= c,\ b_{6}= cy,  \ b_7 = ya,\ b_8 = yab,
       \\ b_9= yac,\ b_{10}= ybc,\ b_{11}= yb,\ b_{12}= ab,\ b_{13}= abc,\ b_{14}= ac,\ b_{15}=  bc, \ b_{16}= yabc.
    \end{gather*}
    Its socle modulo any prime is generated by $b_{16}$.
    
    Similarly, it follows from \Cref{appendixA} that the ring $R^{\uni}_{\bbZ}\colonequals\bbZ[x,y,t,z,a,b,c]/(\widetilde{\theta}_1,\widetilde{\theta}_2,\widetilde{\theta}_3,\widetilde{r}_1,\cdots,\widetilde{r}_9)$ is free over $\mathbb{Z}$ of rank $6$ with basis $b_1,\dots,b_6$. Its socle modulo any prime is generated by $b_6$.
    
    Since $\mathcal{R}/(\varpi,\theta_1,\theta_2,\theta_3,s_1,\dots,s_4)\cong \widetilde{R}_{\bbZ}\otimes k\cong k^{16}$, we see that $(\varpi,\theta_1,\theta_2,\theta_3,s_1,\dots,s_4)$ is a regular sequence. In particular, the ring $\widetilde{R}_v\colonequals \mathcal{R}/(s_1,\dots,s_4)$ is a complete intersection. The evaluation at $\lambda$ of the ideal generated by the $4\times 4$-minors of the Jacobian matrix $(\partial s_i/\partial x, \partial s_i/\partial y,\dots,\partial s_i/\partial c)_{1\le i \le 4}$ gives the ideal
    $$ \big(\alpha(\alpha-\beta)(\alpha-\eta)^2(\alpha\beta+\alpha\eta-2\eta^2)\big)\subseteq \mathcal{O}. $$
    It is non-zero by the assumption $\beta\alpha^{-1}\not\in \{1,q_v,q_v(q_v-2)^{-2}\}$, and so $\widetilde{R}_v$ is generically formally smooth at $\lambda$.

    We now reduce to working with the rings $$\widetilde{R}_{v,\theta}=\widetilde{R}_v/(\theta_1,\theta_2,\theta_3), \ \quad \ R^{\uni}_{v,\theta}=R^{\uni}_v/(\theta_1,\theta_2,\theta_3),$$
    and denote by $I_\theta$ the kernel of the projection map $\pi_\theta\colon \widetilde{R}_{v,\theta}\to R^{\uni}_{v,\theta}$. Since $\varpi$ is a regular element in $\widetilde{R}_{v,\theta}$, the $b_i$'s form an $\mathcal{O}$-basis of $\widetilde{R}_{v,\theta}$, and $I_{\theta}$ is generated as an $\mathcal{O}$-module by $b_7,\dots,b_{16}$. Let $(b_i^*)_{1\le i\le 16}$ denote the dual basis in $\Hom_{\mathcal{O}}(\widetilde{R}_{v,\theta},\mathcal{O})$. By \cite[Proposition 5.39]{BKM2}, $b_{16}^*$ is a generator of $\Hom_{\mathcal{O}}(\widetilde{R}_{v,\theta},\mathcal{O})$ as a free $\widetilde{R}_{v,\theta}$-module, and $b_6^*$ is a generator of $\Hom_{\mathcal{O}}(R^{\uni}_{v,\theta},\mathcal{O})$ as a free $R^{\uni}_{v,\theta}$-module. We consider the isomorphism 
    \begin{align*}
        \Theta\colon \widetilde{R}_{v,\theta} &\xrightarrow{\sim} \Hom_{\mathcal{O}}(\widetilde{R}_{v,\theta},\mathcal{O})
        \\ r&\mapsto (r'\mapsto b_{16}^*(rr')).
    \end{align*}
    Then as in \cite[p.2510]{BKM2}, we have an isomorphism $ \Hom_{\mathcal{O}}(R^{\uni}_{v,\theta},\mathcal{O})\cong \widetilde{R}_{v,\theta}[I_{\theta}], b_6^*\mapsto \Theta^{-1}(b_6^*)$. Let us write $\Theta^{-1}(b_6^*)=\sum_i \mu_i b_i$, and $b_ib_j=\sum_{k}c_{ijk}b_k$ for some $\mu_i,c_{ijk}\in \mathcal{O}$. Then by definition of $\Theta$, we have 
    $$ b_6^*(b_j)=b_{16}^*(\sum_i \mu_i b_ib_j)= \sum_{i}\mu_i c_{ij16}. $$
Considering the matrix $C=(c_{ij16})_{i,j}$, we see that the row vector $(\mu_i)_i$ is the $6$th row vector of the matrix $C^{-1}$. Using Macaulay2 \cite{M2}, we are able to compute $C$. 
But in order to make neater computations, we reduce to the case $\eta=1$ by making the following change of variables 
$$ x\mapsto \eta x'+\eta-1,\quad y,z,t\mapsto \eta y',\eta z',\eta t',\quad a,b,c\mapsto a',b',c',\quad \alpha=\eta\alpha',\beta= \eta \beta'. $$
We find that $\det C=1$, and that 
$$\Theta^{-1}(b^*_6)= (1-\alpha')(\alpha'+\alpha'\beta'-2)b_1+ (-2\alpha'\beta'-\alpha'+3)b_5+b_{12}.$$ Evaluating at $\lambda$, we get that $\lambda(\widetilde{R}_{v,\theta}[I_\theta])=(\alpha-\eta)(\alpha\eta+\alpha\beta-2\eta^2)\cdot\mathcal{O}$.

Let us write $\Lambda^{\uni}$ (resp. $\widetilde{\Lambda}$) for the lattice generated by $\{\dd r_i|_{\lambda}\}_{1\le i\le 9}$ (resp. $\{\dd s_i\}_{1\le i\le 4}$) inside the free $\mathcal{O}$-lattice $\Lambda^{\unr}=\mathcal{O}\cdot \dd a \oplus \mathcal{O}\cdot \dd b \oplus \mathcal{O}\cdot \dd c \oplus \mathcal{O}\cdot \dd r_4$. The matrices of the generators of $\Lambda^{\uni}$ and $\widetilde{\Lambda}$ in the basis of $\Lambda^{\unr}$ are given by
\begin{align*}
    \begin{pmatrix}
        \alpha-\eta & 0 & 0 & 0 & 0 & 0 & 0 & (\alpha\beta-\eta^2) & 0
        \\ 0 & \alpha-\eta & 0 & 0 & 0 & 0 & 0 & 0 & 0
        \\ 0 & 0 & \alpha-\eta & 0 & 0 & (\alpha\beta-\eta^2) & 0 & 0 & 0
        \\ 0 & 0 & 0 & 1 & 0 & 0 & 0 & 0 & 0
    \end{pmatrix}
\end{align*}
and,
\begin{align*}
    \begin{pmatrix}
        0 & \alpha\eta+\alpha\beta-2\eta^2 & 0 & 0
        \\ 0 & 0 & \alpha-\eta & \eta^2(\alpha-\eta)
        \\ \alpha-\eta & 0 & 0 & 0 
        \\ 0 & 0 & 0 & -1
    \end{pmatrix}.
\end{align*}
Since we have  $ \beta-\eta=\alpha^{-1}(\alpha\beta-\eta^2)-\alpha^{-1}\eta(\alpha-\eta)$, we get that 
\begin{align*}
    \ell_{\mathcal{O}}(\Lambda^{\unr}/\Lambda^{\uni}) &=2 \min(\ord_{\varpi}(\alpha-\eta),\ord_{\varpi}(\beta-\eta))+\ord_{\varpi}(\alpha-\eta),
    \\ \ell_{\mathcal{O}}(\Lambda^{\unr}/\widetilde{\Lambda})&= 2 \ord_{\varpi}(\alpha-\eta)+\ord_{\varpi}( \alpha\eta+\alpha\beta-2\eta^2).
\end{align*}
It follows that
\begin{align*}
    \ell_{\mathcal{O}}(\Phi_{\lambda,\widetilde{R}_{v,\theta}/R^{\uni}_{v,\theta}})&=\ell_{\mathcal{O}}(\Lambda^{\uni}/\widetilde{\Lambda})
    \\ &=\ord_{\varpi}(\alpha-\eta)+ \ord_{\varpi}( \alpha\eta+\alpha\beta-2\eta^2) -2\min(\ord_{\varpi}(\alpha-\eta),\ord_{\varpi}(\beta-\eta)). \qedhere
\end{align*} 
\end{proof}
\begin{remark}
    If the representation $\rho_{\lambda}$ of \Cref{cotangentuni} arises from an abelian variety of $\GL_2$-type (in particular, from a weight-two Hilbert modular form that is Steinberg at some finite place), then the condition $\beta\alpha^{-1}\not\in \{q_v,q_v(2-q_v)^{-2}\}$ is satisfied since $\iota^{-1}(\alpha)$ and $\iota^{-1}(\beta)$ are Weil numbers of the same weight, and so is the condition $\beta\neq \alpha$ (see \cite[Theorem 2.1]{CE98}).  
\end{remark}
\subsection{Global deformation ring and $R=\mathbb{T}$}\label{globaldefringandr=t}
We now define the relevant global deformation rings and record the consequences of the Taylor–Wiles–Kisin patching method, as worked out in \cite{BKM2}. This gives isomorphisms between these rings and the relevant Hecke algebras, and a presentation in terms of the local deformation rings.

We define the functor $\mathcal{D}\colon \CNL_{\mathcal{O}}\to \Sets$ (resp. $\mathcal{D}^\Sigma\colon \CNL_{\mathcal{O}}\to \Sets$) of deformations (resp. $\Sigma$-framed liftings) of $\overline{\rho}$, which associates to every $A\in\CNL_{\mathcal{O}}$ the set of strict equivalence classes of liftings $\rho$ (resp. tuples $(\rho,\{t_v\}_{v\in \Sigma})$), where
\begin{itemize}
    \item $\rho\colon G_F\to \GL_2(A)$ is a lift of $\overline{\rho}$ which is unramified outside of $\Sigma$, finite flat at every place $v\mid p$ and such that $\det \rho= \epsilon_p\cdot \psi_{0,\iota}\circ \Art_F^{-1}$.
    \item $t_v$ is an element of $\ker(\GL_2(A)\to \GL_2(k))$.
\end{itemize}
Here, two liftings $\rho$ and $\rho'$ (resp. two $\Sigma$-framed liftings $(\rho,\{t_v\}_{v\in \Sigma})$ and $(\rho',\{t_v'\}_{v\in \Sigma})$) are said to be strictly equivalent if there exists an element $a\in \ker(\GL_2(A)\to \GL_2(k))$ such that $\rho'=a\rho a^{-1}$ (resp. $\rho'=a\rho a^{-1}$ and $t'_v=a t_v$ for each $v\in \Sigma$). 

It is well known that the functors $\mathcal{D}$ and $\mathcal{D}^\Sigma$ are represented by an element of $\CNL_{\mathcal{O}}$ that we denote by $R$ and $R^{\Sigma}$ (see \cite[(3.2.1)]{Kis09}). Moreover, there is a non-canonical isomorphism $R^\Sigma\cong R[[x_1,\dots,x_{4|\Sigma|-1}]]$, and so we can see $R$ as a quotient of $R^\Sigma$.

Now consider the ring
\begin{align*}
    R_{\text{loc}}^\square&= \Big(\widehat{\bigotimes}_{v\in \Sigma_p} R_v^{\text{fl}}\Big)\widehat{\otimes} \Big(\widehat{\bigotimes}_{v\in \Sigma\setminus \Sigma_p}R_v^{\square}\Big).
\end{align*}
The natural transformation $(\rho,\{t_v\}_{v\in \Sigma})\mapsto (t_v^{-1}\rho|_{G_v}t_v)_{v\in \Sigma}$ induces a morphism $R_{\text{loc}}^\square\to R^{\Sigma}$.

Let us also consider the following rings
\begin{align*}
    R_{\text{loc}}&= \Big(\widehat{\bigotimes}_{v\in \Sigma_p} R_v^{\text{fl}}\Big)\widehat{\otimes} \Big(\widehat{\bigotimes}_{v\in \Sigma_{\min}}R_v^{\min}\Big)\widehat{\otimes} \Big(\widehat{\bigotimes}_{v\in \Sigma_{\st}}R_v^{\st}\Big)\widehat{\otimes} \Big(\widehat{\bigotimes}_{v\in \Sigma_{\uni}}R_v^{\uni}\Big)\widehat{\otimes} \Big(\widehat{\bigotimes}_{v\in \Sigma_\square \cup \Sigma_\triangle}R_v^{\square}\Big),
     \\  R_{\text{loc}}^{\min} &= \Big(\widehat{\bigotimes}_{v\in \Sigma_p} R_v^{\text{fl}}\Big) \widehat{\otimes}\Big(\widehat{\bigotimes}_{v\in \Sigma_{\st}}R^{\st}_v\Big)\widehat{\otimes} \Big(\widehat{\bigotimes}_{v\in \Sigma\setminus \Sigma_p\cup \Sigma_{\st}}R_v^{\min}  \Big),
     \\  \overline{R}_{\text{loc}}^{\min} &= \Big(\widehat{\bigotimes}_{v\in \Sigma_p} R_v^{\text{fl}}\Big) \widehat{\otimes}\Big(\widehat{\bigotimes}_{v\in \Sigma_{\st}}R^{\st}_v\Big)\widehat{\otimes} \Big(\widehat{\bigotimes}_{v\in \Sigma\setminus \Sigma_p\cup \Sigma_{\st}}\overline{R}_v^{\min}  \Big),
\end{align*}
which depend on $Q$ since $\Sigma_{\st}=Q$. These allow us to define the global deformation rings with prescribed local properties: 
$$R^{Q}\colonequals R \widehat{\otimes}_{R^\square_\text{loc}} R_\text{loc}\quad  \quad \text{and,}\quad  \quad R^{Q,\min}\colonequals R \widehat{\otimes}_{R^\square_\text{loc}} R_\text{loc}^{\min}.$$
By \cite[Lemma 2.4]{Man21} and \cite[Lemmma 6.3]{BKM2}, the Galois representation in \Cref{Tgaloisrep} induces surjective maps 
$$R^{Q}\twoheadrightarrow \overline{\mathbb{T}}{}^{Q}_{\mathcal{O},\psi_0}(U^{\mathbf{f}})_{\mathfrak{m}} \quad \text{   and,   } \quad R^{Q,\min}\twoheadrightarrow \mathbb{T}^{Q}_{\mathcal{O},\psi_0}(U^{\min})_{\mathfrak{m}},$$
such that for $v\in \Sigma_{\uni}$, the induced map $R^{\uni}_v\to \overline{\mathbb{T}}{}^{Q}_{\mathcal{O},\psi_0}(U^{\mathbf{f}})_{\mathfrak{m}}$ sends $a_v^{\uni}$ to $\mathbf{U}_v$, where $(\rho^{\uni}_v,a_v^{\uni})$ is the universal pair for $R^{\uni}_v$.

The Taylor--Wiles--Kisin patching method in the proof of \cite[Theorem 6.4]{BKM2} and \cite[Theorem 6.3]{BKM1} produces a diagram
\begin{equation}
            \begin{tikzcd}
        & R_\infty = R_{\text{loc}}[[x_1,\dots x_g]] \arrow[r,twoheadrightarrow] \arrow[d,twoheadrightarrow] &  R^Q 
        \\ S_\infty=\mathcal{O}[[y_1,\dots,y_d]] \arrow[ur,hookrightarrow] \arrow[r,hookrightarrow] \arrow[dr,hookrightarrow] & \overline{R}_\infty^{\min}=\overline{R}^{\min}_{\text{loc}}[[x_1,\dots,x_g]]
        \\ & R_\infty^{\min} = R^{\min}_{\text{loc}}[[x_1,\dots,x_g]]  \arrow[u] \arrow[r,twoheadrightarrow] & R^{Q,\min}
    \end{tikzcd}
    \label{PatchingDiagram}
\end{equation}
with $R^Q=R_\infty/(y_1,\dots,y_d)$, and $R^{Q,\min}=R_\infty^{\min}/(y_1,\dots,y_d)$. A key consequence of the patching construction is that
 $(y_1,\dots,y_d,\varpi)$ forms a regular sequence in $R_{\infty}$ and $R_{\infty}^{\min}$. The method also yields the following result.
\begin{theorem}\label{R=T}
The surjections
\begin{equation*}
    R^{Q}\twoheadrightarrow \overline{\mathbb{T}}{}^{Q}_{\mathcal{O},\psi_0}(U^{\mathbf{f}})_{\mathfrak{m}} \quad \quad \text{and} \quad \quad R^{Q,\min}\twoheadrightarrow \mathbb{T}^{Q}_{\mathcal{O},\psi_0}(U^{\min})_{\mathfrak{m}}
\end{equation*}
are in fact isomorphisms. Moreover, if $Q$ does not contain any trivial place, then $\overline{\mathbb{T}}{}^{Q}_{\mathcal{O},\psi_0}(U^{\mathbf{f}})_{\mathfrak{m}}$ 
is Gorenstein, $\mathbb{T}^{Q}_{\mathcal{O},\psi_0}(U^{\min})_{\mathfrak{m}}$ is a complete intersection, and $M_Q^{\psi_0}(U^{\min})$ is free of rank one over $\mathbb{T}^{Q}_{\mathcal{O},\psi_0}(U^{\min})_{\mathfrak{m}}$.
\end{theorem}
\begin{proof}
    The first statement is \cite[Theorem 6.4]{BKM2}. For the second and third statements, we use the presentation of \eqref{PatchingDiagram}, \cite[Lemma 4.4 (3-4)]{BKM1}, and \Cref{propertiesoflocaldefrings}.
\end{proof}
\subsection{Freeness of some Hecke modules}\label{freenessofM}
 We turn to proving the freeness of $M_Q(U)$ over the corresponding Hecke algebra, assuming $Q$ contains no trivial places. The argument follows \cite[\S 3.2]{Diamond97} closely, using the augmentation associated to the form $\mathbf{g}$ of \Cref{formg} (of nebentypus $\psi_0$) to level raise from the minimal level.

For a place $v\in \Sigma_{\square}\cup \Sigma_{\uni}$, let $(A_v,B_v)$ be the roots of the polynomial $X^2-\lambda_{\mathbf{g}}(T_v)X+\psi_{0,\iota}(\varpi_v)q_v$ such that $B_v\equiv q_vA_v \bmod \varpi$ (we need to make a choice if $q_v\equiv \pm1 \bmod p$). For $v\in \Sigma_{\triangle}$, we note that $\mathbf{g}$ is Steinberg at $v$.

We define a map $\gamma_v\colon M^{\psi_0}_Q(U^{\min})\rightarrow M^{\psi_0}_Q(U^{\mathbf{f}})$ given for $v\in \Sigma_{\uni}$ by
$$\gamma_{v}\colonequals A_v\big[\Gamma_0(v) \GL_2(\mathcal{O}_{F_v})\big]-\Big[ \Gamma_0(v) \begin{pmatrix}
    1 & 0 \\ 0 & \varpi_v
\end{pmatrix} 
\GL_2(\mathcal{O}_{F_v})\Big],$$
for $v\in \Sigma_{\triangle}$ by 
$$ \gamma_v\colonequals \psi_{0,\iota}(\varpi_v)q_v \big[ \Gamma_0(v^2)\Gamma_0(v)\big]- \Big[ \Gamma_0(v^2)\begin{pmatrix}
    1 & 0 \\ 0 & \varpi_v
\end{pmatrix}\Gamma_0(v) \Big]\mathbf{U}_v,  $$
and for $v\in \Sigma_{\square}$ by 
$$ 
\gamma_v \colonequals
\psi_{0,\iota}(\varpi_v) q_v
\big[\Gamma_0(v^2)\GL_2(\mathcal{O}_{F_v})\big]
-
\Big[
\Gamma_0(v^2)
\begin{pmatrix}
1 & 0\\
0 & \varpi_v
\end{pmatrix}
\GL_2(\mathcal{O}_{F_v})
\Big] T_v
+
\Big[
\Gamma_0(v^2)
\begin{pmatrix}
1 & 0\\
0 & \varpi_v^2
\end{pmatrix}
\GL_2(\mathcal{O}_{F_v})
\Big]. $$
The level-raising map is then the product 
$ \gamma\colonequals \prod_{v\in \Sigma_{\uni}\cup \Sigma_\square \cup \Sigma_{\triangle}}\gamma_v \colon M^{\psi_0}_Q(U^{\text{min}})\rightarrow M^{\psi_0}_Q(U^{\mathbf{f}}).$ 

There is a unique oldform $\mathbf{g'}$ of level $\mathfrak{n}$ and nebentypus $\psi_0$ characterized by the fact that $T_v\mathbf{g'}=\lambda_{\mathbf{g}}(T_v)\mathbf{g'}$ for $v\nmid p\mathfrak{n}$, $\mathbf{U}_v\mathbf{g'} =0$ for $v\in \Sigma_\triangle$, $\mathbf{U}_v\mathbf{g'}=0$ for $v\in \Sigma_{\square}$, and $\mathbf{U}_v \mathbf{g'}=A_v \mathbf{g'}$ for $v\in \Sigma_{\uni}$. 

The Hecke-equivariant perfect pairing of \Cref{setting} descends to the $\psi_{0,\iota}$-part, giving a perfect pairing
$$\langle \cdot,\cdot\rangle_U\colon M_Q^{\psi_0}(U)\times M_Q^{\psi_0}(U) \to \mathcal{O}, \quad U^{\min}\subseteq U\subseteq U^{\mathbf{f}}$$ 
(see \cite[Lemma 4.11]{Man21}). Similarly to \cite[Proposition 1.4]{DFG04}, we have the following lemma.
\begin{lemma}\label{levelraising} Write $M^{\psi_0}_Q(U^{\min})[\lambda_\mathbf{g}]=\mathcal{O}\cdot x$. Then
$$ M^{\psi_0}_Q(U^{\mathbf{f}})[\lambda_\mathbf{g'}]=\mathcal{O}\cdot \gamma(x).$$
Moreover, let $\gamma^t$ be the adjoint of $\gamma$ with respect to the Hecke-equivariant pairing on both modules, then
\begin{align*}
            \gamma^t\circ\gamma(x) &= \prod_{v\in \Sigma_{\square}} \psi_{0,\iota}(\varpi_v)q_v(1-q_v)(\lambda_{\mathbf{g}}(T_v)^2-(1+q_v)^2\psi_{0,\iota}(\varpi_v))\prod_{v\in \Sigma_{\uni}}(A_v-B_v)(A_v^2-\psi_{0,\iota}(\varpi_v))
            \\ &\times \prod_{v\in \Sigma_{\triangle}} \psi_{0,\iota}(\varpi_v)^2q_v(1-q_v^2)\lambda_{\mathbf{g}}(\mathbf{U}_v)\cdot x.
\end{align*}
    \end{lemma}
    \begin{proof}
         For the first assertion, we need to show that $\gamma$ is injective with $\mathcal{O}$-torsion free cokernel. We therefore verify that $\gamma_v \bmod \varpi$ is injective for all $v$. For $v\in \Sigma_{\uni}$, this follows from the injectivity of 
         \begin{align*}
         (M_Q^{\psi_0}(U^{\min})\otimes_{\mathcal{O}} k)^2 &\longrightarrow M_Q^{\psi_0}(U^{\min}_0(v))\otimes_{\mathcal{O}} k
         \\ (f_1,f_2) &\mapsto \big[\Gamma_0(v) \GL_2(\mathcal{O}_{F_v})\big] f_1 +\Big[ \Gamma_0(v) \begin{pmatrix}
    1 & 0 \\ 0 & \varpi_v
\end{pmatrix} 
\GL_2(\mathcal{O}_{F_v})\Big] f_2,
         \end{align*}
         which is (essentially) the statement of Ihara's lemma, proved in the definite case in \cite[Proof of Lemma 3.1]{Tay06} and in the indefinite case in \cite{MS21}.
         
         For $v\in \Sigma_\triangle$, it is shown in the proof of \cite[Lemma 3.1]{Tay06} that the kernel of the map $(M_Q^{\psi_0}(U^{\min})\otimes_{\mathcal{O}} k)^2 \to M_Q^{\psi_0}(U^{\min}_0(v^2))\otimes_{\mathcal{O}} k$, defined by the same formula as above, consists of forms unramified at $v$. Since $\overline{\rho}$ is ramified at $v$ by assumption, this kernel is necessarily trivial. Finally, the case $v\in \Sigma_{\square}$ is handled in \cite[Lemma 3.1]{Tay06}.
         
          For the second assertion, we note that for a congruence subgroup $U$, and \emph{diagonal} $g\in G_Q(\mathbb{A}_{F}^{a,\infty})$, we have that 
         $$ \langle [ U_0(v)g U]-,-\rangle_{U_0(v)}=\langle-,[Ug\begin{pmatrix}
             \varpi_v & 0 \\ 0 & 1 
         \end{pmatrix}U_0(v)]-\rangle_U. $$
         Therefore, we get that for $v\in \Sigma_{\uni},$
         $$\gamma_v^t= A_v\left[ \GL_2(\mathcal{O}_{F_v}) \begin{pmatrix}
    \varpi_v & 0 \\ 0 & 1
\end{pmatrix} 
\Gamma_0(v)\right]- \psi_{0,\iota}(\varpi_v)  \left[
\GL_2(\mathcal{O}_{F_v})\Gamma_0(v)\right], $$
and so,
 \begin{align*}
     \gamma_v^t\circ \gamma_v(x)&= \big(A_v^2\lambda_{\mathbf{g}}(T_v)-2A_v\psi_{0,\iota}(\varpi_v)(q_v+1)+\psi_{0,\iota}(\varpi_v)\lambda_{\mathbf{g}}(T_v)\big) (x)
     \\ &= \big(A_v^2(A_v+B_v)-2A_v^2B_v-2A_v\psi_{0,\iota}(\varpi_v)+\psi_{0,\iota}(\varpi_v)(A_v+B_v)\big)\cdot x
    \\ &= (A_v-B_v)(A_v^2-\psi_{0,\iota}(\varpi_v))\cdot x
 \end{align*}
 (see \cite[proof of Lemma 6.4]{CH18} for a similar calculation). For the case $v\in \Sigma_{\triangle}$, we similarly get that
 $$ \gamma_v^t = \psi_{0,\iota}(\varpi_v)q_v\left[\Gamma_0(v) \begin{pmatrix}
    \varpi_v & 0 \\ 0 & 1
\end{pmatrix} \Gamma_0(v^2) \right]- \psi_{0,\iota}(\varpi_v) \mathbf{U}_v [\Gamma_0(v)\Gamma_0(v^2)], $$
 and so 
\begin{align*}
    \gamma_v^t\circ\gamma_v(x)= \psi_{0,\iota}(\varpi_v)\Big(\psi_{0,\iota}(\varpi_v)q_v^2\mathbf{U}_v- \psi_{0,\iota}(\varpi_v)q_v^2\mathbf{U}_v-\psi_{0,\iota}(\varpi_v)q_v^2\mathbf{U}_v+ \mathbf{U}_v \mathbf{V}_v\mathbf{U}_v\Big)(x),
\end{align*} 
where $\mathbf{V}_v$ is the double coset operator $[\Gamma_0(v) \begin{pmatrix}
    1 & 0 \\ 0 & \varpi_v
\end{pmatrix}\Gamma_0(v)]$. Given that $\mathbf{U}_v\mathbf{V}_v(x)=\psi_{0,\iota}(\varpi_v)\cdot x$, we get that 
$$\gamma^t_v\circ \gamma_v(x)= \psi_{0,\iota}(\varpi_v)^2q_v(1-q_v^2)\mathbf{U}_v(x).$$
Finally, the case of $v\in \Sigma_\square$ follows from the same calculations as in \cite[\S 3]{Tay06}.
 \end{proof}
 \begin{theorem}\label{freenesswithouttrivialplaces}
     If $Q$ does not contain any trivial place, then $M_Q(U^{\mathbf{f}})$ is free of rank one over $\overline{\mathbb{T}}{}^{Q}_{\mathcal{O}}(U^{\mathbf{f}})_{\mathfrak{m}}$.
 \end{theorem}
 \begin{proof}
Write $M_Q^{\psi_0}(U^{\min})[\lambda_{\mathbf{g}}]=\mathcal{O}\cdot x$, and set $\eta=\psi_{0,\iota}(\varpi_v)^{1/2}$. By \cite[Lemma 3.5]{BKM1} together with \Cref{levelraising}, we have 
\begin{align*}
    \ell_{\mathcal{O}}&\big(\Psi_{\lambda_{\mathbf{g'}}}(M^{\psi_0}_Q(U^{\mathbf{f}}))\big)= \ord_{\varpi}(\langle \gamma(x),\gamma(x)\rangle_{U^{\mathbf{f}}})
    \\ &= \ord_{\varpi}(\langle x,\gamma^t\circ\gamma(x)\rangle_{U^{\min}})
    \\ &= \sum_{v\in \Sigma_\square} \ord_{\varpi}((1-q_v)(\lambda_{\mathbf{g}}(T_v)^2-(1+q_v)^2\eta^2)
     +\sum_{v\in \Sigma_{\uni}}\ord_{\varpi}((A_v-B_v)(A_v-\eta))
     \\ &+ \sum_{v\in\Sigma_{\triangle}}\ord_{\varpi}(q_v^2-1)+\ord_{\varpi}(\langle x,x\rangle_{U^{\min}})
    \\ &= \sum_{v\in \Sigma_\square} \ord_{\varpi}((1-q_v)(\lambda_{\mathbf{g}}(T_v)^2-(1+q_v)^2\eta^2)
    + \sum_{v\in \Sigma_{\uni}}\ord_{\varpi}((A_v-B_v)(A_v-\eta))
     \\ &+ \sum_{v\in\Sigma_{\triangle}}\ord_{\varpi}(q_v^2-1)+ \ell_{\mathcal{O}}\big(\Psi_{\lambda_{\mathbf{g}}}(M^{\psi_0}_Q(U^{\min}))\big).
\end{align*}
Here, we use that for $v\in \Sigma_{\uni}$, one has $A_v+\eta\in \mathcal{O}^\times$. For $v\in \Sigma_{\triangle}$, we use that $\lambda_{\mathbf{g}}(\mathbf{U}_v)\in \mathcal{O}^\times$ since $\mathbf{g}$ is Steinberg at $v$. 

Using the diagram \eqref{PatchingDiagram}, and arguing as in the proof of \cite[Theorem 7.14]{BKM1}, we get the following diagram with exact rows
\begin{center}
    \begin{tikzcd}
        0 \arrow[r] & \mathcal{O}^d \arrow[r] \arrow[d, equal] & \Phi_{\lambda_{\mathbf{g}'},R_\infty}\arrow[r] \arrow[d,twoheadrightarrow] & \Phi_{\lambda_{\mathbf{g}'},R^Q} \arrow[r] & 0 
        \\ 0 \arrow[r] & \mathcal{O}^d \arrow[d,equal] \arrow[r]& \Phi_{\lambda_{\mathbf{g}'},\overline{R}_\infty^{\min}} 
        \\ 0 \arrow[r] & \mathcal{O}^d \arrow[r] & \Phi_{\lambda_{\mathbf{g}},R^{\min}_\infty}\arrow[u,hookrightarrow]\arrow[r] & \Phi_{\lambda_{\mathbf{g}},R^{Q,\min}} \arrow[r] & 0.
    \end{tikzcd}
\end{center}
A simple diagram chase gives
\begin{align*}
     \ell_{\mathcal{O}}(\Phi_{\lambda_{\mathbf{g}'},R^Q})-\ell_{\mathcal{O}}(\Phi_{\lambda_{\mathbf{g}},R^{Q,\min}}) &= \ell_{\mathcal{O}}(\Phi_{\lambda_{\mathbf{g}'},R_{\infty}/\overline{R}_{\infty}^{\min}})-\ell_{\mathcal{O}}\big(\coker( \Phi_{\lambda_{\mathbf{g}},R^{\min}_\infty}\hookrightarrow \Phi_{\lambda_{\mathbf{g}'},\overline{R}^{\min}_\infty})\big).
\end{align*}
Now we have 
$$\ell_{\mathcal{O}}(\Phi_{\lambda_{\mathbf{g}'},R_{\infty}/\overline{R}_{\infty}^{\min}})=\sum_{v\in \Sigma_\square\cup\Sigma_{\triangle}}\ell_{\mathcal{O}}(\Phi_{\lambda_{\mathbf{g}},R^\square_v/R^{\min}_v})+ \sum_{v\in \Sigma_{\uni}} \ell_{\mathcal{O}}(\Phi_{\lambda_{\mathbf{g}},R^{\uni}_v/R^{\min}_v}),$$ and a direct computation gives \begin{equation*}
\scalebox{0.97}{$\displaystyle \ell_{\mathcal{O}}\big(\coker( \Phi_{\lambda_{\mathbf{g}},R^{\min}_\infty}\hookrightarrow \Phi_{\lambda_{\mathbf{g}'},\overline{R}^{\min}_\infty})\big)=\sum_{v\in \Sigma_{\uni}} \ell_{\mathcal{O}}\big( \coker(\Phi_{\lambda_{\mathbf{g}},R^{\min}_v}\to \Phi_{\lambda_{\mathbf{g}'},\overline{R}^{\min}_v})\big)= \sum_{v\in \Sigma_{\uni}} \ord_{\varpi}(A_v-B_v).
$}
\end{equation*}
On the other hand, since $Q$ contains no trivial places, \Cref{R=T} applies. Combined with \Cref{WDcompleteIntersection}, this yields 
\begin{align*}
    \ell_{\mathcal{O}}\big(\Psi_{\lambda_{\mathbf{g}'}}(\overline{\mathbb{T}}{}^{Q}_{\mathcal{O},\psi_0}(U^{\mathbf{f}})_{\mathfrak{m}})\big)-\ell_{\mathcal{O}}\big(\Psi_{\lambda_{\mathbf{g}}}(M^{\psi_0}_Q(U^{\min}))\big) &= \ell_{\mathcal{O}}(\Phi_{\lambda_{\mathbf{g}'},R^Q})-\ell_{\mathcal{O}}(\Phi_{\lambda_{\mathbf{g}},R^{Q,\min}})- \delta_{\lambda_{\mathbf{g}'}}(R^Q).
\end{align*}
By the discussion surrounding \eqref{PatchingDiagram}, \Cref{propertiesoflocaldefrings}, and the invariance properties \Cref{invunderregseq,invtensorprod} of the Wiles defect, we have
$$\delta_{\lambda_{\mathbf{g}'}}(R^Q)=\delta_{\lambda_{\mathbf{g}'}}(R_{\mathrm{loc}})=\sum_{v\in \Sigma_{\uni}}\delta_{\lambda_{\mathbf{g}'}}(R_v^{\uni}).$$
Applying \Cref{cotangentsquare1}, \Cref{cotangentsquare2}, and \Cref{cotangentuni}, we finally get 
$$ \ell_{\mathcal{O}}\big(\Psi_{\lambda_{\mathbf{g}'}}(\overline{\mathbb{T}}{}^{Q}_{\mathcal{O},\psi_0}(U^{\mathbf{f}})_{\mathfrak{m}})\big)= \ell_{\mathcal{O}}\big(\Psi_{\lambda_{\mathbf{g'}}}(M^{\psi_0}_Q(U^{\mathbf{f}})\big),$$
which by \Cref{freenessoverGorenstein} and \Cref{R=T} implies that $M^{\psi_0}_Q(U^{\mathbf{f}})$ is free of rank one over $\overline{\mathbb{T}}{}^{Q}_{\mathcal{O},\psi_0}(U^{\mathbf{f}})_{\mathfrak{m}}$.

By \Cref{mult1overE},  $M_Q(U^{\mathbf{f}})[\tfrac{1}{\varpi}]$ is free of rank one over $\overline{\mathbb{T}}{}^{Q}_{\mathcal{O}}(U^{\mathbf{f}})_{\mathfrak{m}}[\tfrac{1}{\varpi}]$, and given the inclusion $\mathfrak{J}_{\psi_0}\subset \mathfrak{m}$, we have
\begin{equation*}
    \dim_k M_Q(U^{\mathbf{f}})/\mathfrak{m}M_Q(U^{\mathbf{f}})=\dim_k M_Q^{\psi_0}(U^{\mathbf{f}})/\mathfrak{m}M_Q^{\psi_0}(U^{\mathbf{f}})=1.
\end{equation*}
Therefore, $M_Q(U^{\mathbf{f}})$ is free of rank one over $\overline{\mathbb{T}}{}^{Q}_{\mathcal{O}}(U^{\mathbf{f}})_{\mathfrak{m}}$.
 \end{proof}
We now establish a quantitative level-lowering result for the ring-theoretic congruence module. This will serve as a key input later, providing an upper bound on the length of the cohomological congruence module by \Cref{upperbound}. Since the patching method requires working with a fixed nebentypus, we must first understand how the congruence module varies upon fixing the nebentypus.
\begin{corollary}\label{congmodule+/-nebentypus}
      Let $C_{U^{\mathbf{f}},p}$ be the $p$-Sylow subgroup of $C_{U^{\mathbf{f}}}$, and let $a_{\psi}= \sum_{g\in C_{U^{\mathbf{f}},p}} \psi_{\iota}(g)^{-1}g \in \mathcal{O}[C_{U^{\mathbf{f}},p}]$. Then we have 
    $$ \ell_{\mathcal{O}}\big(\Psi_{\lambda_{\mathbf{f}}}(\overline{\mathbb{T}}{}^{Q}_{\mathcal{O}}(U^{\mathbf{f}})_{\mathfrak{m}})\big)=\ord_{\varpi}(\lambda_{\mathbf{f}}(a_\psi))+\ell_{\mathcal{O}}\big(\Psi_{\lambda_{\mathbf{f}}}(\overline{\mathbb{T}}{}^{Q}_{\mathcal{O},\psi}(U^{\mathbf{f}})_{\mathfrak{m}})\big). $$
\end{corollary}
\begin{proof}
    % Let $t$ be a non-trivial place in $Q$. By \Cref{freenesswithouttrivialplaces}, $M_{\{t\}}(U^{\mathbf{f}})$ is free of rank one over $\overline{\mathbb{T}}^{\{t\}}_{\mathcal{O}}(U^{\mathbf{f}})_{\mathfrak{m}}$. Let $M_{\{t\}}(U^{\mathbf{f}})^{Q\text{-new}}$ denote the $Q$-new submodule, defined as the orthogonal complement, with respect to $\langle \cdot,\cdot\rangle_{U^{\mathbf{f}}}$, of the subspace spanned by the images of the level raising maps at $v\in Q$. Then  $$M_{\{t\}}(U^{\mathbf{f}})^{Q\text{-new}}=M_{\{t\}}(U^{\mathbf{f}})[I],$$ where $I=\ker(\overline{\mathbb{T}}^{\{t\}}_{\mathcal{O}}(U^{\mathbf{f}})_{\mathfrak{m}}\to \overline{\mathbb{T}}^{Q}_{\mathcal{O}}(U^{\mathbf{f}})_{\mathfrak{m}})$, and so is isomorphic to $(\overline{\mathbb{T}}^{Q}_{\mathcal{O}}(U^{\mathbf{f}})_{\mathfrak{m}})^*$. Since $M_{\{t\}}(U^{\mathbf{f}})$ is a free $\mathcal{O}[C_{U^{\mathbf{f}},p}]$-module by \cite[Lemma 4.4]{Man21}, the summand $M_{\{t\}}(U^{\mathbf{f}})^{Q\text{-new}}$ is also free over $\mathcal{O}[C_{U^{\mathbf{f}},p}]$. Since $\mathcal{O}[C_{U^{\mathbf{f}},p}]$ is complete intersection, we deduce that  $\overline{\mathbb{T}}^{Q}_{\mathcal{O}}(U^{\mathbf{f}})_{\mathfrak{m}}$ is free over $\mathcal{O}[C_{U^{\mathbf{f}},p}]$. 
    By the same argument as in \cite[Lemma 4.16]{Man21}, we have an isomorphism $\overline{\mathbb{T}}{}^{Q}_{\mathcal{O}}(U^{\mathbf{f}})_{\mathfrak{m}}\cong \overline{\mathbb{T}}{}^{Q}_{\mathcal{O},\psi_0}(U^{\mathbf{f}})_{\mathfrak{m}}\otimes_{\mathcal{O}} \mathcal{O}[C_{U^{\mathbf{f}},p}]$. In particular, $\overline{\mathbb{T}}{}^{Q}_{\mathcal{O}}(U^{\mathbf{f}})_{\mathfrak{m}}$ is free over $\mathcal{O}[C_{U^{\mathbf{f}},p}]$.
    It follows by the proof of \cite[Lemma 4.11]{Man21} that there is an isomorphism 
     $$ \overline{\mathbb{T}}{}^{Q}_{\mathcal{O},\psi}(U^{\mathbf{f}})_{\mathfrak{m}}= \overline{\mathbb{T}}{}^{Q}_{\mathcal{O}}(U^{\mathbf{f}})_{\mathfrak{m}}/\mathfrak{J}_\psi \xrightarrow{\times a_\psi} \overline{\mathbb{T}}{}^{Q}_{\mathcal{O}}(U^{\mathbf{f}})_{\mathfrak{m}}[\mathfrak{J}_\psi]  $$
     which is clearly $\overline{\mathbb{T}}{}^{Q}_{\mathcal{O}}(U^{\mathbf{f}})_{\mathfrak{m}}$-equivariant. Taking the $\lambda_{\mathbf{f}}$-part, and noting that $\mathfrak{J}_{\psi}\subset \mathfrak{p}_{\lambda_{\mathbf{f}}}$, we get an isomorphism $\overline{\mathbb{T}}{}^{Q}_{\mathcal{O},\psi}(U^{\mathbf{f}})_{\mathfrak{m}}[\mathfrak{p}_{\lambda_{\mathbf{f}}}]\xrightarrow{\times a_{\psi}}\overline{\mathbb{T}}{}^{Q}_{\mathcal{O}}(U^{\mathbf{f}})_{\mathfrak{m}}[\mathfrak{p}_{\lambda_{\mathbf{f}}}]$. By \Cref{congmoduleforR}, this gives the desired result.
\end{proof}
\begin{corollary}\label{changeofcongmoduleforheckealg}
    Let $v\in Q$, set $\widetilde{Q}=Q\setminus\{v\}$, and assume that $\widetilde{Q}$ contains a non-trivial place. Then we have
    $$ \ell_{\mathcal{O}}\big(\Psi_{\lambda_{\mathbf{f}}}(\overline{\mathbb{T}}{}^{\widetilde{Q}}_{\mathcal{O}}(U^{\mathbf{f}})_{\mathfrak{m}})\big)=\ell_{\mathcal{O}}\big(\Psi_{\lambda_{\mathbf{f}}}(\overline{\mathbb{T}}{}^{Q}_{\mathcal{O}}(U^{\mathbf{f}})_{\mathfrak{m}})\big)+c_v. $$
\end{corollary}
\begin{proof}
    Arguing analogously to \cite[Proposition 7.14]{BKM1} and using \Cref{R=T}  (which admits an analogous formulation with $\psi$ in place of $\psi_0$), we have  
     $$\ell_{\mathcal{O}}(\Phi_{\lambda_{\mathbf{f}},\overline{\mathbb{T}}{}^{\widetilde{Q}}_{\mathcal{O},\psi}(U^{\mathbf{f}})_{\mathfrak{m}}})- \ell_{\mathcal{O}}(\Phi_{\lambda_{\mathbf{f}},\overline{\mathbb{T}}{}^{Q}_{\mathcal{O},\psi}(U^{\mathbf{f}})_{\mathfrak{m}}})= \ell_{\mathcal{O}}(\Phi_{\lambda_{\mathbf{f}}, R^{\widetilde{Q}}_\infty/R^{Q}_\infty})= \ell_{\mathcal{O}}(\Phi_{\lambda_{\mathbf{f}}, R^{\uni}_v/R^{\st}_v}). $$
     Here we decorate the ring $R_\infty$ with a with a superscript $Q$ to emphasize its dependence on $Q$.
     
    Let $n_v$ is the largest integer $n$ such that $\rho_{\mathbf{f}}(G_v) \bmod \varpi^n$ is trivial. As noted in \cite[Remark 7.8]{BKM2}, and using \cite[Proposition 7.9 (3.)]{BKM1}, we can compute that $\ell_{\mathcal{O}}(\Phi_{\lambda_{\mathbf{f}}, R^{\uni}_v/R^{\st}_v})=(c_v-n_v)+2n_v=c_v+n_v$. By \cite[Theorem 6.5]{BKM2}, we have
    $$ \delta_{\lambda_{\mathbf{f}}}(\overline{\mathbb{T}}{}^{\widetilde{Q}}_{\mathcal{O},\psi}(U^{\mathbf{f}})_{\mathfrak{m}}) -\delta_{\lambda_{\mathbf{f}}}(\overline{\mathbb{T}}{}^{Q}_{\mathcal{O},\psi}(U^{\mathbf{f}})_{\mathfrak{m}})=\delta_{\lambda_{\mathbf{f}}}(R^{\uni}_v)-\delta_{\lambda_{\mathbf{f}}}(R^{\st}_v)=3n_v-2n_v=n_v.$$
   Combining these two equalities with the definition of the Wiles defect yields
    $$\ell_{\mathcal{O}}(\Psi_{\lambda_{\mathbf{f}}}(\overline{\bbT}{}^{\widetilde{Q}}_{\mathcal{O},\psi}(U^{\mathbf{f}})_{\mathfrak{m}}))- \ell_{\mathcal{O}}(\Psi_{\lambda_{\mathbf{f}}}(\overline{\bbT}{}^{Q}_{\mathcal{O},\psi}(U^{\mathbf{f}})_{\mathfrak{m}}))=c_v.$$
    Finally, \Cref{congmodule+/-nebentypus} gives
    $$ \ell_{\mathcal{O}}\big(\Psi_{\lambda_{\mathbf{f}}}(\overline{\bbT}{}^{\widetilde{Q}}_{\mathcal{O}}(U^{\mathbf{f}})_{\mathfrak{m}})\big)- \ell_{\mathcal{O}}\big(\Psi_{\lambda_{\mathbf{f}}}(\overline{\bbT}{}^{Q}_{\mathcal{O}}(U^{\mathbf{f}})_{\mathfrak{m}})\big)= \ell_{\mathcal{O}}\big(\Psi_{\lambda_{\mathbf{f}}}(\overline{\bbT}{}^{\widetilde{Q}}_{\mathcal{O},\psi}(U^{\mathbf{f}})_{\mathfrak{m}})\big)- \ell_{\mathcal{O}}\big(\Psi_{\lambda_{\mathbf{f}}}(\overline{\bbT}{}^{Q}_{\mathcal{O},\psi}(U^{\mathbf{f}})_{\mathfrak{m}})\big). \qedhere  $$ 
\end{proof}
 \subsection{Non-Gorensteinness of the Hecke algebra}\label{nonGorenstein}
In this section, we propose a definition of a Gorensteiness defect similar to the one given in \cite{WWE21} and \cite{KW08}. We show that it is invariant under quotienting by a regular sequence, and is multiplicative on completed tensor products. These properties are analogous to those in \Cref{invunderregseq} and \Cref{invtensorprod} of the Wiles defect. Consequently, we deduce that the Hecke algebra $\overline{\bbT}{}^{Q}_{\mathcal{O}}(U^{\mathbf{f}})_{\mathfrak{m}}$ is \emph{not} Gorenstein if $Q$ contains trivial places. An immediate consequence is that the module $M_Q(U^{\mathbf{f}})$ is not free over $\overline{\bbT}{}^{Q}_{\mathcal{O}}(U^{\mathbf{f}})_{\mathfrak{m}}$, given that it is self-dual.
\begin{notation}
    Let $R\in \CNL_{\mathcal{O}}$, and $M$ a finitely generated $R$-module. We denote by $g_R(M)=\dim_k(M/\mathfrak{m}_RM)$ the minimal number of generators of $M$.
\end{notation}
\begin{definition}
     Let $R\in \CNL_{\mathcal{O}}$ which is Cohen--Macaulay and flat over $\mathcal{O}$ of relative dimension $d$. Let $\omega_R$ be the dualizing module of $R$. We define the type of $R$ to be the integer
     $$ r(R) \colonequals g_R(\omega_R).$$
\end{definition}
\begin{remark}
     We found this terminology in \cite[3.6.13]{BH98}.
\end{remark}
Note that $r(R)=1$ if and only if $R$ is Gorenstein, and thus $r$ is a measure of the defect of Gorensteinness. Moreover, the type is invariant under quotienting by regular sequences (see \cite[\S 21.3]{Eisen}). From this, we get the following lemma.
\begin{lemma}\label{Gdefectdim0}
   Let $(\varpi,\theta_1,\dots,\theta_d)$ be a regular sequence in $R$, and let $\overline{R}_\theta\colonequals R/(\varpi,\theta_1,\dots,\theta_d)$. Then we have that 
    $$ r(R)=r(\overline{R}_{\theta})= g_{\overline{R}_{\theta}}(\Hom_{k}(\overline{R}_\theta,k)).$$
\end{lemma}
\begin{proof}
    By \cite[Theorem 21.15]{Eisen}, $\Hom_{k}(\overline{R}_\theta,k)$ is the dualizing module of $\overline{R}_\theta$.
\end{proof}
\begin{lemma}\label{prodoftypes}
    Let $R,R'\in \CNL_{\mathcal{O}}$ be Cohen-Macaulay and flat $\mathcal{O}$-algebras of relative dimension $d$ and $d'$. Then we have
    $$ r(R\widehat{\otimes}_{\mathcal{O}}R')=r(R)\cdot r(R').$$
\end{lemma}
\begin{proof}
    First, let us note that if $(\varpi,\theta_1,\dots,\theta_d)$ is a regular sequence in $R$, and $(\varpi,\theta_1',\dots,\theta_{d'}')$ is a regular sequence in $R'$, then $(\varpi,\theta_1\otimes 1,\dots,\theta_d\otimes 1,1\otimes\theta_1',\dots,1\otimes\theta_{d'}')$ is a regular sequence in $R\widehat{\otimes}_{\mathcal{O}}R'$. Indeed, by \cite[Lemma 4.4]{BKM1} $R\widehat{\otimes}_{\mathcal{O}}R'$ is Cohen--Macaulay of relative dimension $d+d'$ over $\mathcal{O}$, and its quotient by this sequence is a finite $k$-algebra, and so we can invoke \cite[Theorem 2.1.2 (c)]{BH98}. By \Cref{Gdefectdim0}, we can reduce to the case where $R$ and $S$ are both zero dimensional. Then the result follows from the isomorphism 
\begin{equation*}
\Hom_{k}(R\otimes_k S,k)\xrightarrow{\sim}
\Hom_k(R,k)\otimes_k\Hom_k(S,k).\qedhere
\end{equation*}
\end{proof}
\begin{corollary}\label{notGorenstein}
    The set $Q$ contains a trivial place if and only if the algebra $\overline{\bbT}{}^{Q}_{\mathcal{O}}(U^{\mathbf{f}})_{\mathfrak{m}}$ is not Gorenstein.
\end{corollary}
\begin{proof} Suppose  $\overline{\bbT}{}^{Q}_{\mathcal{O}}(U^{\mathbf{f}})_{\mathfrak{m}}$ is Gorenstein, we claim that $\overline{\bbT}{}^{Q}_{\mathcal{O},\psi_0}(U^{\mathbf{f}})_{\mathfrak{m}}$ is Gorenstein as well. Since $\mathcal{O}[C_{U^{\mathbf{f}},p}]$ is a complete intersection, \cite[Theorem 21.15]{Eisen} yields an isomorphism
\begin{equation}\label{HeckeGorensteinDuality}
    \overline{\bbT}{}^{Q}_{\mathcal{O}}(U^{\mathbf{f}})_{\mathfrak{m}}\cong \Hom_{\mathcal{O}[C_{U^{\mathbf{f}},p}]}\big(\overline{\bbT}{}^{Q}_{\mathcal{O}}(U^{\mathbf{f}})_{\mathfrak{m}},\mathcal{O}[C_{U^{\mathbf{f}},p}]\big).
\end{equation}
 As noted in the proof of \Cref{congmodule+/-nebentypus}, $\overline{\bbT}{}^{Q}_{\mathcal{O}}(U^{\mathbf{f}})_{\mathfrak{m}}$ is free over $\mathcal{O}[C_{U^{\mathbf{f}},p}]$, so quotienting both sides of \eqref{HeckeGorensteinDuality} by $\mathfrak{J}_{\psi_0}$ gives the Gorensteinness of $\overline{\bbT}{}^{Q}_{\mathcal{O},\psi_0}(U^{\mathbf{f}})_{\mathfrak{m}}$.

    By \Cref{R=T}, $\overline{\bbT}{}^{Q}_{\mathcal{O},\psi_0}(U^{\mathbf{f}})_{\mathfrak{m}}$ is isomorphic to the universal deformation ring $R^Q$. Moreover, the discussion surrounding \eqref{PatchingDiagram} gives $R^Q=R_{\text{loc}}[[x_1,\dots,x_g]]/(y_1,\dots,y_d)$, where $y_1,\dots,y_d$ forms a regular sequence, so that   $r(R^Q)=r(R_{\text{loc}})$. By \Cref{prodoftypes} and \Cref{propertiesoflocaldefrings}, we have 
    $$r(R_{\text{loc}})= \prod_{v\in \Sigma\setminus \Sigma_p} r(R_v^{\tau_v})=\prod_{v\in Q} r(R_v^{\st}),$$
    which is strictly greater than $1$ if and only if $Q$ contains a trivial place. Therefore, $\overline{\bbT}{}^{Q}_{\mathcal{O},\psi_0}(U^{\mathbf{f}})_{\mathfrak{m}}$ is Gorenstein if and only if $Q$ contains no trivial places. 

    Finally, if $\overline{\bbT}{}^{Q}_{\mathcal{O},\psi_0}(U^{\mathbf{f}})_{\mathfrak{m}}$ is Gorenstein, then by \Cref{freenesswithouttrivialplaces}, $M_Q(U^{\mathbf{f}})$ is free of rank one over $\overline{\bbT}{}^{Q}_{\mathcal{O},}(U^{\mathbf{f}})_{\mathfrak{m}}$. Since $M_Q(U^{\mathbf{f}})$ is self-dual, $\overline{\bbT}{}^{Q}_{\mathcal{O}}(U^{\mathbf{f}})_{\mathfrak{m}}$ is Gorenstein.
\end{proof}
\section{Recollection of facts on abelian varieties}\label{factsonAB}
In this section, we develop the algebra necessary for proving our main theorem in \Cref{complambdashimura}. We start by setting the notation and  recalling some useful facts on abelian varieties with semistable reduction. Specifically, we introduce the monodromy pairing on the character group of the toric part of the abelian variety, which plays an essential role in the proof. Moreover, we recall a few results on the $p$-adic uniformization of abelian varieties. 

In the second part of this section, following the approach of \cite{Helm}, we study the abstract situation of a faithful action of a finite flat $\bbZ$-algebra $\bbT$ acting on an abelian variety. The central result we prove here is \Cref{tatemoduletochargrp}, which allows us to relate the congruence module of the character group to the cohomological congruence module. Finally, in the third part, we introduce $\lambda$-Shimura degrees and relate them to congruence modules.

Let $H$ is a smooth commutative algebraic group over $F$. For an integer $n\in \mathbb{N}_{\ge 1}$, we write $H[p^n]$ for the algebraic subgroup of $H$ consisting of $p^n$-torsion elements, and we set $H[p^{\infty}]=\cup_{n\ge 1}H[p^n]$. We can then form the $p$-adic Tate module $T_pH\colonequals \varprojlim_{n\ge 1} H[p^n](\overline{F})$, and the $p$-adic \textit{contravariant} Tate module $$T^\vee_pH\colonequals \Hom(H[p^{\infty}](\overline{F}),\bbQ_p/\bbZ_p(1)).$$ Both $\bbZ_p$-modules equipped with a $G_F$-action, and they are related by a duality pairing
$$ T_p H\times T_p^\vee H\to \bbZ_p(1), $$
which is perfect if $H$ is an abelian variety.

\subsection{Abelian varieties with semistable reduction}\label{ssred} In this section, we let $\ell\neq p$ be a prime, and $K/\mathbb{Q}_{\ell}$ be a finite extension with ring of integers $\mathcal{O}_K$, residue field $k_K$. We let $G_K$ denote the absolute Galois group of $K$, and $I_K\subset G_K$ be its inertia subgroup.

Let $A$ be an abelian variety over $K$ of dimension $g$, and $\mathcal{A}$ be its Néron model over $\mathcal{O}_{K}$. Its special fibre $\mathcal{A}_{k_K}$ is not necessarily connected, and we let $\mathcal{A}^0_{k_K}$ be the connected component containing the identity. We define the component group $\Phi(A)$ as the finite flat group scheme fitting in the following exact sequence
\begin{equation*}
    0 \rightarrow \mathcal{A}^0_{k_K} \rightarrow \mathcal{A}_{k_K} \rightarrow \Phi(A)\rightarrow 0.
\end{equation*}

By Chevalley's structure theorem, there is a short exact sequence of group schemes
\begin{equation*}
    0 \rightarrow \mathcal{C}\rightarrow \mathcal{A}^0_{k_K} \rightarrow \mathcal{B}\rightarrow 0,
\end{equation*}
where $\mathcal{B}$ is an abelian variety over $k_v$, and $\mathcal{C}$ is the largest connected normal affine subgroup of $\mathcal{A}^0_{k_K}$. Furthermore, there is a unique short exact sequence
\begin{align*}
    0 \rightarrow \mathcal{T}(A) \rightarrow \mathcal{C}\rightarrow \mathcal{U}\rightarrow 0,
\end{align*}
where $\mathcal{T}(A)$ is a torus, and $\mathcal{U}$ is unipotent. We call $\mathcal{T}(A)$ the toric part of $A$. 

Let $\overline{k}_K$ be a fixed algebraic closure of $k_K$. We define the character group of $A$ to be
\begin{equation*}
    \mathcal{X}(A)=\Hom_{\overline{k}_K}(\mathcal{T}_{\overline{k}_K}(A),\mathbb{G}_{m,\overline{k}_K}).
\end{equation*}
Note that the formations of $\Phi(A)$, $\mathcal{X}(A)$, and $\mathcal{T}(A)$ are functorial in $A$. 

From now on, we assume that $A$ has semistable reduction, which means that $\mathcal{U}=0$. Let $\hat{A}$ be the dual abelian variety, and consider the $p$-adic Weil pairing
\begin{equation*}
    (\cdot \  ,\cdot)_p\colon T_p(A)\times T_p(\hat{A})\rightarrow \bbZ_p(1).
\end{equation*}
There is a canonical isomorphism $T_p(A)^{I_{K}}\cong T_p(\mathcal{A}_{k_K}^0)$. Letting $T_p(A)^{\text{tor}}\colonequals T_p(\mathcal{T}(A))$, we obtain a filtration of $\bbZ_p$-modules
\begin{equation*}
    T_p(A)^{\text{tor}}\subseteq T_p(A)^{I_{K}}\subseteq T_p(A).
\end{equation*}
By Grothendieck's orthogonality theorem \cite[ Exposé IX, \S 9.2]{SGA7I}, the submodule $T_p(A)^{I_{K}}$ is the orthogonal complement to $T_p(\hat{A})^{\text{tor}}$ with respect to the $p$-adic Weil pairing. It follows that there exists a morphism
$$ \varphi\colon I_{K} \longrightarrow \Hom_{\bbZ_p}\big(T_p(A)/T_p(A)^{I_{K}},T_p(A)^{\text{tor}}\big) $$
given by $g\mapsto (x\mapsto (g-1)x)$. Since the target is a pro-$p$ group, $\varphi$ must factor through the quotient map $t_p\colon I_{K}\twoheadrightarrow \bbZ_p(1)$.
Noting that $T_p(A)^{\text{tor}}\cong \Hom_{\bbZ_p}(\mathcal{X}(A)\otimes_\bbZ \bbZ_p,\bbZ_p)(1)$, and that $T_p(A)/T_p(A)^{I_{K}}\cong \Hom_{\bbZ_p}(T_p(\hat{A})^{\text{tor}},\bbZ_p)(1)\cong \mathcal{X}(\hat{A})\otimes_\bbZ\bbZ_p$, we get using $\varphi$ a pairing 
\begin{equation}\label{monodromypairing}
    u_{A}\colon \mathcal{X}(A)\otimes\bbZ_p\times \mathcal{X}(\hat{A})\otimes \bbZ_p\rightarrow \bbZ_p
\end{equation}
called the monodromy pairing, and a short exact sequence:
\begin{equation*}
    0 \rightarrow \mathcal{X}(\hat{A})\otimes \bbZ_p \xrightarrow{\alpha} \Hom(\mathcal{X}(A)\otimes \bbZ_p,\bbZ_p)\rightarrow \Phi(A)[p^\infty]\rightarrow 0,
\end{equation*}
where $(\alpha(x))(y)=u_{A}(x,y)$.

Let us now recall some results on the $p$-adic uniformization of abelian varieties (see \cite[\S 3]{Em03}). Since $A$ has semistable reduction, it has a uniformization in the category of rigid analytic spaces over $K$. This is the data of an exact sequence
    $$ 0 \to \mathcal{X}(\hat{A}) \to G_A^{\text{an}} \to A^{\text{an}} \to 0 $$
which is $G_{K}$-equivariant, and where $G_A$ is a semiabelian variety over $K$. A morphism $\pi\colon B\to A$ of abelian varieties with semistable reduction induces a morphism $G_B\to G_A$ of algebraic groups which restricts to the morphism $\mathcal{X}(\hat{B})\to \mathcal{X}(\hat{A})$ induced by $\hat{\pi}$. 

Suppose that we have an exact sequence $0\to C \to B \xrightarrow{\pi} A \to 0$. Following the discussion in \cite[\S 3]{Em03}, we have by functoriality of the uniformization a diagram
    \begin{center}
\begin{tikzcd}
            & 0 \arrow[d] & 0 \arrow[d] & 0 \arrow[d]
        \\ 0 \arrow[r] & \overline{\mathcal{X}}(A)\arrow[r] \arrow[d] &  \overline{G}_{\hat{A}}^{\text{an}} \arrow[d] \arrow[r] & \hat{A}^{\text{an}}  \arrow[d] \arrow[r] & 0
        \\ 0\arrow[r] & \mathcal{X}(B) \arrow[d] \arrow[r]  & G_{\hat{B}}^{\text{an}}  \arrow[r] \arrow[d] & \hat{B}^{\text{an}} \arrow[r] \arrow[d] & 0
             \\ 0\arrow[r] & \mathcal{X}(C) \arrow[d] \arrow[r]  & G_{\hat{C}}^{\text{an}}  \arrow[r] \arrow[d] & \hat{C}^{\text{an}}  \arrow[r] \arrow[d] & 0
             \\ & 0 & 0 & 0
 \end{tikzcd}
    \end{center}
which has exact rows and columns (by \cite[Lemma 3.1]{Em03}). Here, the lattice $\overline{\mathcal{X}}(A)$ and the algebraic group $\overline{G}_{\hat{A}} $ are defined by the property of ensuring the exactness of the two leftmost columns. In fact, the top row need not be the rigid analytic uniformization. However by \cite[Lemma 3.2]{Em03}, we have that $G_{\hat{A}}=(\overline{G}_{\hat{A}})^0$, that $\mathcal{X}(A)=\overline{\mathcal{X}}(A)\cap G_{\hat{{A}}}^{\mathrm{an}}$, and that the map $\overline{\mathcal{X}}(A)/\mathcal{X}(A)\to \overline{G}_{\hat{{A}}}/G_{\hat{A}}$ is an isomorphism. We remark that by torsion freeness of $\mathcal{X}(C)$, we have  
\begin{equation}\label{saturation}
\overline{\mathcal{X}}(A)= \mathcal{X}(B) \cap \mathcal{X}(A)\otimes\bbQ.
\end{equation}

\subsection{A result of Helm}\label{resultofHelm} The main goal of this section is to establish \Cref{tatemoduletochargrp} which is the analogue of \cite[Proposition 4.14]{Helm}. This theorem asserts that under favorable conditions, the contravariant Tate module is isomorphic to two copies of the character group. This enables us to work interchangeably between these two objects.

Let $J$ be an abelian variety defined over $F$, and assume that it is equipped with a faithful action of a finite flat $\bbZ$-algebra $\bbT$, where every element of $\bbT$ acts by an endomorphism defined over $F$. We write $\bbT_{\mathcal{O}}\colonequals \bbT\otimes_{\bbZ}\mathcal{O}$, then this algebra admits a factorization $\bbT_{\mathcal{O}}\cong \prod_{\mathfrak{m}} \bbT_{\mathfrak{m}}$ as the product over all of its maximal ideals $\mathfrak{m}$ of its localization at those ideals. We write $T^\vee_{\mathfrak{m}}J\colonequals(T_p^\vee J\otimes_{\bbZ_p}\mathcal{O})_{\mathfrak{m}}$ for the localized contravariant Tate-module. Then there is a $\bbT[G_F]$-equivariant perfect duality 
$$ T_{\mathfrak{m}}^\vee J \times (J[p^\infty]\otimes_{\bbZ_p}\mathcal{O})_{\mathfrak{m}}\longrightarrow E/\mathcal{O}(1) $$
since this is true before localizing, and the localization at $\mathfrak{m}$ is just the multiplication by an idempotent. Note that we have an isomorphism $T^\vee_{\mathfrak{m}}(J)\cong H^1(J_{\overline{F}},\mathcal{O})_{\mathfrak{m}}(1)$ of $\bbT_{\mathfrak{m}}[G_F]$-modules.

Let $\mathcal{A}$ be the category whose objects are abelian varieties $A$ defined over $F$, which are equipped with an action of $\bbT$ defined over $F$, and such that there exists a $\bbT$-equivariant isogeny $\phi\colon J\to A$ defined over $F$. The morphisms $\Hom_{\mathcal{A}}(A,B)$ are $\bbT$-equivariant maps $A\to B$ of abelian varieties defined over $F$.

\begin{notation}
    If $A$ is an abelian variety over $F$ and $v\nmid p$ is a finite place of $F$, we write $\mathcal{X}_v(A)$, $\mathcal{T}_v(A)$, $\Phi_v(A)$, and $u_{A,v}$, for the objects $\mathcal{X}(A_{F_v})$, $\mathcal{T}(A_{F_v})$, $\Phi(A_{F_v})$, and $u_{A_{F_v}}$ introduced in \Cref{ssred}.
\end{notation}

By \cite[Lemma 4.5]{Helm}, there exists an abelian variety  $J^{\min}\in \mathcal{A}$ such that $T^\vee_{\mathfrak{m}}J^{\min}$ is free of rank $2$ over $\mathbb{T}_{\mathfrak{m}}$. We will construct such a variety explicitly in our application. 

We have the following Lemma whose proof is exactly the same as in \cite[Lemma 4.6]{Helm}.

\begin{lemma}\label{existsideal}
    If $M$ is a $G_{F}$-stable submodule of $T^\vee_{\mathfrak{m}}J^{\min}$ of finite index, then there exists an ideal $I$ of $\bbT_{\mathfrak{m}}$ such that $M=I T^\vee_{\mathfrak{m}}J^{\min}$.
\end{lemma}
\begin{lemma}\label{annihilatorisogeny}
Let $\phi\colon J^{\min}\to A$ be a $\bbT$-equivariant isogeny, and let $I=\Ann_{\bbT}(\ker\phi)$. Then the induced map
$$ \phi^* \colon T_{\mathfrak{m}}^\vee A\longrightarrow T_{\mathfrak{m}}^\vee J^{\min} $$
identifies $T_{\mathfrak{m}}^\vee A$ with $IT_{\mathfrak{m}}^\vee J^{\min}$.
\end{lemma}
\begin{proof}
    By \cite[Lemma 2.1]{Em03}, we have an exact sequence 
    \begin{equation*}
        0 \to T_{\mathfrak{m}}^\vee A \to T_{\mathfrak{m}}^\vee J^{\min} \to T_{\mathfrak{m}}^\vee \ker(\phi) \to 0.
    \end{equation*}
    By \Cref{existsideal}, there exists an ideal $I'$ of $\bbT_{\mathfrak{m}}$ such that $T_{\mathfrak{m}}^\vee A = I'T_{\mathfrak{m}}^\vee J^{\min}$. Since $T^\vee_{\mathfrak{m}}J^{\min}$ is free over $\bbT_{\mathfrak{m}}$, it follows that $$I'=\Ann_{\bbT_{\mathfrak{m}}}(T_{\mathfrak{m}}^\vee \ker(\phi))=\Ann_{\bbT_{\mathfrak{m}}}((\ker(\phi)[p^\infty]\otimes_{\bbZ_p}\mathcal{O})_{\mathfrak{m}})=\Ann_{\bbT_{\mathfrak{m}}}((\ker(\phi)\otimes_{\bbZ}\mathcal{O})_{\mathfrak{m}})$$ 
    Note that $(\ker(\phi)\otimes_{\bbZ}\mathcal{O})_{\mathfrak{m}}=\ker(\phi)\otimes_{\bbT}\bbT_{\mathfrak{m}}$. By stability of annihilators of finite modules under flat base-change (\stackcite{07T8}), we get that $I'=I\bbT_{\mathfrak{m}}$.
\end{proof}
Let $t$ be a finite place of $F$ such that $J$ has purely toric reduction at $t$. Then every $A\in \mathcal{A}$ has purely toric reduction at $t$, and by \cite[Lemma 4.13]{Helm}, the module $\mathcal{X}_t(A)\otimes \bbQ$ a free of rank one over $\bbT\otimes \bbQ$. Set $\mathcal{X}_t(A)_{\mathfrak{m}}\colonequals (\mathcal{X}_t(A)\otimes_{\bbZ}\mathcal{O})_{\mathfrak{m}}$. The natural surjection $\Hom_{\mathbb{Z}_p}(T_p(A),\mathbb{Z}_p)\to \Hom_{\mathbb{Z}_p}(T_p(A)^{\mathrm{tor}},\mathbb{Z}_p) $ can be identified after localization at $\mathfrak{m}$ with a surjection $T^\vee_{\mathfrak{m}}A\twoheadrightarrow \mathcal{X}_t(A)_{\mathfrak{m}}$ by the discussion in \Cref{ssred}.  We then have the following result.
\begin{theorem}\label{tatemoduletochargrp}
    Let $A\in \mathcal{A}$. The horizontal arrows in the diagram
\begin{center}
    \begin{tikzcd}
    T_{\mathfrak{m}}^\vee J^{\min} \otimes_{\bbT_{\mathfrak{m}}} (\Hom_{\mathcal{A}}(A,J^{\min})\otimes_{\bbZ}\mathcal{O})_{\mathfrak{m}}\arrow[r] \arrow[d,twoheadrightarrow] & T_{\mathfrak{m}}^\vee A \arrow[d,twoheadrightarrow]
            \\    \left(\mathcal{X}_t(J^{\min})_{\mathfrak{m}}\otimes_{\bbT_{\mathfrak{m}}} (\Hom_{\mathcal{A}}(A,J^{\min})\otimes_{\bbZ}\mathcal{O})_{\mathfrak{m}}\right)^{\tf}\arrow[r]& \mathcal{X}_t(A)_{\mathfrak{m}}
    \end{tikzcd}
\end{center}
defined by $x\otimes \eta \mapsto \eta^*(x)$, are isomorphisms. In particular, if $\mathcal{X}_t(J^{\min})_{\mathfrak{m}}$ is free of rank one over $\bbT_{\mathfrak{m}}$, then $T^\vee_{\mathfrak{m}}A\cong (\mathcal{X}_t(A)_{\mathfrak{m}})^2$.
\end{theorem}
\begin{proof}
    Let us fix a $\bbT$-equivariant isogeny $\phi\colon J^{\min}\to A$. This induces an injective $\bbT_{\mathfrak{m}}[G_F]$-equivariant map $\phi^*\colon T_{\mathfrak{m}}^\vee A \to T^\vee_{\mathfrak{m}}J^{\min}$ such that by \Cref{annihilatorisogeny}, $\mathrm{Im}(\phi^*)=I T^\vee_{\mathfrak{m}}J^{\min}$ for $I=\Ann_{\bbT}(\ker(\phi))$. Now for $x\in T_{\mathfrak{m}}^\vee A$, we can write $\phi^*(x)=\sum_i \nu_i x_i$ with $\nu_i\in I$ and $x_i\in T_{\mathfrak{m}}^\vee J^{\min}$. Since the map $\nu_i\colon J^{\min}\to J^{\min}$ annihilates $\ker(\phi)$, it admits a factorization $\nu_i=\eta_i \circ \phi$ with $\eta_i\colon A\to J^{\min}$. Therefore, we see that $x$ is the image of $\sum_i x_i\otimes \eta_i$ by the map 
    $$ T_{\mathfrak{m}}^\vee J^{\min} \otimes_{\bbT_{\mathfrak{m}}} (\Hom_{\mathcal{A}}(A,J^{\min})\otimes_{\bbZ}\mathcal{O})_{\mathfrak{m}}\longrightarrow  T_{\mathfrak{m}}^\vee A $$
    defined in the statement of the proposition. Thus, we have shown that it is surjective. To get injectivity, note that this map is an isomorphism after applying $-\otimes_{\mathcal{O}} K$, since both sides would be free of rank two over $\bbT_{\mathfrak{m}}\otimes_{\mathcal{O}} K$ by Faltings' isogeny theorem. Thus the injectivity follows from $\mathcal{O}$-torsion freeness of $T_{\mathfrak{m}}^\vee J^{\min} \otimes_{\bbT_{\mathfrak{m}}} (\Hom_{\mathcal{A}}(A,J^{\min})\otimes_{\bbZ}\mathcal{O})_{\mathfrak{m}}\cong (\Hom_{\mathcal{A}}(A,J^{\min})\otimes_{\bbZ}\mathcal{O})_{\mathfrak{m}}^2$.
\end{proof}
\subsection{$\lambda$-Shimura degrees}\label{lambdashimuradeg}
Let $J$ and $\bbT$ be as in \Cref{resultofHelm}, and let us fix a maximal ideal $\mathfrak{m}$ of $\bbT$ and an augmentation $\lambda\colon \bbT_{\mathfrak{m}}\to \mathcal{O}$. We suppose that there is a $\bbT$-equivariant principal polarization $\varphi\colon \widehat{J}\xrightarrow{\sim} J$ possibly defined over a finite extension of $F$. We also fix a finite place $t$ of $F$, and we assume that for an abelian variety $J^{\min}$ as in \Cref{resultofHelm}, $\mathcal{X}_t(J^{\min})_{\mathfrak{m}}$ is free of rank one over $\bbT_{\mathfrak{m}}$.

Set $I=\bbT \cap \ker(\lambda)$, which is a saturated ideal of $\bbT$, and consider the abelian variety $A=J/IJ$ so that we have the following diagram
\begin{center}
    \begin{tikzcd}
        0 \arrow[r] &  I J \arrow[r] &  J \arrow[r,"\pi"] & A \arrow[r] & 0
        \\ 0 \arrow[r] &  \widehat{A} \arrow[r,"\widehat{\pi}"] &  \hat{J} \arrow[u,"\varphi"] \arrow[r] & \widehat{I J}\arrow[r] & 0.
    \end{tikzcd}
\end{center}
We define the Galois representation $\rho\colon G_F\to \GL_2(\mathcal{O})$ by the $G_F$-module $(T^\vee_{\mathfrak{m}}(\hat{A})\otimes_{\bbT_{\mathfrak{m}},\lambda}\mathcal{O})^{\tf}$, and we assume that $A$ has toric reduction at a finite place $v$ of $F$.

\begin{definition}
    We define the $\lambda$-Shimura degree of $J$ to be the half-integer $\delta_J \in \tfrac{1}{2}\mathbb{Z}_{\ge 0}$ given by
    $$ \delta_J \colonequals \frac{1}{2}\ell_{\mathcal{O}}\Big(\coker\Big( T^\vee_{\mathfrak{m}}(A)[\mathfrak{p}_\lambda]\xrightarrow{(\pi\circ\varphi\circ \hat{\pi})^*}(T^\vee_{\mathfrak{m}}(\hat{A}) \otimes_{\bbT_{\mathfrak{m}},\lambda}\mathcal{O})^{\tf}\Big)\Big). $$
    We also define the integers $j_v,u_v\in \bbZ_{\ge 0}$ by
    \begin{align*}
        \coker\big( T^\vee_{\mathfrak{m}}(A)[\mathfrak{p}_\lambda] \to \mathcal{X}_v(A)_{\mathfrak{m}}[\mathfrak{p}_\lambda]\big)&\cong \mathcal{O}/\varpi^{u_v},
        \\ \coker\big( \mathcal{X}_v(A)_{\mathfrak{m}}[\mathfrak{p}_\lambda]\to \mathcal{X}_v(J)_{\mathfrak{m}}[\mathfrak{p}_\lambda]\big) &\cong \mathcal{O}/\varpi^{j_v}.
    \end{align*}
\begin{remark}
    All the cokernels in the definition are of the form given above since the maps considered are non-zero and are between free $\mathcal{O}$-modules.
\end{remark}
\end{definition}
\begin{lemma}\label{shimuradegonX}
    The $\lambda$-Shimura degree $\delta_J$ of $J$ is in fact an integer, and we have an isomorphism 
    $$ \coker\Big( \mathcal{X}_{v}(A)_{\mathfrak{m}} [\mathfrak{p}_\lambda]\xrightarrow{(\pi\circ\varphi\circ \hat{\pi})^*} (\mathcal{X}_v(\hat{A})_{\mathfrak{m}} \otimes_{\bbT_{\mathfrak{m}},\lambda}\mathcal{O})^{\tf}\Big)\cong \mathcal{O}/\varpi^{\delta_J-u_v}. $$
\end{lemma}
\begin{proof}
    Since the Weil pairing is Hecke equivariant, it induces an isomorphism of $\bbT_{\mathfrak{m}}$-modules $T_{\mathfrak{m}}^\vee(A) \cong T_{\mathfrak{m}}^\vee (\hat{A})^*$. In particular, \Cref{quotientsiso} yields isomorphisms 
    \begin{align*}
        T^\vee_{\mathfrak{m}}(A)[\mathfrak{p}_\lambda]\cong (T_{\mathfrak{m}}^\vee(\hat{A})\otimes_{\bbT_{\mathfrak{m}},\lambda}\mathcal{O})^* \quad \text{ and, }\quad (T^\vee_{\mathfrak{m}}(A)\otimes_{T_{\mathfrak{m}},\lambda}\mathcal{O})^{\tf}\cong (T^\vee_{\mathfrak{m}}(\hat{A})[\mathfrak{p}_\lambda])^*.
    \end{align*}
    Since $A$ has toric reduction, by Grothendieck's orthogonality theorem, we have an exact sequence of $\bbT_{\mathfrak{m}}$-modules
    $$ 0 \to \mathcal{X}_v(\hat{A})_\mathfrak{m}^*\to T^\vee_{\mathfrak{m}}(A)\to \mathcal{X}_v(A)_{\mathfrak{m}}\to 0, $$
    and we also have the same exact sequence if we replace $A$ with $\hat{A}$. Therefore, we have the following short exact sequence of complexes (the complexes are seen in the vertical direction, and the exact sequence is in the horizontal direction)
    \begin{center}
        \begin{tikzcd} \vspace{-6pt}
          & 0 \arrow[d] & 0 \arrow[d] & 0 \arrow[d]
            \\ 0 \arrow[r] & \mathcal{X}_v(\hat{A})^*_{\mathfrak{m}}[\mathfrak{p}_\lambda] \arrow[r] \arrow[d] & (\mathcal{X}_v(A)^*_{\mathfrak{m}} \otimes_{\bbT_{\mathfrak{m}},\lambda}\mathcal{O})^{\tf} \arrow[r] \arrow[d] & N' \arrow[r] \arrow[d] & 0
            \\ 0 \arrow[r] & T^\vee_{\mathfrak{m}}(A)[\mathfrak{p}_\lambda] \arrow[r] \arrow[d] & (T^\vee_{\mathfrak{m}}(\hat{A}) \otimes_{\bbT_{\mathfrak{m}},\lambda}\mathcal{O})^{\tf} \arrow[r] \arrow[d] & N \arrow[r] \arrow[d] & 0
            \\ 0 \arrow[r] & \mathcal{X}_v(A)_{\mathfrak{m}}[\mathfrak{p}_\lambda] \arrow[r] \arrow[d]  & (\mathcal{X}_v(\hat{A})_{\mathfrak{m}}\otimes_{\bbT_\mathfrak{m},\lambda}\mathcal{O})^{\tf} \arrow[r] \arrow[d] & N'' \arrow[r] \arrow[d] & 0
            \\ & 0 & 0 & 0
        \end{tikzcd}
    \end{center}
    Taking the long exact sequence associated to cohomology, we get the exact sequence
    \begin{align*}
        0 \to \ker( N'\to N) \to 0 &\to \coker( T^\vee_{\mathfrak{m}}(A)[\mathfrak{p}_\lambda] \to \mathcal{X}_v(A)_{\mathfrak{m}}[\mathfrak{p}_\lambda])^* \to \ker(N\to N'')/\im(N'\to N)
        \\ &\to \coker( T^\vee_{\mathfrak{m}}(A)[\mathfrak{p}_\lambda] \to \mathcal{X}_v(A)_{\mathfrak{m}}[\mathfrak{p}_\lambda]) \to 0.
    \end{align*}
    Therefore, we have 
$$2\delta=\ell_{\mathcal{O}}(N)=\ell_{\mathcal{O}}(N')+\ell_{\mathcal{O}}(\ker(N\to N'')/\im(N'\to N))+\ell_{\mathcal{O}}(N'')= \ell_{\mathcal{O}}(N')+2u_v+\ell_{\mathcal{O}}(N'').$$
We finish the proof by noting that
taking the dual of the first row, we have that $N''\cong \Ext^1_{\mathcal{O}}(N',\mathcal{O})\cong \Hom_{\mathcal{O}}(N',E/\mathcal{O})$. In particular, we have that $\ell_{\mathcal{O}}(N')=\ell_{\mathcal{O}}(N'')$.
\end{proof}

\begin{definition}
    We define the module $\Phi^{\tf}_v(A)$ by the following exact sequence coming from the monodromy pairing
\begin{equation}\label{mondromypairingatf}
     0 \to (\mathcal{X}_v(\hat{A})_{\mathfrak{m}}\otimes_{\bbT_{\mathfrak{m}},\lambda}\mathcal{O})^{\tf}\to (\mathcal{X}_v(A)_{\mathfrak{m}}^*\otimes_{\bbT_{\mathfrak{m}},\lambda}\mathcal{O})^{\tf}\to \Phi^{\tf}_v(A) \to 0.
\end{equation}
It is a quotient of $\Phi_v(A)_{\mathfrak{m}}\otimes_{\bbT_{\mathfrak{m}},\lambda}\mathcal{O}$.
\end{definition}

\begin{lemma}\label{someequalities}
    We have the following equalities
    $$ c_v(\rho)=\ell_{\mathcal{O}}(\Phi_v^{\tf}(A))+u_v \quad \text{ and, }\quad j_v=\ell_{\mathcal{O}}\big(\coker(\Phi^{\tf}_v(J)\to \Phi^{\tf}_v(A))\big). $$
\end{lemma}
\begin{proof}
For the first equality, recall that $\rho=(T^\vee_{\mathfrak{m}}(\hat{A})\otimes_{\bbT_{\mathfrak{m}},\lambda}\mathcal{O})^{\tf}$. We have 
$$ ((T^\vee_{\mathfrak{m}}(\hat{A})\otimes_{\bbT_{\mathfrak{m}},\lambda}\mathcal{O})^{\tf})^{I_{v}}= \ker\big( (T^\vee_{\mathfrak{m}}(\hat{A})\otimes_{\bbT_{\mathfrak{m}},\lambda}\mathcal{O})^{\tf} \to (\mathcal{X}_v(\hat{A})_{\mathfrak{m}}\otimes_{\bbT_{\mathfrak{m}},\lambda}\mathcal{O})^{\tf}\big).$$
Recall that $\sigma_v\in I_v$ is a fixed lift a generator of the $\bbZ_p$-quotient of $I_v$. 
In the exact sequence \eqref{mondromypairingatf} associated to the monodromy pairing, the leftmost map is given by applying $\sigma_v-1$.  And since $$\coker\left[(\mathcal{X}_v(A)^*_{\mathfrak{m}}\otimes_{\bbT_{\mathfrak{m}},\lambda}\mathcal{O})^{\tf}\to \ker\Big( (T^\vee_{\mathfrak{m}}(\hat{A})\otimes_{\bbT_{\mathfrak{m}},\lambda}\mathcal{O})^{\tf} \to (\mathcal{X}_v(\hat{A})_{\mathfrak{m}}\otimes_{\bbT_{\mathfrak{m}},\lambda}\mathcal{O})^{\tf}\Big)\right]\cong \mathcal{O}/\varpi^{u_v}, $$
we find that the action of $\sigma_v\in I_v$ is given by the matrix   $$  \begin{pmatrix}
        1 & s\varpi^{u_v+k_v}
        \\ 0 & 1 
    \end{pmatrix},$$ where $k_v=\ell_{\mathcal{O}}(\Phi^{\tf}_v(A))$ and $s\in \mathcal{O}^\times$. This gives the first equality.

For the second equality, it follows from \cite[Lemma 3.1 (iii)]{Em03} that the leftmost vertical map in the diagram 
\begin{center}
    \begin{tikzcd}
        0 \arrow[r] & (\mathcal{X}_v(\hat{J})_{\mathfrak{m}}\otimes_{\bbT_{\mathfrak{m}},\lambda}\mathcal{O})^{\tf}\arrow[r] \arrow[d] & (\mathcal{X}_v(J)_{\mathfrak{m}}[\mathfrak{p}_\lambda])^*\arrow[r] \arrow[d] & \Phi_v^{\tf}(J)\arrow[r] \arrow[d] & 0
        \\0 \arrow[r] & (\mathcal{X}_v(\hat{A})_{\mathfrak{m}} \otimes_{\bbT_{\mathfrak{m}},\lambda}\mathcal{O})^{\tf}\arrow[r]  &(\mathcal{X}_v(A)_{\mathfrak{m}}[\mathfrak{p}_\lambda])^*\arrow[r]  & \Phi_v^{\tf}(A) \arrow[r]  & 0,
    \end{tikzcd}
\end{center}
is surjective. Applying the snake lemma, we find that  
$$\Ext^1_\mathcal{O}\big(\mathcal{X}_v(J)_{\mathfrak{m}}[\mathfrak{p}_\lambda]/\mathcal{X}_v(A)_{\mathfrak{m}}[\mathfrak{p}_\lambda],\mathcal{O}\big)\cong \coker\big(\Phi^{\tf}_v(J)\to \Phi^{\tf}_v(A)\big).$$
We conclude by noting that 
\begin{equation*}
    \ell_{\mathcal{O}}\big(\Ext^1_\mathcal{O}(\mathcal{X}_v(J)_{\mathfrak{m}}[\mathfrak{p}_\lambda]/\mathcal{X}_v(A)_{\mathfrak{m}}[\mathfrak{p}_\lambda],\mathcal{O})\big)=\ell_{\mathcal{O}}\big((\mathcal{X}_v(J)_{\mathfrak{m}}[\mathfrak{p}_\lambda]/\mathcal{X}_v(A)_{\mathfrak{m}}[\mathfrak{p}_\lambda])_{\tors}\big). \qedhere
\end{equation*}
\end{proof}

\begin{lemma}\label{shdegtocongmodule}
    We have the following equality
    $$ \delta_J = j_v+u_v+ \ell_{\mathcal{O}}(\Psi_\lambda(\mathcal{X}_v(J)_{\mathfrak{m}})). $$
    For $v=t$, we have $j_t=u_t=0$, and hence $\delta_J =\ell_{\mathcal{O}}(\Psi_\lambda(\mathcal{X}_t(J)_{\mathfrak{m}}))$.
\end{lemma}
\begin{proof}
    The statement follows from observing the following diagram
\begin{center}
    \begin{tikzcd}
     & \mathcal{X}_v(A)_{\mathfrak{m}}[\mathfrak{p}_\lambda] \arrow[d,hookrightarrow]
    \\ 0 \arrow[r] & \mathcal{X}_v(J)_{\mathfrak{m}}[\mathfrak{p}_\lambda] \arrow[r]  & (\mathcal{X}_v(J)_{\mathfrak{m}}\otimes_{\bbT_{\mathfrak{m}},\lambda}\mathcal{O})^{\tf}\arrow[r] \arrow[d,"\cong"] & \Psi_\lambda(\mathcal{X}_v(J)_{\mathfrak{m}})  \arrow[r] & 0
    \\ & & (\mathcal{X}_v(\hat{A})_{\mathfrak{m}}\otimes_{\bbT_{\mathfrak{m}},\lambda}\mathcal{O})^{\tf}.
    \end{tikzcd}
\end{center}
For $v=t$, consider the following diagram with exact rows
\begin{center} 
    \begin{tikzcd} 
        0 \arrow[r] & T^\vee_{\mathfrak{m}}(A) \arrow[r] \arrow[d] & T^\vee_{\mathfrak{m}}(J) \arrow[r] \arrow[d] & T^\vee_{\mathfrak{m}}(IJ) \arrow[r] & 0
        \\ 0 \arrow[r] &\mathcal{X}_t(A)_{\mathfrak{m}} \arrow[r] &\mathcal{X}_t(J)_{\mathfrak{m}}.
    \end{tikzcd} 
\end{center} 
Taking the $\lambda$-equivariant part, we get
\begin{center} \vspace{-6pt}
    \begin{tikzcd}
        & T^\vee_{\mathfrak{m}}(A)[\mathfrak{p}_\lambda] \arrow[r,"\sim"] \arrow[d] & T^\vee_{\mathfrak{m}}(J)[\mathfrak{p}_\lambda]   \arrow[d,twoheadrightarrow]
        \\ 0 \arrow[r] &\mathcal{X}_t(A)_{\mathfrak{m}}[\mathfrak{p}_\lambda] \arrow[r]&\mathcal{X}_t(J)_{\mathfrak{m}}[\mathfrak{p}_\lambda]. 
    \end{tikzcd} 
\end{center} \vspace{-6pt}
The top map is an isomorphism since its cokernel is torsion, and $T^\vee_{\mathfrak{m}}(IJ)$ is $\varpi$-torsion free. The fact that $T^\vee_{\mathfrak{m}}(J)[\mathfrak{p}_\lambda]\to \mathcal{X}_t(J)_{\mathfrak{m}}[\mathfrak{p}_\lambda]$ is surjective follows from \Cref{tatemoduletochargrp} since in this case $\mathcal{X}_t(J_{\min})_{\mathfrak{m}}$ is free over $\bbT_{\mathfrak{m}}$, and so the horizontal map in the statement of the theorem admits a section. Therefore, we get that the other two maps are surjective, which give that $u_t=j_t=0$. 
\end{proof}

\begin{lemma}\label{scalartoShimuradeg}
    Let us write $\mathcal{X}_v(J)_{\mathfrak{m}}[\mathfrak{p}_\lambda]=\mathcal{O}\cdot z$. Then we have
    $$ \ord_{\varpi}(u_J(z,\varphi^*z))=\delta_J +c_v(\rho)- 2j_v-2u_v.$$
\end{lemma}
\begin{proof}
    Let us write $\mathcal{X}_v(A)_{\mathfrak{m}}[\mathfrak{p}_\lambda]=\mathcal{O}\cdot x$ and $(\mathcal{X}_v(\hat{A})_{\mathfrak{m}}\otimes_{\bbT_{\mathfrak{m}},\lambda}\mathcal{O})^{\tf}=\mathcal{O}\cdot y$. Then by \Cref{shimuradegonX}, we have
    $$ u_A(x,(\pi\circ \varphi \circ \hat{\pi})^*x)=\varpi^{\delta_J -u_v}u_A(x,y). $$
    Moreover by \eqref{mondromypairingatf}, we have $\ord_{\varpi}(u_A(x,y))=\ell_{\mathcal{O}}(\Phi^{\tf}_v(A))$. On the other hand, we have by equivariance of the monodromy pairing that 
    $$ u_A(x,(\pi\circ \varphi \circ \hat{\pi})^*x)=u_J(\pi^*x,\varphi^*\pi^*x).$$
    We also have $\ord_{\varpi}(u_A(\pi^*x,\varphi^*\pi^*x))=2j_v+ \ord_{\varpi}(u_J(z,\varphi^*z))$. The result now follows by combining these equalities with the first equality of \Cref{someequalities}.
\end{proof}
\section{Comparison of $\lambda$-Shimura degrees}\label{complambdashimura}
This section is dedicated to the proof of \Cref{maintheorem}. Analogous to the approach of Ribet and Takahashi, the key component in the proof is Ribet’s exact sequence, which realizes the Jacquet--Langlands correspondence on the character groups of the toric parts of the Jacobians of Shimura curves. To establish this exact sequence, we review in \Cref{integralmodelsofshimuracurves} the study of the integral models of Shimura curves both at places dividing and not dividing the discriminant. Furthermore, to connect the definite case and indefinite case, we establish in \Cref{chargroupvnotinQ} a relation between character groups at places not dividing the discriminant, and spaces of definite quaternionic forms. 
\subsection{Character group of Jacobians of curves}\label{chargrpJac}
Let $v$ be a finite place of $F$, and let $\mathcal{C}$ be an admissible curve over $\mathcal{O}_{F_v}$ with generic fiber $\mathcal{C}_{F_v}=C$ (we refer to \cite[\S 3]{JL85} for the definition of an admissible curve). Let $J=\Pic^0(C/F_v)$ be the Jacobian of the curve $C$, and $\mathcal{J}$ its Néron model over $\mathcal{O}_{F_v}$. The special fibre $\mathcal{C}_{k_v}$ can be described using the graph $\mathcal{G}$ defined as follows
\begin{itemize}
    \item The set of vertices of $\mathcal{G}$ is the set $\text{Ver}(\mathcal{G})$ of irreducible components of $\mathcal{C}_{k_v}$.
    \item The set of edges of $\mathcal{G}$ is the set $\text{Ed}(\mathcal{G})$ of singular points of $\mathcal{C}_{k_v}$.
    \item The edge corresponding to a singular point $e\in \text{Ed}(\mathcal{G})$ connects the two vertices corresponding to the two irreducible components of $\mathcal{C}_{k_v}$ meeting at the point $e$.
\end{itemize}
By \cite[Exposé IX 12.3.7]{SGA7I} and the argument in \cite[§2]{Ribet90}, there is an isomorphism
\begin{equation*}
    \mathcal{T}_v(J)\cong H^1(\mathcal{G},\bbZ)\otimes \mathbb{G}_m.
\end{equation*}
which is equivalent to the isomorphism $    \mathcal{X}_v(J)\cong H_1(\mathcal{G},\bbZ)$.
This allows us to give a combinatorial description of the character group as follows.

Let $\overline{\mathcal{G}}$ be the bouquet of circles obtained by collapsing all vertices of $\mathcal{G}$ to a single point. The map $\mathcal{G}\to \overline{\mathcal{G}}$ induces an inclusion $\mathcal{X}_v(J)\cong H_1(\mathcal{G},\mathbb{Z})\to H_1(\overline{\mathcal{G}},\mathbb{Z})$. We choose an orientation on $\mathcal{G}$ so that for every $e\in \text{Ed}(\mathcal{G})$, we write $i(e)$ and $j(e)$ for the source and target of $e$. This orientation determines an isomorphism $H_1(\overline{\mathcal{G}},\mathbb{Z})\cong \mathbb{Z}^{\mathrm{Ed}(\mathcal{G})}$, and hence a non-canonical inclusion $\mathcal{X}_v(J)\hookrightarrow \mathbb{Z}^{\mathrm{Ed}(\mathcal{G})}$. Let $(\bbZ^{\text{Ver}(\mathcal{G})})_0$ and $(\mathbb{Z}^{\text{Ed}(\mathcal{G})})_0$ denote the kernels of the augmentation maps $\bbZ^{\text{Ver}(\mathcal{G})}\rightarrow \bbZ$ and $\bbZ^{\text{Ed}(\mathcal{G})}\rightarrow \bbZ$, respectively. Then by \cite[Proposition 2.2]{Ribet90}, we obtain an exact sequence
\begin{equation}\label{dualgraphexactsequence}
    0\rightarrow \mathcal{X}_v(J) \rightarrow \bbZ^{\text{Ed}(\mathcal{G})}\xrightarrow{\alpha} (\bbZ^{\text{Ver}(\mathcal{G})})_0 \rightarrow 0,
\end{equation}
where the left-most map is the augmentation map, and $\alpha(e)=j(e)-i(e)$. Moreover, given that $J$ is self-dual, the monodromy pairing \eqref{monodromypairing} on $\mathcal{X}_v(J)$ is the restriction of the Euclidean pairing on $\bbZ^{\text{Ed}(\mathcal{G})}$.
\subsection{Integral models of Shimura curves}\label{integralmodelsofshimuracurves}
Let $Q$ be a finite set of finite places of $F$ of indefinite type. We assume that we have a decomposition $Q=\overline{Q}\sqcup \{v,w\}$, and we set $\widetilde{Q}=\overline{Q}\cup\{v\}$.
\subsubsection{ Integral model at a place not dividing the discriminant}\label{integralmodel1} 

\begin{theorem}[{\cite[Theorem 4.3]{Tho16}}]
    Let $\leftidx{_v}{\mathcal{J}}{_{\overline{Q}}}$ denote the subset of $\mathcal{J}_{\overline{Q}}$ consisting of subgroups $\overline{U}=\prod_{v'} \overline{U}_{v'}\subset G_{\overline{Q}}(\mathbb{A}_F^\infty)$ such that $\overline{U}_v=\GL_2(\mathcal{O}_{F_v})$. Let $\overline{U}\in \leftidx{_v}{\mathcal{J}}{_{\overline{Q}}}$, then the morphisms
    \begin{align*}
        &X_{\overline{Q}}(\overline{U})_{F_v}\longrightarrow \Spec F_v
        \\ &X_{\overline{Q}}(\overline{U}_0(v))_{F_v}\longrightarrow \Spec F_v
    \end{align*}
    extend canonically to flat projective morphisms 
       \begin{align*}
        &{\mathfrak{X}_{\overline{Q}}}(\overline{U})\longrightarrow \Spec \mathcal{O}_{F_v}
        \\ &{\mathfrak{X}_{\overline{Q}}}(\overline{U}_0(v))\longrightarrow \Spec \mathcal{O}_{F_v}.
    \end{align*} 
    Moreover, ${\mathfrak{X}}{_{\overline{Q}}}(\overline{U})$ is smooth over $\mathcal{O}_{F_v}$, and ${\mathfrak{X}}{_{\overline{Q}}}(\overline{U}_0(v))$ is regular and semi-stable over $\mathcal{O}_{F_v}$. The action of the group $G_{\overline{Q}}(\mathbb{A}_F^{a,v,\infty})$ on the projective systems $\{X_{\overline{Q}}(\overline{U})_{F_v}\}_{U\in \leftidx{_v}{\mathcal{J}}{_{\overline{Q}}}} $ and $\{X_{\overline{Q}}(\overline{U}_0(v))_{F_v}\}_{\overline{U}\in \leftidx{_v}{\mathcal{J}}{_{\overline{Q}}}} $ extend canonically to actions on the projective systems $\{ {\mathfrak{X}_{\overline{Q}}}(\overline{U})\}_{\overline{U}\in \leftidx{_v}{\mathcal{J}}{_{\overline{Q}}}} $ and $\{{\mathfrak{X}_{\overline{Q}}}(\overline{U}_0(v))\}_{\overline{U}\in \leftidx{_v}{\mathcal{J}}{_{\overline{Q}}}} $.
\end{theorem}
Let $\overline{U}=\prod_{v'} \overline{U}_{v'}\in \leftidx{_{v}}{\mathcal{J}}_{\overline{Q}}$. We now define a good subgroup $\widetilde{U}=\prod_{v'} \widetilde{U}_{v'}\in \mathcal{J}_{\widetilde{Q}}$ by setting $\widetilde{U}_{v'}=\overline{U}_{v'}$ if $v\neq v'$, and $\widetilde{U}_v$ is the unique maximal compact subgroup of $G_{\widetilde{Q}}(F_v)$.

 Let $\kappa(v)$ be the residue field of a fixed algebraic closure $\overline{F}_v$. Denote by ${\widetilde{\mathfrak{X}}}_{\overline{Q}}(\overline{U}_0(v))_{\kappa(v)}$ the normalization of the 
 fibre ${\mathfrak{X}}_{\overline{Q}}(\overline{U}_0(v))_{\kappa(v)}$, and by $\text{Sing}_{\overline{Q}}(\overline{U})$ its set of singular points. Then there is an isomorphism of projective systems of smooth $\kappa(v)$-schemes with a right $G_{\overline{Q}}(\mathbb{A}_F^{a,v,\infty})$-action:
\begin{equation}\label{ver(g)}
    \big\{{\widetilde{\mathfrak{X}}}_{\overline{Q}}(\overline{U}_0(v))_{\kappa(v)}\big\}_{\overline{U}\in \leftidx{_v}{\mathcal{J}}{_{\overline{Q}}}} \cong \big\{{\mathfrak{X}}_{\overline{Q}}(\overline{U})_{\kappa(v)} \sqcup {\mathfrak{X}}_{\overline{Q}}(\overline{U})_{\kappa(v)} \big\}_{\overline{U}\in \leftidx{_v}{\mathcal{J}}{_{\overline{Q}}}}.
\end{equation}
Moreover, there is an isomorphism of projective systems of sets with a right $G_{\overline{Q}}(\mathbb{A}_F^{a,v,\infty})$-action
\begin{equation}\label{sspoints}
    \big\{\text{Sing}_{\overline{Q}}(\overline{U})\big\}_{\overline{U}\in \leftidx{_v}{\mathcal{J}}{_{\overline{Q}}}} \cong \big\{Y_{\widetilde{Q}}(\widetilde{U})\big\}_{\overline{U}\in \leftidx{_v}{\mathcal{J}}{_{\overline{Q}}}}.
\end{equation}
We recall that the definition of $Y_{\widetilde{Q}}(\widetilde{U})$ is given in \Cref{Shimura set}.

\begin{lemma}\label{chargroupvnotinQ} There is an isomorphism  $$\mathcal{X}_v\big(J_{\overline{Q}}(\overline{U}_0(v))\big)\xrightarrow{\sim} \big(\bbZ^{Y_{\widetilde{Q}}(\widetilde{U})}\big)_0.$$ Moreover, if $\mathfrak{m}$ is a non-Eisenstein maximal ideal of $\overline
{\bbT}^{\overline{Q}}_{\mathcal{O}}(\overline{U}_0(v))$, then there is an isomorphism of $\overline{\bbT}{}^{\overline{Q}}_{\mathcal{O}}(\overline{U}_0(v))_{\mathfrak{m}}$-modules
    $$ \mathcal{X}_v(J_{\overline{Q}}(\overline{U}_0(v)))_{\mathfrak{m}}=\big(\mathcal{X}_v(J_{\overline{Q}}(\overline{U}_0(v)))\otimes_{\bbZ}\mathcal{O}\big)_{\mathfrak{m}}\xrightarrow{\sim} M_{\widetilde{Q}}(\widetilde{U}), $$
    where $\overline{\bbT}{}^{\overline{Q}}_{\mathcal{O}}(\overline{U})_{\mathfrak{m}}$ acts on the second module via its quotient map $\overline{\bbT}{}^{\overline{Q}}_{\mathcal{O}}(\overline{U})_{\mathfrak{m}}\twoheadrightarrow \overline{\bbT}{}^{\widetilde{Q}}_{\mathcal{O}}(\widetilde{U})_{\mathfrak{m}}$. 
\end{lemma}
\begin{proof}
Let $\mathcal{G}$ be the graph associated to $\mathfrak{X}_{\overline{Q}}(\overline{U}_0(v))$.  By \eqref{ver(g)} and \eqref{sspoints}, we have that $\text{Ed}(\mathcal{G})=Y_{\widetilde{Q}}(\widetilde{U})$, and $\text{Ver}(\mathcal{G})=\{C_1,C_2\}$. We choose the orientation on $\mathcal{G}$ such that $i(e)=C_1$ and $j(e)=C_2$ for every $e\in \text{Ed}(\mathcal{G})$. The exact sequence \eqref{dualgraphexactsequence} then translates to
$$ 0 \to \mathcal{X}_v(J_{\overline{Q}}(\overline{U}))\to \bbZ^{Y_{\widetilde{Q}}(\widetilde{U})}\xrightarrow{\deg} \bbZ \to 0.  $$
This gives the first part of the proof. For the second part of the proof, we have that $H_{\widetilde{Q}}(\widetilde{U})=(\bbZ^{Y_{\widetilde{Q}}(\widetilde{U})})\otimes \mathcal{O}$, and $\mathcal{X}_v(J_{\overline{Q}}(\overline{U}))\otimes \mathcal{O}$ can be identified with the subspace $H_{\widetilde{Q}}(\widetilde{U})_0$ of $H_{\widetilde{Q}}(\widetilde{U})$ consisting of functions $f$ satisfying $\sum_y f(y)=0$.
We recall that there is a Hecke-equivariant perfect pairing
\begin{align*}
    H_{\widetilde{Q}}(\widetilde{U})&\times H_{\widetilde{Q}}(\widetilde{U})\longrightarrow \mathcal{O},
\end{align*}
and since $H_{\widetilde{Q}}(\widetilde{U})_0$ is stable under the involution $W_{\widetilde{U}}$, the element $f^{\text{triv}}\colon Y_{\widetilde{Q}}(\widetilde{U})\to \mathcal{O},y\mapsto 1$ is orthogonal to $H_{\widetilde{Q}}(\widetilde{U})_0$. Therefore, the above pairing restricts to a Hecke-equivariant perfect pairing
$$ \mathcal{O}\cdot f^{\text{triv}} \times \big(H_{\widetilde{Q}}(\widetilde{U})/H_{\widetilde{Q}}(\widetilde{U})_0\big)\longrightarrow \mathcal{O}. $$
Since $f^{\text{triv}}$ factors through the reduced norm, every maximal ideal $\mathfrak{m}$ of $\bbT^{\widetilde{Q}}_{\mathcal{O}}(\widetilde{U})$ in the support of $f^{\text{triv}}$ is Eisenstein (see for example \cite[proof of lemma 3.1]{Tay06}). This gives the result.
\end{proof}

\subsubsection{ Integral model at a place dividing the discriminant}

Let $\breve{F}$ be the completion of the maximal unramified extension of $F$ inside $\overline{F}_w$, and let $\breve{\mathcal{O}}$ be its ring of integers. We write $\Omega^2\rightarrow \Spf \mathcal{O}_{F_w}$ for $F_w$-adic Drinfeld upper-half plane, $\Omega^2_{\breve{\mathcal{O}}}\rightarrow \Spf \breve{\mathcal{O}}$ for its base extension to $\breve{\mathcal{O}}$, and we define:
\begin{equation*}
    \mathcal{M} \cong \Omega^2_{\breve{\mathcal{O}}} \times \left[(D_Q\otimes_F F_w)^\times/(\mathcal{O}_Q\otimes_{\mathcal{O}_{F}}\mathcal{O}_{F_w})^\times \right].
\end{equation*}
Recall that $\widetilde{Q}=Q\setminus\{w\}$. We fix an isomorphism $G_Q(\mathbb{A}_F^{w,\infty})\cong G_{\widetilde{Q}}(\mathbb{A}_F^{w,\infty})$ compatible with the identifications $G_Q(\mathbb{A}_F^{Q,\infty})\cong \GL_2(\mathbb{A}_F^{Q,\infty})$ and $G_{\widetilde{Q}}(\mathbb{A}_F^{\widetilde{Q},\infty})\cong\GL_2(\mathbb{A}_F^{\widetilde{Q},\infty})$.

The reduced norm gives an isomorphism $(D_Q\otimes_F F_w)^\times/(\mathcal{O}_Q\otimes_{\mathcal{O}_{F}}\mathcal{O}_{F_w})^\times\cong F_w^\times/\mathcal{O}_{F_w}^\times$. We let $G_{\widetilde{Q}}(F_w)\cong \GL_2(F_w)$ act on $\mathcal{M}$ on the left via its usual action via $\PGL_2(F_w)$ on the leftmost factor and by the determinant and this isomorphism on the rightmost factor.
\begin{theorem}[{\cite[Theorem 4.5]{Tho16}}]\label{integralmodelvinQ}
    Let $U\in \mathcal{J}_Q$. The morphism $X_Q(U)\rightarrow \Spec \breve{F}$ extends canonically to a flat projective morphism ${\mathfrak{X}_{Q}}(U)\rightarrow \breve{\mathcal{O}}$ with ${\mathfrak{X}_{Q}}(U)$ regular and semi-stable over $\breve{\mathcal{O}}$, and there is an isomorphism of formal schemes over $\Spf \breve{\mathcal{O}}$
    \begin{equation}\label{CD}
       { \widehat{\mathfrak{X}}_{Q}}(U)\cong G_{\widetilde{Q}}(F)\backslash \big(\mathcal{M}\times G_{\widetilde{Q}}(\mathbb{A}_F^{w,\infty})/U^w\big),
    \end{equation}
where ${ \widehat{\mathfrak{X}}_{Q}}(U)$ is the $\varpi_w$-adic completion of $ { \mathfrak{X}_{Q}}(U)$. Moreover, there is a canonical right $G_Q(\mathbb{A}_F^{a,\infty})$-action on the projective system $\{ { \widehat{\mathfrak{X}}_{Q}}(U)\}_{U\in \mathcal{J}_Q}$ of $\breve{\mathcal{O}}$-schemes, extending the action on $\{X_Q(U)_{\breve{F}}\}_{U\in \mathcal{J}_Q}$, and making the system of isomorphisms \eqref{CD} $G_Q(\mathbb{A}_F^{w,\infty})$-equivariant.
\end{theorem}
\subsection{A geometric realization of Jacquet--Langlands}\label{chargrpofshimuracurve}
Let $Q$ be a finite set of finite places of $F$ as in \Cref{integralmodelsofshimuracurves}. Let  $\overline{U}=\prod_{v'} \overline{U}_{v'}\in \mathcal{J}_{\overline{Q}}$ be a good subgroup that is unramified at $w$ and such that $\overline{U}_v\cong \Gamma_0(v)$. We associate to $\overline{U}$ the good subgroups $\widetilde{U}=\prod_{v'} \widetilde{U}_{v'}\in \mathcal{J}_{\widetilde{Q}}$, $U=\prod_{v'} U_{v'}\in \mathcal{J}_{Q}$ defined by the following: 
\begin{itemize}
    \item $\widetilde{U}_{v'}\cong \overline{U}_{v'}$ for all $v'\neq v$, and $\widetilde{U}_{v}$ is the maximal open compact subgroups of $G_{\widetilde{Q}}(F_{v})$.
    \item $U_{v'}\cong \overline{U}_{v'}$ for all $v'\not\in \{v,w\}$, and  $U_{v}$ and $U_{w}$ are the maximal open compact subgroups of $G_Q(F_{v})$ and $G_Q(F_w)$ respectively.
\end{itemize}

Let $\text{BT}$ be the Bruhat-Tits tree of the group $\PGL_2(F_{w})$. Its set of vertices is equal to the set of $\mathcal{O}_{F_{w}}$-lattices $M\subseteq F_{w}^2$ taken up to $F_{w}^\times$-multiple, which is in a natural bijection with $\PGL_2(F_{w})/\PGL_2(\mathcal{O}_{F_{w}})$. Two vertices $[M]$ and $[M']$ are joined by an edge if and only if we can choose representatives $M,M'$ such that $M\subset M'$ and $[M':M]=q_{w}$. Hence, the set of edges is in bijection with $\GL_2(F_{w})^+/\Gamma_0(w)F_{w}^\times$, where $ \GL_2(F_{w})^+$ denotes the invertible matrices whose determinant has even $w$-adic valuation. Note that $\text{BT} $ is equipped with an action of the group $\PGL_2(F_{w})$. 

 By \Cref{integralmodelvinQ}, the dual graph of $J_{Q}(U)_{F_{w}}$ is
\begin{equation*}
    G_{\widetilde{Q}}(F)\backslash\big(\text{BT}\times F_{w}^\times/\mathcal{O}_{F_{w}}^\times \times G_{\widetilde{Q}}(\mathbb{A}_F^{w,\infty})/\widetilde{U}^{w}\big).
\end{equation*}
Hence, its set of vertices is equal to
\begin{equation*}
    G_{\widetilde{Q}}(F)\backslash\big(\text{Ver}(\text{BT})\times F_{w}^\times/\mathcal{O}_{F_{w}}^\times \times G_{\widetilde{Q}}(\mathbb{A}_F^{w,\infty})/\widetilde{U}^{w}\big).
\end{equation*}
The group $\GL_2(F_{w})$ has two orbits on the set $\text{Ver}(\text{BT})\times F_{w}^\times/\mathcal{O}_{F_{w}}^\times$ whose representatives are $(\mathcal{O}_{F_{w}}^2,1)$ and $(\mathcal{O}_{F_{w}}^2,\varpi_{w})$. Both of these have stabilizers isomorphic to $\GL_2(\mathcal{O}_{F_{w}})$, and therefore the set of vertices of the dual graph is equal to
\begin{equation*}
     G_{\widetilde{Q}}(F)\backslash\big( (\GL_2(F_{w})/\GL_2(\mathcal{O}_{F_{w}}))^2 \times G_{\widetilde{Q}}(\mathbb{A}_F^{w,\infty})/\widetilde{U}^{w}\big)\cong Y_{\widetilde{Q}}(\widetilde{U})\sqcup Y_{\widetilde{Q}}(\widetilde{U}).
\end{equation*}
Similarly, the set of edges of the dual graph of $J_{Q}(U)_{F_{w}}$ is
\begin{equation*}
    G_{\widetilde{Q}}(F)\backslash\big(\text{Ed}(\text{BT})\times F_{w}^\times/\mathcal{O}_{F_{w}}^\times \times G_{\widetilde{Q}}(\mathbb{A}_F^{w,\infty})/\widetilde{U}^{w}\big).
\end{equation*}
Via the isomorphism  $\GL_2(F_{w})^+/\Gamma_0(w)F_{w}^\times \times F_{w}^\times/\mathcal{O}_{F_{w}}^\times\cong \GL_2(F_{w})/\Gamma_0(w)$, we get that the set of edges of the dual graph is isomorphic to $Y_{\widetilde{Q}}(\widetilde{U}_{0}(w))$. We get by (\ref{dualgraphexactsequence}) the following exact sequence
\begin{equation*}
    0 \rightarrow \mathcal{X}_{w}(J_{Q}(U))\rightarrow \mathbb{Z}^{Y_{\widetilde{Q}}(\widetilde{U}_0(w))}\xrightarrow{\alpha} (\mathbb{Z}^{Y_{\widetilde{Q}}(\widetilde{U})})^2 \rightarrow 0.
\end{equation*} 
The map $\alpha$ is the one induced by the two degeneracy maps $$Y_{\widetilde{Q}}(\widetilde{U}_0(w))\rightarrow Y_{\widetilde{Q}}(\widetilde{U}),\  y\mapsto y \quad \text{ and } \quad Y_{\widetilde{Q}}(\widetilde{U}_0(w))\rightarrow Y_{\widetilde{Q}}(\widetilde{U}),\  y\mapsto \begin{pmatrix}
    1 & 0 \\ 0 & \varpi_w
\end{pmatrix}\cdot y.$$
Noting that $\ker(\alpha)\subseteq (\mathbb{Z}^{\widetilde{Q}}(\widetilde{U}_0(w))_0$ and using \Cref{chargroupvnotinQ}, we finally obtain Ribet's exact sequence
\begin{equation}\label{ribetexseq}
    0 \rightarrow \mathcal{X}_{w}(J_Q(U))\rightarrow \mathcal{X}_{v}(J_{\overline{Q}}(\overline{U}_0(w)))\rightarrow \mathcal{X}_{v}(J_{\overline{Q}}(\overline{U}))^2 \rightarrow 0;
\end{equation}
 see also \cite[Proposition 2.7.6]{Dei14} and \cite[\S 3.2]{Raj01}.
\subsection{Proof of the main Theorem}\label{proofofmaintheorem}
Let $Q_0=\{t\}$ if $[F:\bbQ]$ is even and $Q_0=\{t,s\}$ if $[F:\bbQ]$ is odd, where $s,t$ are non-trivial places. Let $Q$ finite set of finite places of $F$ of indefinite type, containing $Q_0$.

We specialize to the setting of \Cref{Section2}. In particular, let $\lambda_{\mathbf{f}}\colon \overline{\bbT}{}^Q_{\mathcal{O}}(U^{\mathbf{f}})_{\mathfrak{m}}\to \mathcal{O}$ be the augmentation corresponding to $\mathbf{f}$, and set $I_{\mathbf{f}}=\overline{\bbT}{}^Q(U^{\mathbf{f}})\cap \ker(\lambda_{\mathbf{f}})$. We define the optimal quotient $A_{\mathbf{f},Q}$ by the following exact sequence of abelian varieties
    \begin{align*}
        0 \rightarrow I_{\mathbf{f}} J_Q(U^{\mathbf{f}})\rightarrow J_Q(U^{\mathbf{f}}) \xrightarrow{\varphi} A_{\mathbf{f},Q} \rightarrow 0.
\end{align*}
 We write $\delta_Q$ for the $\lambda_{\mathbf{f}}$-Shimura degree of $J_Q(U^{\mathbf{f}})$.
Finally, let  $I_Q=\ker(\overline{\bbT}{}^{Q_0}(U^{\mathbf{f}})\twoheadrightarrow \overline{\bbT}{}^{Q}(U^{\mathbf{f}}))$, and define the abelian variety
    \begin{align*}
         J^{\min}_{Q}\colonequals \big(J_{Q_0}(U^{\mathbf{f}})[I_{Q}]\big)^\circ.
    \end{align*}
 By \cite[Corollary 2.5]{Em03}, we have an isomorphism $T^\vee_{\mathfrak{m}}(J^{\min}_{Q})\cong T^\vee_{\mathfrak{m}}(J_{Q_0}(U^{\mathbf{f}}))/I_{Q}\cong (\overline{\bbT}{}^{Q}_{\mathcal{O}}(U^{\mathbf{f}})_{\mathfrak{m}})^2.$
\begin{proposition}\label{preliminaryfinalresult} Assume that $\mathcal{X}_t(J^{\min}_Q)_{\mathfrak{m}}$ is free of rank one over $\overline{\bbT}{}^{Q}_{\mathcal{O}}(U^{\mathbf{f}})_{\mathfrak{m}}$, and that we have a decomposition $Q=\overline{Q}\sqcup \{v,w\}$ with $Q_0\subset \overline{Q}$. If  $\delta_{\overline{Q}}=\ell_{\mathcal{O}}(\Psi_{\lambda_{\mathbf{f}}}(\overline{\bbT}{}^{\overline{Q}}_{\mathcal{O}}(U^{\mathbf{f}})_{\mathfrak{m}}))$, then the same holds for $Q$, that is
    $$ \delta_{Q}= \ell_{\mathcal{O}}\big(\Psi_{\lambda_{\mathbf{f}}}(\overline{\bbT}{}^{Q}_{\mathcal{O}}(U^{\mathbf{f}})_{\mathfrak{m}})\big),$$
and moreover, $$\delta_Q= \ell_{\mathcal{O}}\big(\Psi_{\lambda_{\mathbf{f}}}(\mathcal{X}_{v}(J_{Q}(U^{\mathbf{f}}))_{\mathfrak{m}})\big)=\ell_{\mathcal{O}}\big(\Psi_{\lambda_{\mathbf{f}}}(\mathcal{X}_{w}(J_{Q}(U^{\mathbf{f}}))_{\mathfrak{m}})\big).$$
\end{proposition}
\begin{proof}
    By Ribet's exact sequence \eqref{ribetexseq}, we have an injection $\mathcal{X}_{w}(J_{Q}(U^{\mathbf{f}}))_{\mathfrak{m}} \hookrightarrow \mathcal{X}_{v}(J_{\overline{Q}}(U^{\mathbf{f}}))_{\mathfrak{m}}$ with torsion free cokernel. Passing to the $\lambda_{\mathbf{f}}$-equivariant parts yields an isomorphism $$\mathcal{X}_{w}(J_{Q}(U^{\mathbf{f}}))_{\mathfrak{m}}[\mathfrak{p}_{\lambda_{\mathbf{f}}}] \xrightarrow{\sim} \mathcal{X}_{v}(J_{\overline{Q}}(U^{\mathbf{f}}))_{\mathfrak{m}}[\mathfrak{p}_{\lambda_{\mathbf{f}}}].$$
    Writing $\mathcal{X}_{w}(J_{Q}(U^{\mathbf{f}}))_{\mathfrak{m}}[\mathfrak{p}_{\lambda_{\mathbf{f}}}]=\mathcal{O}\cdot z$ and applying \Cref{scalartoShimuradeg} twice, we get
    $$ \ord_{\varpi}(u_{J_{\overline{Q}}(U^{\mathbf{f}})}(z,\varphi^*z))=\delta_{\overline{Q}}+c_{v}-2j'_{v}-2u'_{v}
    = \ord_{\varpi}(u_{J_{Q}(U^{\mathbf{f}})}(z,\varphi^*z))=\delta_Q+c_{w}-2j_{w}-2u_{w}.$$
    Since $\mathfrak{m}$ is non-Eisenstein and $v\not\in \overline{Q}$, \cite[Proposition 5]{Raj01} gives $\Phi_{v}(J_{\overline{Q}}(U^{\mathbf{f}}))_{\mathfrak{m}}=0$. By \Cref{someequalities}, we may write $c_{v}=j'_{v}+u'_{v}$, and so \Cref{shdegtocongmodule} yields
    $$ \ell_{\mathcal{O}}\big(\Psi_{\lambda_{\mathbf{f}}}(\overline{\bbT}{}^{\overline{Q}}_{\mathcal{O}}(U^{\mathbf{f}})_{\mathfrak{m}})\big) -\ell_{\mathcal{O}}\big(\Psi_{\lambda_{\mathbf{f}}}(\mathcal{X}_t(J_{Q}(U^{\mathbf{f}}))_{\mathfrak{m}})\big) =\delta_{\overline{Q}}-\delta_Q=c_{v}+c_{w}-2j_{w}-2u_{w}. $$
    Now \Cref{changeofcongmoduleforheckealg} gives
    $$ \ell_{\mathcal{O}}\big(\Psi_{\lambda_{\mathbf{f}}}(\overline{\bbT}{}^{\overline{Q}}_{\mathcal{O},\psi}(U^{\mathbf{f}})_{\mathfrak{m}})\big)- \ell_{\mathcal{O}}\big(\Psi_{\lambda_{\mathbf{f}}}(\overline{\bbT}{}^{Q}_{\mathcal{O},\psi}(U^{\mathbf{f}})_{\mathfrak{m}})\big)= c_{v}+c_{w}. $$
    As noted in \Cref{upperbound}, the inequality $\ell_{\mathcal{O}}(\Psi_{\lambda_{\mathbf{f}}}(\overline{\bbT}{}^{Q}_{\mathcal{O},\psi}(U^{\mathbf{f}})_{\mathfrak{m}}))\ge \ell_{\mathcal{O}}(\Psi_{\lambda_{\mathbf{f}}}(\mathcal{X}_t(J_{Q}(U^{\mathbf{f}}))_{\mathfrak{m}})) $ always holds, and so we must have an equality. In particular, we get that $j_{w}+u_{w}=0$. We can then conclude using the hypothesis of the proposition (note that $v$ and $w$ are interchangeable).
\end{proof}

\begin{proposition}\label{freenessofXt}
    If $Q$ does not contain any trivial place and $v\in Q$, then the module $\mathcal{X}_v(J_Q(U^{\mathbf{f}}))_{\mathfrak{m}}$ is free of rank one over $\overline{\bbT}{}^{Q}_{\mathcal{O}}(U^{\mathbf{f}})_{\mathfrak{m}}$.
\end{proposition}
\begin{proof}
    The proof is identical to \cite[Lemma 6.5]{Helm}. Since the action of $\overline{\bbT}{}^{Q}_{\mathcal{O}}(U^{\mathbf{f}})_{\mathfrak{m}}$ on $\mathcal{X}_v(J_Q(U^{\mathbf{f}}))_{\mathfrak{m}}$ is faithful, it suffices to show that the quotient $X\colonequals \mathcal{X}_v(J_Q(U^{\mathbf{f}}))_{\mathfrak{m}}/\mathfrak{m}\mathcal{X}_v(J_Q(U^{\mathbf{f}}))_{\mathfrak{m}}$ is at most one dimensional over $k$. Recall from \Cref{resultofHelm} that there is a $\overline{\bbT}{}^{Q}_{\mathcal{O}}(U^{\mathbf{f}})_{\mathfrak{m}}[G_F]$-equivariant surjection 
    $$ T^\vee_{\mathfrak{m}}(J_Q(U^{\mathbf{f}}))\twoheadrightarrow \mathcal{X}_v(J_Q(U^{\mathbf{f}}))_{\mathfrak{m}}, $$
     and by \Cref{freenessoverGorenstein} we have $T^\vee_{\mathfrak{m}}(J_Q(U^{\mathbf{f}}))\cong (\overline{\bbT}{}^{Q}_{\mathcal{O}}(U^{\mathbf{f}})_{\mathfrak{m}})^2$, which reduces to $T^\vee_{\mathfrak{m}}(J_Q(U^{\mathbf{f}}))/\mathfrak{m}\cong \overline{\rho}^*(1)$ as $G_F$-modules. Reducing the surjection modulo $\mathfrak{m}$, we get a surjection $\overline{\rho}^*(1)\twoheadrightarrow X$. If $\dim_k(X)>1$, then $X\cong \overline{\rho}^*(1)$ as $G_F$-modules. As $I_v$ acts trivially on $\mathcal{X}_v(J_Q(U^{\mathbf{f}}))_{\mathfrak{m}}$, we get that $\overline{\rho}$ is unramified at $v$. We know by \cite[Théorème (A), \S 5-6]{Car} (see also \cite[Proposition 7]{Raj01}) that $\Frob_{v}$ acts as $\mathbf{U}_{v}$ on $\mathcal{X}_v(J_{Q}(U^{\mathbf{f}}))_{\mathfrak{m}}$. Here, $\mathbf{U}_{v}$ is regarded as an element of $\overline{\bbT}{}^{Q}_{\mathcal{O}}(U^{\mathbf{f}})_{\mathfrak{m}}$ through the surjection $\overline{\bbT}{}^{Q'}_{\mathcal{O}}(U^{\mathbf{f}})_{\mathfrak{m}}\twoheadrightarrow \overline{\bbT}{}^{Q}_{\mathcal{O}}(U^{\mathbf{f}})_{\mathfrak{m}}$ coming from the Jacquet--Langlands correspondence with  $Q'=Q\setminus\{v\}$. Thus $\Frob_v$ acts on $X$ as a scalar, which is a contradiction since $v$ is not a trivial place for $\overline{\rho}$.
\end{proof}

\begin{corollary}\label{freeToricModule}
    The module $\mathcal{X}_t(J_{Q}^{\min})_{\mathfrak{m}}$ is free of rank one over $\overline{\bbT}{}^{Q}_{\mathcal{O}}(U^{\mathbf{f}})_{\mathfrak{m}}$.
\end{corollary}
\begin{proof}
    By \Cref{freenessofXt}, the module $\mathcal{X}_{t}(J_{Q_0}(U^{\mathbf{f}}))_{\mathfrak{m}}$ is free of rank one over $\overline{\bbT}{}^{Q_0}_{\mathcal{O}}(U^{\mathbf{f}})_{\mathfrak{m}}$. We recall that $I_Q=\ker(\overline{\bbT}{}^{Q_0}(U^{\mathbf{f}})\twoheadrightarrow \overline{\bbT}{}^{Q}(U^{\mathbf{f}}))$, and 
$
         J^{\min}_{Q}\colonequals \big(J_{Q_0}(U^{\mathbf{f}})[I_{Q}]\big)^\circ.
$
   The injection $J^{\min}_{Q}\hookrightarrow J_{Q_0}(U^{\mathbf{f}})$ induces a $\overline{\bbT}{}^{Q_0}_{\mathcal{O}}(U^{\mathbf{f}})_{\mathfrak{m}}$-equivariant surjective map
    $$ \mathcal{X}_{t}(J_{Q_0}(U^{\mathbf{f}}))_{\mathfrak{m}}\twoheadrightarrow \mathcal{X}_{t}(J^{\min}_{Q})_{\mathfrak{m}}.$$
    Since $ \mathcal{X}_{t}(J_{Q_0}(U^{\mathbf{f}}))_{\mathfrak{m}}$ is free of rank one over $\overline{\bbT}{}^{Q_0}_{\mathcal{O}}(U^{\mathbf{f}})_{\mathfrak{m}}$, the same is true of $\mathcal{X}_{t}(J^{\min}_{Q})_{\mathfrak{m}}$ over $\overline{\bbT}{}^{Q}_{\mathcal{O}}(U^{\mathbf{f}})_{\mathfrak{m}}$.
\end{proof}
\begin{corollary}\label{finalresultoverindefinite}
We have the following equalities  
    $$ \delta_Q= \ell_{\mathcal{O}}\big(\Psi_{\lambda_{\mathbf{f}}}(M_{Q}(U^{\mathbf{f}}))\big) =  \ell_{\mathcal{O}}\big(\Psi_{\lambda_{\mathbf{f}}}(\overline{\bbT}{}^{Q}_{\mathcal{O}}(U^{\mathbf{f}})_{\mathfrak{m}})\big). $$
    Moreover, let $v\not\in Q$ is a finite place at which $\mathbf{f}$ is Steinberg, and set $\widetilde{Q}=Q\cup\{v\}$. Then we have
$$ \ell_{\mathcal{O}}\big(\Psi_{\lambda_{\mathbf{f}}}(M_{\widetilde{Q}}(U^{\mathbf{f}}))\big) =  \ell_{\mathcal{O}}\big(\Psi_{\lambda_{\mathbf{f}}}(\overline{\bbT}{}^{\widetilde{Q}}_{\mathcal{O}}(U^{\mathbf{f}})_{\mathfrak{m}})\big).$$
\end{corollary}

\begin{proof}
    The first part follows by induction using \Cref{freeToricModule}, \Cref{tatemoduletochargrp}, and \Cref{preliminaryfinalresult}.
    
     For the second part, since $\mathfrak{m}$ is non-Eisenstein and $v\not\in Q$, \cite[Proposition 5]{Raj01} gives $\Phi_v(J_{Q}(U^{\mathbf{f}}))_{\mathfrak{m}}=0$. By \Cref{someequalities}, we may therefore write $c_v=j_v+u_v$, and  \Cref{shdegtocongmodule} yields $$\delta_Q=c_v+ \ell_{\mathcal{O}}\big(\Psi_{\lambda_{\mathbf{f}}}(\mathcal{X}_v(J_Q(U^{\mathbf{f}}))_{\mathfrak{m}})\big).$$ Since $\mathcal{X}_v(J_Q(U^{\mathbf{f}}))_{\mathfrak{m}}\cong M_{\widetilde{Q}}(U^{\mathbf{f}})$ by \Cref{chargroupvnotinQ}, the claim follows from the first part together with \Cref{changeofcongmoduleforheckealg}.
\end{proof}

\appendix
\section{Gröbner bases Calculations}\label{appendixA}

We work within the setting of the proof of \Cref{cotangentuni}. The purpose of this appendix is to construct bases for the rings $\widetilde{R}_{\mathbb{Z}}$ and $R^{\uni}_{\mathbb{Z}}$ from \emph{loc. cit.}, using the theory of Gröbner bases.

Let $K$ be a field, and consider the polynomial ring $K[y,a,b,c]=K[x,y,t,z,a,b,c]/(\widetilde{\theta}_1,\widetilde{\theta}_2,\widetilde{\theta}_3)$, where $\widetilde{\theta}_1=t-z$, $\widetilde{\theta}_2=b-z$, $\widetilde{\theta}_3=c-x$. Throughout, monomials are ordered via the homogeneous lexicographic order satisfying $1\prec c \prec b \prec a \prec y$. Given $f\in K[y,a,b,c]$, we define its initial term $\mathrm{in}(f)$ to be its greatest term with respect to $\prec$, and given an ideal $J\subset K[y,a,b,c]$, we let $\mathrm{in}(J)$ denote the ideal generated by the initial terms of the elements of $J$ (see \cite[\S 15]{Eisen}). Given $f\in K[y,a,b,c]$, and $G=\{g_i,\dots,g_t\}\subset K[y,a,b,c] $, we say that $f$ reduces to $0$ modulo $G$ if we can write
$$ f=\sum_i f_i g_i, \quad \text{with }\mathrm{in}(f_ig_i)\prec \mathrm{in}(f)\text{ for all }i. $$
 Given $g_i,g_j\in G$, we define their $S$-polynomial to be
$$ S(g_i,g_j)\colonequals \frac{\mathrm{in}(g_j)g_i-\mathrm{in}(g_i)g_j}{\gcd(\mathrm{in}(g_i),\mathrm{in}(g_j))}.  $$

\subsection{} We consider the ring $R_K=K[y,a,b,c]/I$, where $I$ is the ideal generated by 
\begin{alignat*}{2}
    \widetilde{s}_1 &= c^2, &\quad\quad\quad
    \widetilde{s}_2 &= ac + b^2 - bc, 
     \\\widetilde{s}_3 &= a^2 + 2bc, &\quad\quad\quad
    \widetilde{s}_4 &= bc - y^2 - y^2c - b^2 - b^2c + 2cy + c^2.
\end{alignat*}
We let $G=\{g_1,\dots,g_7\}\subset I$, where 
\begin{align*}
g_1 &= y^2-2cy+b^2-bc = (y^2 - 2y + b^2 - b - c + 1)\,\widetilde{s}_1 + (c-1)\,\widetilde{s}_4, \quad g_2 = \widetilde{s}_3, \\
g_3 &= ab^2+b^3 = (-a+3b)\,\widetilde{s}_1 + (a+b+c)\,\widetilde{s}_2 - c\,\widetilde{s}_3, \quad g_4 = \widetilde{s}_2, \\
g_5 &= b^4 = (a-b)^2\,\widetilde{s}_1 + (-ac+b^2+bc)\,\widetilde{s}_2, \quad g_6 = b^2c = (b-a)\,\widetilde{s}_1 + c\,\widetilde{s}_2, \quad g_7 = \widetilde{s}_1.
\end{align*}
To verify that $G$ is a Gröbner basis, i.e. that the initial terms $\mathrm{in}(g_i)$ generate $\mathrm{in
}(I)$, 
%and $\mathrm{in}(g_i)\nmid g_j$ for $i\neq j$, 
it suffices by Buchberger's criterion \cite[Theorem 15.8]{Eisen} to check that every $S$-polynomial $S(g_i,g_j)$ reduced to $0$ modulo $G$ . 
We may skip checking the pairs $g_i$ and $g_j$ whose initial terms are coprime by \cite[Exercise 15.20]{Eisen}, as well as any pair of monomials since in that case $S(g_i,g_j)=0$. Therefore, it remains to check the following $S$-polynomials.
{\begin{table}[h]
\centering
\footnotesize
{\renewcommand{\arraystretch}{0.96}
\begin{tabular}{|c|l|}
\hline
Pair & $S(g_i,g_j)  $ \\
\hline
$(g_2,g_3)$ & $-ab^3 + 2b^3c = -b g_3 + g_5 + 2b g_6$ \\
$(g_2,g_4)$ & $-ab^2 + abc + 2bc^2 = -g_3 + b g_4 + g_6 + 2b g_7$ \\
$(g_3,g_4)$ & $-b^4 + 2b^3c = -g_5 + 2b g_6$ \\
$(g_3,g_5)$ & $b^5 = b g_5$ \\
$(g_3,g_6)$ & $b^3c = b g_6$ \\
$(g_4,g_6)$ & $b^4 - b^3c = g_5 - b g_6$ \\
$(g_4,g_7)$ & $b^2c - bc^2 = g_6 - b g_7$ \\
\hline
\end{tabular}}
\end{table}}

Consequently, by Macaulay's theorem \cite[Theorem 15.3]{Eisen}, a basis of $R_K$ is given by the monomials that are not divisibly by any $\mathrm{in}(g_i)$, which are
$\{1,y,a,b,c,yc,ya,yab,yb^2,ybc,yb,ab,b^3,b^2,bc,yb^3\}.$

Using the relations $\widetilde{s}_2$ and $g_6$, we find that $ac\equiv bc-b^2\bmod I$, and so $abc\equiv -b^3 \bmod I$. It follows that $yac\equiv ybc-yb^2\bmod I$ and $yabc\equiv-b^3 \bmod I$. Hence the elements
    \begin{gather*}
               b_{1}= 1,\ b_{2}= y,\ b_{3} = a, \ b_{4}=b,\ b_{5}= c,\ b_{6}= yc,  \ b_7 = ya,\ b_8 = yab,
       \\ b_9= yac,\ b_{10}= ybc,\ b_{11}= yb,\ b_{12}= ab,\ b_{13}= abc,\ b_{14}= ac,\ b_{15}=  bc, \ b_{16}= yabc.
    \end{gather*}
form a basis of $R_K$. 

\subsection{} We now consider the ring $R^{\uni}_K=K[y,a,b,c]/J$, where $J$ is the ideal generated by 
\begin{gather*}\widetilde{r}_1=ac, \ \widetilde{r}_2=bc, \ \widetilde{r}_3=c^2, \ 
\widetilde{r}_4=y+(y^2+b^2)(1+c)-y(1+c)^2,\\   \widetilde{r}_5=a^2+bc,\  \widetilde{r}_6=ab-cy, \  \widetilde{r}_7=ay+cb, \ \widetilde{r}_8=ay+b^2, \ \widetilde{r}_9=b(y-a).
\end{gather*}
We let $G=\{g_1,\cdots,g_9\}$, where
\begin{align*}
g_1 &= y^2 - 2ab = \widetilde{r}_4 - \widetilde{r}_8 + \widetilde{r}_7 - \widetilde{r}_2 - a\widetilde{r}_9 - b\widetilde{r}_5 + b\widetilde{r}_2 + y\widetilde{r}_6 - c\widetilde{r}_8 + c\widetilde{r}_7 - c\widetilde{r}_2 - 2\widetilde{r}_6 + y\widetilde{r}_3, \\
g_2 &= ay = \widetilde{r}_7 - \widetilde{r}_2, \quad g_3 = by - ab = \widetilde{r}_9, \quad g_4 = yc - ab = -\widetilde{r}_6, \quad
g_5= a^2 = \widetilde{r}_5 - \widetilde{r}_2, \\ g_6 &= ac = \widetilde{r}_1, \quad g_7 = b^2 = \widetilde{r}_8 - \widetilde{r}_7 + \widetilde{r}_2, \quad g_8 = bc = \widetilde{r}_2, \quad  g_9 = c^2 = \widetilde{r}_3.
\end{align*}
As before, we verify that the $S$-polynomials reduce to 0 modulo $G$. The only non-trivial pairs are the following
\begin{table}[h]
\centering
\footnotesize
{\renewcommand{\arraystretch}{0.96}
\begin{tabular}{|c|l|c|l|}
\hline
Pair & $S(g_i,g_j)$ & Pair & $S(g_i,g_j)$ \\
\hline
$(g_1,g_2)$ & $-2a^2b = -2b\cdot g_5$ & $(g_2,g_3)$ & $a^2b = b\cdot g_5$ \\
$(g_1,g_3)$ & $aby - 2ab^2 = a\cdot g_3 + b\cdot g_5 - 2a\cdot g_7$ & $(g_2,g_4)$ & $a^2b = b\cdot g_5$ \\
$(g_1,g_4)$ & $aby - 2abc = a\cdot g_3 + b\cdot g_5 - 2a\cdot g_8$ & $(g_3,g_4)$ & $ab^2 - abc = a\cdot g_7 - a\cdot g_8$ \\
$(g_4,g_6)$ & $-a^2b = -b\cdot g_5$ & $(g_3,g_7)$ & $-ab^2 = -a\cdot g_7$ \\
$(g_4,g_8)$ & $-ab^2 = -a\cdot g_7$ & $(g_3,g_8)$ & $-abc = -a\cdot g_8$ \\
$(g_4,g_9)$ & $-abc = -a\cdot g_8$ & & \\
\hline
\end{tabular}}
\end{table}

Hence we find that $\mathrm{in}(J)=\langle y^2,ay,by,cy,a^2,ac,b^2,bc,c^2\rangle$, and Macaulay's theorem gives a basis $\{1,y,a,b,c,ab\}$ of $R^{\uni}_K$. Since $yc=ab+g_4$, we deduce that $b_1,\dots,b_6$ form a basis of $R^{\uni}_K$.

\bibliographystyle{alpha}
%\bibliography{literature}

\end{document}

%% file: operators.tex
\newcommand{\bbA}{\mathbb A}

\newcommand{\bbC}{\mathbb C}

\newcommand{\bbF}{\mathbb F}

\newcommand{\bbH}{\mathbb H}

\newcommand{\bbN}{\mathbb N}

\newcommand{\bbQ}{\mathbb Q}
\newcommand{\bbR}{\mathbb R}

\newcommand{\bbT}{\mathbb T}

\newcommand{\bbZ}{\mathbb Z}

\renewcommand{\sl}{\mathfrak{sl}}

\newcommand{\GL}{\operatorname{GL}} % general linear group
\newcommand{\SL}{\operatorname{SL}} % special linear group
\newcommand{\CNL}{\operatorname{CNL}}

\newcommand{\Ext}{\operatorname{Ext}} % Ext functor
\newcommand{\Spec}{\operatorname{Spec}} % Spectrum of a ring
\newcommand{\Spf}{\operatorname{\mathsf{Spf}}} % formal spectrum

\newcommand{\Ann}{\operatorname{Ann}} % Annihilator
\newcommand{\Aut}{\operatorname{Aut}} % Automorphism group
\newcommand{\coker}{\operatorname{coker}} % cokernel
\newcommand{\Der}{\operatorname{Der}} % Derivations
\newcommand{\Fitt}{\operatorname{Fitt}} % Fitting ideal

\newcommand{\Art}{\operatorname{Art}} % artin map

\newcommand{\dd}{\operatorname{d}} % Derive

\newcommand{\End}{\operatorname{End}} % endomorphism set
\newcommand{\Frob}{\operatorname{Frob}} % Frobenius endomorphism
\newcommand{\Gal}{\operatorname{Gal}} % Galois group
\newcommand{\Hom}{\operatorname{Hom}} % homomorphism set
\newcommand{\id}{\operatorname{id}} % identity
\newcommand{\im}{\operatorname{im}} % image
\newcommand{\Ind}{\operatorname{Ind}} % Induced representation
\newcommand{\ord}{\operatorname{ord}} % Order
\newcommand{\PGL}{\operatorname{PGL}} % projective general linear group
\newcommand{\Pic}{\operatorname{Pic}} % Picard group of a ringed space
\newcommand{\rank}{\operatorname{rank}} % rank
\newcommand{\Res}{\operatorname{\mathrm{Res}}} % Restricted representation
\newcommand{\tr}{\operatorname{tr}} % trace
\newcommand{\Sets}{\operatorname{Sets}} % Sets
\newcommand{\tf}{\operatorname{tf}} %torsionfree
\newcommand{\tors}{\operatorname{tors}} %torsionpart
\newcommand{\Nrd}{\operatorname{Nrd}}
\newcommand{\uni}{\operatorname{uni}}
\newcommand{\unr}{\operatorname{unr}}

\newcommand{\st}{\operatorname{st}}
\newcommand{\proj}{\operatorname{proj}}
\newcommand{\alt}{\operatorname{alt}}